\documentclass[paper=a4,bibtotoc,parskip=half]{scrartcl}

\usepackage{mathrsfs}  
\usepackage{amssymb}
\usepackage{amsmath}
\usepackage{amsthm}
\usepackage{graphicx}
\usepackage{float}
\usepackage{placeins}
\usepackage{xargs}  % Use more than one optional parameter in a new commands
\usepackage[pdftex,dvipsnames]{xcolor} 
\usepackage{array}
\usepackage{calc}
\usepackage{hf-tikz}
\usepackage{multicol}
\usepackage[ruled,noline]{algorithm2e}
\SetKwInOut{Input}{Input} 
\SetKwInOut{Output}{Output}
\usepackage[font=small,labelfont=bf]{caption}
\usepackage{subcaption}
\usepackage{rotating}
\DeclareGraphicsExtensions{.png,.eps}
\usepackage[outdir=./]{epstopdf}

\definecolor{orange}{rgb}{1,0.4,0}
\definecolor{darkjunglegreen}{rgb}{0.1, 0.14, 0.13}
\definecolor{asparagus}{rgb}{0.53, 0.66, 0.42}
\definecolor{burntumber}{rgb}{0.54, 0.2, 0.14}
\definecolor{arsenic}{rgb}{0.23, 0.27, 0.29}
\definecolor{brightmaroon}{rgb}{0.76, 0.13, 0.28}
\definecolor{oxfordblue}{rgb}{0.0, 0.13, 0.28}
\definecolor{outerspace}{rgb}{0.25, 0.29, 0.3}
\definecolor{bluebell}{rgb}{0.64, 0.64, 0.82}

\usepackage{pgfplots}
\usepackage{tikz}
\usepackage{xypic}
\usepackage{tcolorbox}

\usetikzlibrary{matrix}
\usetikzlibrary{decorations.pathreplacing}

\tikzset{
	mybrace/.style={decorate,decoration={brace,aspect=#1}}
}
\usetikzlibrary{fit}
\tikzset{%
	highlight/.style={rectangle,fill=red!15,draw,fill opacity=0.5,thick,inner sep=0pt}
}
\tikzset{%
	highlightb/.style={rectangle,fill=blue!15,draw,fill opacity=0.5,thick,inner sep=0pt}
}
\tikzset{%
	highlightw/.style={rectangle,fill=none,draw,fill opacity=0.5,thick,inner sep=0pt}
}
\newcommand{\tikzmark}[2]{\tikz[overlay,remember picture,baseline=(#1.base)] \node (#1) {#2};}
\newcommand{\Highlight}[1][submatrix]{%
	\tikz[overlay,remember picture]{
		\node[highlight,fit=(left.north west) (right.south east)] (#1) {};}
}

\newcommand{\Highlightb}[1][submatrix]{%
	\tikz[overlay,remember picture]{
		\node[highlightb,fit=(left.north west) (right.south east)] (#1) {};}
}

\usepackage{colortbl}

\usepackage{xcolor}

\usepackage{booktabs}

\usepackage{yhmath}% http://ctan.org/pkg/yhmath
\usepackage{mathdots}% http://ctan.org/pkg/mathdots

\usepackage{microtype}

\usepackage{epsfig} %% add this line and the next only if you have pictures
\usepackage{graphics} %% pictures should be in esp format

\usepackage{MnSymbol}% http://ctan.org/pkg/mnsymbol

\usepackage{mathtools}

\newcommand{\rectangle}[3][0pt]{{%
		\fboxsep=-\fboxrule\sbox0{}\wd0=#2\ht0=#3\dp0=#1\relax\fbox{\box0}}}

\newcolumntype{?}{!{\vrule width 1.5pt}}

\usepackage[colorlinks=true,
           linkcolor=red,
            urlcolor=blue,
            citecolor=black]{hyperref}

\numberwithin{equation}{section}

\theoremstyle{plain}
\newtheorem{theorem}{Theorem}[section]
\newtheorem{lemma}[theorem]{Lemma}
\newtheorem{prop}[theorem]{Proposition}

\theoremstyle{definition}
\newtheorem{definition}{Definition}[section]

\theoremstyle{remark}
\newtheorem{rem}{Remark}[section]

\newcommand{\mc}[1]{\mathcal{#1}}
\newcommand{\wh}[1]{\widehat{#1}}

\newcommand{\deff}{:=}

\renewcommand{\vec}[1]{\mbox{vec}(#1)}

\newcommand{\lsim}{{\;\raise0.3ex\hbox{$<$\kern-0.75em\raise-1.1ex\hbox{$\sim$}}\;}}

\newcommand{\kron}{\otimes}

\newcommand{\bfP}{\mathbf{P}}
\newcommand{\bfD}{\mathbf{D}}
\newcommand{\bfZ}{\mathbf{Z}}
\newcommand{\bfL}{\mathbf{L}}
\newcommand{\bfM}{\mathbf{M}}
\newcommand{\bfX}{\mathbf{X}}
\newcommand{\bfY}{\mathbf{Y}}
\newcommand{\bfA}{\mathbf{A}}
\newcommand{\bfB}{\mathbf{B}}
\newcommand{\bfG}{\mathbf{G}}
\newcommand{\bfS}{\mathbf{S}}
\newcommand{\bfI}{\mathbf{I}}
\newcommand{\bfE}{\mathbf{E}}

\newcommand{\bfC}{\mathbf{C}}
\newcommand{\bfR}{\mathbf{R}}
\newcommand{\bfU}{\mathbf{U}}
\newcommand{\bfV}{\mathbf{V}}
\newcommand{\bfQ}{\mathbf{Q}}
\newcommand{\bfJ}{\mathbf{J}}

\newcommand{\bbM}{\mathbb{M}}

\newcommand{\bfW}{\mathbf{W}}

\newcommand{\bel}{\mathsf{bel}}

\newcommand{\bfzero}{\mathbf{0}}

\newcommand{\RR}{\mathbb{R}}
\newcommand{\CC}{\mathbb{C}}

\newcommand{\Mat}[2]{\mathcal{M}_{#1 , #2}}

\newcommand{\gaussd}{\mc{N}}

\newcommand{\p}{\partial}
\renewcommand{\Re}{\operatorname{Re}}             % real part
\newcommand{\Laplacian}{\Delta}

\newcommand{\figref}[1]{\figurename~\ref{#1}}

\newcommand{\myointersection}{
	\mathbin{
		\mathchoice
		{\ointersection{\displaystyle}}
		{\ointersection{\textstyle}}
		{\ointersection{\scriptstyle}}
		{\ointersection{\scriptscriptstyle}}
	}
}
\newcommand{\ointersection}[1]{\tikz[baseline=(X.base), inner sep=0, outer sep=0]\node[draw,circle] (X) {$#1\cap$};}

\newcommand*{\arr}{\longrightarrow}

\author{{{Lorenzo Baldassari}} \and {{Andrea Scapin}}}
\title{Multi-scale classification for electro-sensing\thanks{\footnotesize This work was supported by the SNF grant 200021-172483.}}
\begin{document}
\maketitle

\begin{abstract}
	This paper introduces premier and innovative (real-time) multi-scale method for target classification in electro-sensing. The intent is that of mimicking the behavior of the weakly electric fish, which is able to retrieve much more information about the target by approaching it. The method is based on a family of transform-invariant shape descriptors computed from generalized polarization tensors (GPTs) reconstructed at multiple scales. The evidence provided by the different descriptors at each scale is fused using Dempster-Shafer Theory. Numerical simulations show that the recognition algorithm we proposed performs undoubtedly well and yields a robust classification.

	\par\vskip\baselineskip\noindent
\end{abstract}

\section{Introduction}

The biological behavior of weakly electric fish has been studied by scholars for years. These fish orient themselves at night in complete darkness by using electrosensory information, which makes these animals an ideal subject for developing bio-inspired imaging techniques. Such interest has motivated a huge number of studies addressing the active electro-sensing problem from many different perspectives since Lissmann and Machin’s work \cite{p5,p6,p7,p8,p9,p10,p11,p12}. One of the most noteworthy potential bio-inspired applications is in underwater robotics. Building autonomous robots with electro-sensing technology may supply unexplored navigation, imaging
and classification capabilities, especially when the sight is unreliable due, for example, to the turbidity of the surrounding waters or the poor lighting conditions \cite{p14,p15}.

From the mathematical point of view, the electro-sensing problem is to detect and locate the dielectric target and to identify its shape and material parameters given the current
distribution over the skin of the fish. 
Ammari et al. \cite{AMM2013} stated a rigorous mathematical model for treating the inverse problem of electro-sensing dielectric objects. They exploited the smallness of targets in order to apply the framework of asymptotic small-volume expansions. The electric current, which contains information on the target, is measured by a discrete number of receptors along the fish body. When enough measurements are collected, it is possible to recover the contracted generalized polarization tensors (CGPTs), which do encode information about the unknown target. One way in which the fish can acquire enough independent measurements is by exploiting the movement, i.e., by collecting several static measurements while swimming around the target \cite{Ammari_2017}. In this way, it creates a synthetic-aperture view of the dielectric object that yields high-resolved reconstruction of its features. Although the inverse problem is severly ill-posed, classification works well. 
In \cite{Ammari2013} new shape descriptors, relying upon the CGPTs, which are invariants under rotations, translations, and scaling of the target, are found. 
In the previous works, a single circular trajectory around the target has been considered.  In the two dimensional case, for small targets, the magnitude of the electric signal due to the presence of the target is of order $\varepsilon^2$, where $\varepsilon$ is the length-scale, which is of the same order throughout the whole trajectory. In this type of setting it is natural to reconstruct the target's features only up to some small order $K^*$, which is called the resolving order. $K^*$ is essentially determined by the signal-to-noise ratio (SNR), that sets a limit to the fineness of the reconstruction we are capable of, see, for instance, \cite{Ammari2014}.  
It has been shown that the reconstruction is accurate enough to perform a dictionary matching approach for both homogeneous and inhomogeneous objects, see \cite{Ammari11652,SCAPIN20191872}.

%In the recent years, the interest in the applications of optimal experimental design theory has considerably increased. 

%
%The aim of this work is to improve the recognition capabilities of the fish by acquiring measurements at different length-scales on multiple circular orbits around the target. Since the fish is able to retrieve much more information about its shape and material parameters when approaching it, the classification is expected to be more robust as soon as multiple circular orbits are considered.
%
%When dealing with expert systems \cite{WALLEY19961} an important problem people face is that of classification. When many classifiers are available, the problem of combining them to enhance the classification capabilities naturally arises, see e.g. \cite{155943}.
%The approach we present in this paper relies on the so-called Transferable Belief Model (TBM), see, for instance, \cite{1591952}.
%The output of a scoring-classifier, i.e., a list of numerical scores (a score for each element of the dictionary) that corresponds to the evidence at hand, can be converted into a belief assignment. Following \cite{Huynh2010AdaptivelyEW,MONDEJARGUERRA201557}, a natural way to translate the scores into beliefs is to consider the Shannon's entropy as a confidence factor associated to the evidence.
%Belief assignments are then combined by means of some combination rule, such as Dempster-Shafer rule, in order to obtain a synthetized new belief that pulls together all the information.
%This approach is particularly suitable for the electro-sensing problem hereabove.

The aim of this work is to improve the recognition capabilities of the fish by acquiring measurements at different length-scales on multiple circular orbits around the target. The main advantage of the multi-scale configuration is that the descriptors introduced in \cite{Ammari2014} can be compared at different orders up to the resolving order, which is increasing with respect to the length-scale. Therefore, selecting different comparison orders produces different classifiers.	
When many classifiers are available, the problem of combining them to enhance the classification capabilities naturally arises, see e.g. \cite{155943}.
The approach we present in this paper relies on the so-called Transferable Belief Model (TBM), see, for instance, \cite{1591952}.
The output of a scoring-classifier, i.e., a list of numerical scores (a score for each element of the dictionary) that corresponds to the evidence at hand, can be converted into a belief assignment. Following \cite{Huynh2010AdaptivelyEW,MONDEJARGUERRA201557}, a natural way to translate the scores into beliefs is to consider the Shannon's entropy as a confidence factor associated to the evidence.
Belief assignments are then combined by means of some combination rule, such as Dempster-Shafer rule, in order to obtain a synthetized new belief that pulls together all the information.
This approach is particularly suitable for the electro-sensing problem hereabove.
Since the fish is able to retrieve much more information about its shape and material parameters when approaching it, the classification is expected to be more robust as soon as multiple circular orbits are considered.

The paper is structured as follows. In Section \ref{sec:model_electrosensing}, a preliminary description of the experimental design for electro-sensing is discussed. We show that the design matrix associated to the forward linear operator $\mc{L}$ defined in \cite{Ammari11652} can be expressed as a generalized block Kronecker product by vectorization, see \cite{regalia1989kronecker,Marco2018}. A reflexive minimum norm g-inverse of the acquisition operator, arising in a natural way from the block Kronecker structure, is also used. This has been recently introduced in \cite{Marco2018} in the context of bivariate polynomial regression.
%Although this form is simple for computations, it is desirable to have the CGPTs rearranged order-wisely and block-wisely. Therefore other forms involving the Tracy-Singh product of matrices are presented. 
%Furthermore, thanks to the particular form of the design matrix, a method to recursively update the reconstruction position after position right away. This online update is based on Greville's well known formulas in \cite{doi:10.1137/1002004}.
Based on Greville's well known formulas in \cite{doi:10.1137/1002004}, this g-inverse provides a method to update recursively the estimates position after position right away.

In Section \ref{sec:rank-analysis}, a detailed analysis of the structure of the design matrix is carried out. In particular, the need of creating a synthetic-aperture view is readily understood by inspecting the rank of the design matrix. Assuming a circular acquisition setting, i.e., the fish collects the data swimming on a circular trajectory around the target, an estimate on the reconstruction error of the CGPTs is derived. The estimate has an upper bound depending on the length-scale, and it is formally equal to that given in \cite{Ammari2014}. Finally, issues related with limited-view data, i.e., data collected by receptors covering a limited angle of view, are discussed. In particular, a study of the spectrum of the matrix of receptors shows the impact of the angle of view on the reconstruction: the closer the angle of view is to $2\pi$, the more informative the estimate becomes. Furthermore, if the reconstruction order is small, the limited-view configuration has a minor impact on the reconstructed CGPTs.

In Section \ref{sec:recognition}, the classification problem based on a multi-scale acquisition setting is addressed. In particular, measurements at different length-scales are used to improve the resolving power in the reconstruction of the CGPTs. As a matter of fact, the closer the orbit is to the target, the higher is the SNR, the higher is the order of CGPTs-based descriptors that can be used in the comparison. 
A matching algorithm, which generalizes the one proposed in \cite{Ammari11652} to the case study we consider, is presented. Firstly, a certain number of concentric orbits around the target are thoroughly chosen. On each orbit, a comparison between the theoretical and measured shape descriptors up to a properly chosen length-scale dependent order is required. Similarly to \cite{Ammari11652}, the comparison is done by means of a given metric, and it yields a list of scores. The normalized list of scores produced on each orbit is converted into an evidence distribution, which is then stored. The Shannon's entropy is used as a confidence factor, see \cite{Huynh2010AdaptivelyEW,MONDEJARGUERRA201557}.
The evidence distributions computed along different orbits are subsequently combined by using the TBM conjunctive rule introduced in \cite{Smets2008}.

In Section \ref{sec:numerical-exp} we perform numerical simulations in order to test the performance of the recognition algorithm, introduced in Section \ref{sec:recognition}, on a particular dictionary of dielectric targets. The reported results show an enhancement of the recognition rate, corroborating the idea that combining descriptors at different length-scales makes the classification more robust. Both the minimum norm reflexive generalized inverse and the Moore-Penrose inverse are used in the reconstruction.

%We also compare the solution given by the minimum-norm reflexive generalized inverse with the one given by recursive least-square (RLS) estimation (see Appendix for RLS). 

\section{Model specification for electro-sensing}

\label{sec:model_electrosensing}

Let us now briefly summarise the model of electro-sensing derived in \cite{AMM2013}: the body of the fish is $\Omega$, an open bounded set in $\RR^2$, with smooth boundary $\partial \Omega$, and with outward normal unit vector denoted by $\nu$. 
The electric organ is a dipole $f(x)$ inside $\Omega$ or a sum of point sources inside $\Omega$ satisfying the charge neutrality condition. The skin of the fish is very thin and highly resistive. Its effective thickness, that is, the skin thickness times the contrast between the water and the skin conductivities, is denoted by $\xi$, and it is much smaller than the fish size. We assume that the conductivity of the background medium is one. We consider a smooth bounded target $D = \delta B$, and $B$ is a smooth bounded domain containing the origin. We assume that the conductivity of $D$ is $0 < k \ne 1$, and define the contrast $\lambda \deff (k+1)/(2(k-1))$. In the presence of $D$, the electric potential emitted by the fish is the solution to the following equations: 
\begin{equation} \label{eq:model_u} \begin{cases} \Delta u = f  & \mbox{in } \Omega ,\\ \nabla \cdot (1 + (k-1) \chi_D ) \nabla u = 0  &  \mbox{in } \RR^2 \setminus \overline{\Omega} ,\\ u|_+ - u |_- = \xi \dfrac{\partial u}{\partial \nu} \biggr |_+  &  \mbox{on } \partial \Omega,\\  \dfrac{\partial u}{\partial \nu} \biggr |_- = 0   & \mbox{on } \partial \Omega ,\\ |u(x)| = O(|x|^{-1})  & \mbox{as } |x| \to \infty.\end{cases} \end{equation}
Here, $\chi_D$ is the characteristic function of $D$, $\partial/\partial \nu$ is the normal derivative, and $|_{\pm}$ denotes the limits from, respectively, outside and inside $\Omega$. Following \cite{Ammari11652}, we introduce the function $H$ defined as
\begin{equation} \label{eq:H_fun} H(x) \deff p(x) + \mc{S}_{\Omega}\left [ \frac{\p u}{\p \nu} \biggr |_+ \right ] - \xi \mc{D}_{\Omega}\left [ \frac{\p u}{\p \nu} \biggr |_+ \right ], \end{equation} 
where $\Laplacian p = f$ on $\RR^2$. $\mc{S}_{\Omega}$ and $\mc{D}_{\Omega}$ are the single- and double-layer potentials, respectively, defined in Appendix \ref{apx:GPTs}. It is readily seen that the following representation formula holds:
\begin{equation} \label{eq:u-H-rep} u(x) - H(x) = \mc{S}_{D} (\lambda {I} - \mc{K}_D^*)^{-1} \left ( \frac{\p H}{\p \nu} \right ), \end{equation} 
where ${I}$ is the identity and $\mc{K}_D^*$ is the Neumann-Poincar\'e operator associated to the target $D$, see Appendix \ref{apx:GPTs}.

\subsection{Data acquisition system} \label{sbsec:setting}

In this section we aim at describing the data acquisition system, i.e., the experimental setting we shall adopt to solve the inverse problem.

\medskip

As we briefly mentioned in the introduction, the fish use the movement in order to swim around the target, creating a synthetic aperture view.

Suppose that the scanning movement consists of a single circular orbit 
$\mc{O}$, with radius $\rho$, the target being located at its center. On each orbit only a discrete number of positions accounts for the data acquisition process. Precisely, $M$ different positions are sampled along $\mc{O}$, and for each position $s$ the corresponding electric signal $u^{(s)} - H^{(s)}$ is measured by $N_r$ receptors on the skin, $\{ x_r^{(s)} \}_{r = 1}^{N_r}$. Here $u^{(s)}$ and $H^{(s)}$ denote the solution to \eqref{eq:model_u} and the function defined by \eqref{eq:H_fun}, associated to the position $s$, respectively.

This type of architecture resambles a multi-static SIMO (Single-Input Multi-Output) system.

\begin{center}
\begin{tabular}{@{} *5l @{}}    \toprule
	\emph{Symbol} & \emph{Meaning}   \\\midrule
	$\Omega_{s}$   & Fish body \\ 
	$\mathbf{p}_{s}$ & dipole moment \\ 
	$\zeta_{s}$ & electric organ \\ 
	$x_r^{(s)}$  & $r$-th receptor \\
	$u^{(s)}$  & electric potential solution to \eqref{eq:model_u} \\
	$H^{(s)}$  & function defined in \eqref{eq:H_fun} \\\midrule
	$N_r$  & number of receptors \\
	$M$  & number of positions \\
	$\mc{O}$ & circular orbit\\\bottomrule
	\hline
\end{tabular}
\captionof{table}{Notation referred to position $s \in \{ 1 , ... , M \}$ on the orbit $\mc{O}$. \label{Tab:Tcr}}
\end{center}

For any orbit $\mc{O}$ we get an $N_r \times M$ matrix of data $\bfQ$, which is called Multi-Static Response (MSR) matrix, whose $(r,s)$-entry is defined as
\begin{equation}(\bfQ)_{r,s} = u^{(s)}(x_r^{(s)}) - H^{(s)}(x_r^{(s)}) . \label{eq:MSRentry}\end{equation}

Henceforth, we shall use the MATLAB colon notation for specifying sub-matrices of a given matrix.
For instance, given matrix $\bfX$, we shall denote by $\bfX_{i,:}$ [resp. $\bfX_{:,j}$] the $i$-th row [resp. the $j$-th column] of $\bfX$.

\subsection{Data acquisition operator}

In order to simplify the notation, without loss of generality, we assume that the dielectric object is centered at the origin, and that the impedance of the fish is $\xi = 0$.

We recall the following theorem which provides an expansion of \eqref{eq:MSRentry},  see \cite{Ammari11652}.

\begin{theorem} \label{thm:MSR} Consider $M$ different positions of the fish along the circular orbit $\mc{O}$ of radius $\rho$, with $\rho$ large enough, indexed by $s = 1 , ... , M$. Let $\{ x_r^{(s)} \}_{r=1}^{N_r}$ be a set of receptors distributed on $\p \Omega_s$, the dipole located at $\zeta_s \in \mc{O}$ with dipole moment $\mathbf{p}_s$, and $K \ge 1$. Then the following expansion holds:	
	\begin{equation} \label{eq:q_sr} 
	u^{(s)}(x_r^{(s)}) - H^{(s)}(x_r^{(s)})=  \sum_{m+n = 1}^{K+1} \underbrace{\left [ \begin{matrix} A_{s,m} & B_{s,m} \end{matrix} \right ]}_{\bfS_{s,m}}\, \underbrace{\begin{bmatrix}
		M_{mn}^{cc} & M_{mn}^{cs}\\ M_{mn}^{sc} & M_{mn}^{ss}
		\end{bmatrix}}_{\bfM_{mn}}\, \underbrace{\begin{bmatrix}  \cos n \theta_{x_r^{(s)}} \\ \sin n \theta_{x_r^{(s)}} \end{bmatrix} \, \frac{-1}{2 \pi n r^n_{x_r^{(s)}}}}_{\bfG_{rn}^{(s)\top}} + {O}(\delta^{K+2}) ,
	\end{equation}
	where $M_{mn}^{cc}, M_{mn}^{cs} , M_{mn}^{sc}$ and $M_{mn}^{ss}$ are as in Definition \ref{def:CGPTs}, $ r = 1 , ... , N_r$,
	\begin{equation} \label{eq:AandB} 
	\begin{split}
	A_{s,m} &= - \frac{(-1)^m}{2 \pi} \mathbf{p}_s \cdot 
 \begin{bmatrix} 
	\phi_{m+1}(\zeta_s) \\ \psi_{m+1} (\zeta_s) \end{bmatrix}
		- \frac{1}{2 \pi m} \int_{\p \Omega_s} \frac{\p u^{(s)}}{\p \nu} \biggr |_{+} (y) \phi_m(y) \text{ d} \sigma_y ,\\
	B_{s,m} &= \frac{(-1)^m}{2 \pi} \mathbf{p}_s \cdot 
	 \begin{bmatrix} 
	- \psi_{m+1}(\zeta_s) \\ \phi_{m+1} (\zeta_s) \end{bmatrix}
	- \frac{1}{2 \pi m} \int_{\p \Omega_s} \frac{\p u^{(s)}}{\p \nu} \biggr |_{+} (y) \psi_m (y) \text{ d} \sigma_y ,\\
	\end{split}
	\end{equation}
	
	\begin{equation*} \label{eq:phiandpsi} 
	\phi_m(x) = \frac{\cos(m \theta_x)}{r^m_x} ,\qquad
	\psi_m(x) = \frac{\sin(m \theta_x)}{r^m_x},
	\end{equation*}
	and
	\begin{equation} \label{eq:Mblocks} 
	\mathbb{M}^{(K)} = \bbM = \begin{bmatrix}
	\bfM_{11} & \bfM_{12} & \dots & \bfM_{1K} \\
	\bfM_{21} &           &  \iddots &    \bfzero       \\
	\vdots &      \iddots     & \iddots &      \vdots     \\
	\bfM_{K1}  &      \bfzero     & \dots &      \bfzero     \\
	\end{bmatrix}
	\end{equation}
	is the upper anti-diagonal block matrix of the CGPTs of order $\le K$. Here, $\top$ denotes the transpose of a matrix.
\end{theorem}

We define $\varepsilon = \delta/\rho$ the length-scale associated to the orbit $\mc{O}$, i.e., the ratio between the size of the target and the distance $\rho$.

A more careful analysis of the reminder in formula \eqref{eq:q_sr} shows that the remainder can be expressed in term of the length-scale $\varepsilon$, and written as $O(\varepsilon^{K+2})$. See Appendix \ref{apx:reminder}.

By Theorem \ref{thm:MSR}, the rows of $\bfQ$ admit the following expansions: 

\begin{equation} \label{eq:regression_s}
(\bfQ)_{\,: \, , s}  = \mathcal{L}^{(s)}(\bbM^{(K)}) + \mathbf{E}_{\,: \, , s}, \qquad \| \mathbf{E}_{\,: \, , s}\|_\infty = O(\varepsilon^{K+2}) , \qquad s = 1 , ... , M,
\end{equation}
where $\mc{L}^{(s)} : \Mat{2K}{2K} \arr \RR^{N_r}$ is the linear map defined by  \eqref{eq:q_sr}, i.e., $\mc{L}^{(s)}(\mathbb{M}) = \mathbf{G}^{(s)} \, \mathbb{M} \,\mathbf{S}_{s,\,:}^\top \,$, $K$ is the truncation order, and $\varepsilon = \delta/\rho$ is the length-scale associated to the orbit $\mc{O}$. Thus, we can write the expansion of the complete MSR matrix as follows:
\begin{equation} \label{eq:regression_eq}
\bfQ = \mathcal{L}(\bbM^{(K)}) + \mathbf{E}, \qquad \|\mathbf{E} \|_\infty = O(\varepsilon^{K+2}) .
\end{equation}
The linear map $\mathcal{L} : \Theta \subseteq \Mat{2K}{2K} \arr \Mat{N_r}{M}$ is the truncated output (or forward) operator.

The acquisition operator $\mathcal{L}(\mathbb{M})$ is defined by \eqref{eq:q_sr}. More precisely, it can be written as
%\begin{equation*} \label{eq:L_acquisition_block} 
%\mathcal{L}(\mathbb{M}) = \left [\begin{matrix} 
%\mathcal{L}^{(1)}(\mathbb{M} ) \\
%\mathcal{L}^{(2)}(\mathbb{M} ) \\
%\vdots \\
%\mathcal{L}^{(M)}(\mathbb{M} )
%  \end{matrix} \right ] =  
%  \left [\begin{matrix} 
%  \mathbf{S}_{1,\,:}\, \mathbb{M} \, \bfG^{(1)\top}\\
%  \mathbf{S}_{2,\,:}\, \mathbb{M} \,\bfG^{(2)\top}\\
%  \vdots \\
%  \mathbf{S}_{M,\,:}\, \mathbb{M} \, \bfG^{(M)\top}
%  \end{matrix} \right ] ,
%\end{equation*}
\begin{equation*} \label{eq:L_acquisition_block} 
\mathcal{L}(\mathbb{M}) = \left [\begin{matrix} 
\mathcal{L}^{(1)}(\mathbb{M} ) &
\mathcal{L}^{(2)}(\mathbb{M} ) &
\dots &
\mathcal{L}^{(M)}(\mathbb{M} )
\end{matrix} \right ] =  \begin{bmatrix} 
\mathbf{G}^{(1)} \, \mathbb{M} \,\mathbf{S}_{1,\,:}^\top &   \dots  & 
\mathbf{G}^{(M)} \, \mathbb{M} \,\mathbf{S}_{M,\,:}^\top 
\end{bmatrix} ,
\end{equation*}
where
%\begin{equation*} \label{eq:Li_acquisition} 
%\mathcal{L}_i(\mathbb{M})  = 
%\rectangle[2pt]{45pt}{7pt} \times \rectangle[15pt]{45pt}{30pt} \times  \rectangle[15pt]{30pt}{30pt} .
%\end{equation*}
\begin{equation*} \label{eq:Li_acquisition} 
\mathcal{L}^{(s)}(\mathbb{M})  =  \rectangle[60pt]{45pt}{30pt} \times \rectangle[15pt]{45pt}{30pt} \times  
\rectangle[15pt]{7pt}{30pt} .
\end{equation*}
We define block matrices $\mathbf{G} \in \Mat{M}{1}(\Mat{N_r}{2K})$, $\mathbf{S} \in \Mat{M}{1}(\Mat{1}{2K})$ by vertically stacking  the matrices, as in \figref{fig:blocks}.

\begin{figure}\centering
$\mathbf{G} = \boxed{ \begin{matrix}\mathbf{G}^{(1)}\\ \mathbf{G}^{(2)}\\\vdots\\ \mathbf{G}^{(M)}\end{matrix} } = \boxed{ \begin{matrix}\rectangle[15pt]{100pt}{30pt}^\top\\ \rectangle[15pt]{100pt}{30pt}^\top \\\vdots\\ \rectangle[15pt]{100pt}{30pt}^\top\end{matrix} } \qquad \mathbf{S} = \boxed{ \begin{matrix}\mathbf{S}_{1,\,:}\\ \mathbf{S}_{2,\,:}\\\vdots\\ \mathbf{S}_{M,\,:}\end{matrix} } = \boxed{ \begin{matrix}\rectangle[0pt]{45pt}{7pt}\\ \rectangle[0pt]{45pt}{7pt} \\\vdots\\ \rectangle[0pt]{45pt}{7pt} \end{matrix} } $
  \caption{}
\label{fig:blocks}
\end{figure}

%\begin{figure}[H]
%	\centering
%\begin{tikzpicture}
%\node[anchor=south west,inner sep=0] (image) at (0,0) {\includegraphics[width=0.5\textwidth]{ee}};
%\begin{scope}[x={(image.south east)},y={(image.north west)}]
%\draw[red, thick] (0.47,0.47) node[left=0.2cm,above=0.2cm] {$\mathsf{\delta}$} rectangle (0.58,0.58);
%\draw[->, dashed, thick, blue,  arrows={-latex}]  (0.52,0.52) -- (0.77,0.89) node[sloped,midway,above=-0.1cm] {$\mathsf{\rho}$};
%\end{scope}
%\end{tikzpicture}
%\end{figure}

We are interested in estimating the matrix parameter $\mathbb{M}$ from the MSR matrix $\bfQ$. Therefore, we aim at solving the following minimization problem
\begin{equation} \label{eq:lst_squares} 
\min_{\mathbb{M} \perp \text{ker} (\mc{L})} \|\mc{L}(\mathbb{M})  - \bfQ \|_{\text{F}},
\end{equation} 
where $\| \cdot \|_{\text{F}}$ denotes the Frobenius norm of a matrix.
%The least-squares estimation problem is 
%\begin{equation*} \label{eq:lst_squares_est} 
%\bbM^{\text{est}} = \arg \min_{\mathbb{M} \perp \text{ker} (\mc{L})} \|\mc{L}(\mathbb{M})  - \bfQ \|_{\text{F}}.
%\end{equation*} 

\subsection{Generalized Kronecker form of the forward operator}

In this section we vectorize the data acquisition operator $\mc{L} : \Mat{2K}{2K} \longrightarrow \Mat{N_r}{M}$ in order to find a matrix representation.

\begin{lemma}\label{lem:kr_design_matrix} The operator $\mc{L}$ defined by \eqref{eq:L_acquisition_block} can be represented in a vectorized form employing the product defined in Definition \ref{def:kronecker_mitra}:
	\begin{equation*} %\label{eq:L_acquisition} 
	\vec{\mc{L}(\mathbb{M})} = (\bfS \kron \{\bfG^{(s)}\}) \, \vec{\mathbb{M}}.
	\end{equation*}
\end{lemma}
\begin{proof} By definition,
	\begin{equation*} %\label{eq:L_acquisition} 
	\mathcal{L}(\mathbb{M}) = \begin{bmatrix} 
	\mathbf{G}^{(1)} \, \mathbb{M} \,\mathbf{S}_{1,\,:}^\top &   \dots  & 
	\mathbf{G}^{(M)} \, \mathbb{M} \,\mathbf{S}_{M,\,:}^\top 
	\end{bmatrix} .
	\end{equation*}
	Therefore
	\begin{equation*} %\label{eq:L_acquisition} 
	\begin{split}
	\vec{\mc{L}(\mathbb{M})}  & = 
	\mbox{vec} \begin{bmatrix} 
	\mathbf{G}^{(1)} \, \mathbb{M} \,\mathbf{S}_{1,\,:}^\top &   \dots  & 
	\mathbf{G}^{(M)} \, \mathbb{M} \,\mathbf{S}_{M,\,:}^\top 
	\end{bmatrix} \\ &=   
	\begin{bmatrix} \bfG^{(1)} \, \mathbb{M} \,\mathbf{S}_{1,\,:}^\top \\ \vdots \\ \bfG^{(M)}\, \mathbb{M} \,\mathbf{S}_{M,\,:}^\top\end{bmatrix}
	\\ & =   
	\begin{bmatrix}(\mathbf{S}_{1,\,:} \kron \bfG^{(1)}) \, \vec{\mathbb{M}} \\ \vdots \\ (\mathbf{S}_{M,\,:} \kron \bfG^{(M)}) \, \vec{\mathbb{M}} 
	\end{bmatrix}
	\\ & =   
	\begin{bmatrix}\mathbf{S}_{1,\,:} \kron \bfG^{(1)}\\ \vdots \\ \mathbf{S}_{M,\,:} \kron \bfG^{(M)}
	\end{bmatrix} \vec{\mathbb{M}}
	\\ & =  ( \,\bfS \kron  \{\bfG^{(s)}\} \,) \, \vec{\mathbb{M}}.
	\end{split}
	\end{equation*}	
\end{proof}

Notice that the matrix $\bfL \in \Mat{M N_r}{4 K^2}$ defined by
\begin{equation}\label{eq:matrix_L}
\bfL \deff \bfL_{\mathcal{L}}  =  \bfS \kron \{\bfG^{(s)}\} 
\end{equation}
is the unique $M N_r \times 4 K^2$ matrix such that $\vec{\mc{L}(\bfX)} = \bfL_{\mc{L}} \vec{\bfX}$, for all $\bfX \in \Mat{2K}{2K}$.

\medskip

Hence, minimization problem \eqref{eq:lst_squares} assumes the following form 
\begin{equation} \label{eq:lst_squares_vec} 
\min_{\bbM} \| \bfL_{\mc{L}} \vec{\bbM} - \vec{\bfQ}\|_2 ,
\end{equation}
We aim at seeking a vector $\vec{\wh{\bbM}}$ which is optimal in the least-squares sense.

As it is well known, the standard least-squares estimator for \eqref{eq:lst_squares_vec} is given by the Moore-Penrose inverse of $\bfL$, denoted by $\bfL^\dagger $. If $\bfL$ is full column rank, than
\begin{equation}\label{eq:MPinverse}
\vec{\wh{\bbM}}_{MP} = \bfL^\dagger \vec{\bfQ} = (\bfL^\top \bfL)^{-1} \bfL^\top \vec{\bfQ} .
\end{equation}

However, the special block Kronecker form of $\bfL$ suggests to employ the following generalized inverse \cite{Marco2018}.

\begin{theorem} \label{thm:generalized_inverse}If $\bfS$ and $\bfG^{(s)}$ for $s = 1, ... , M$, are full column rank, then 
	\begin{equation} \label{eq:minnormsol} \mathfrak{L} \deff \bfS^\emph{\dagger} \kron_C \{\bfG^{(s)\emph{\dagger}}\} \end{equation}
	is a reflexive minimum norm g-inverse of $\bfL_{\mc{L}} = \bfS \kron \{\bfG^{(s)}\}$. Here $\kron_C$ denotes the column-wise generalized Kronecker product defined in \ref{def:col_gkp}, and \emph{$^\dagger$} denotes the Moore–Penrose inverse.	
\end{theorem}

\begin{proof} The proof is readily obtained by noticing that
	\begin{equation*}  (\bfS^{\dagger} \kron_C \{\bfG^{(s){\dagger}}\}) (\bfS \kron \{\bfG^{(s)}\}) = \bfI_{4 K^2}. \end{equation*}
\end{proof}

This particular generalized inverse is useful for solving \eqref{eq:lst_squares_vec} when $\vec{\bfQ}$ lies in the range of $\bfL_{\mc{L}}$ \cite{Marco2018}. Notice that $\mathfrak{L}$ is not the same as $\bfL^{\dagger}$ in general.

\begin{figure}[ht]
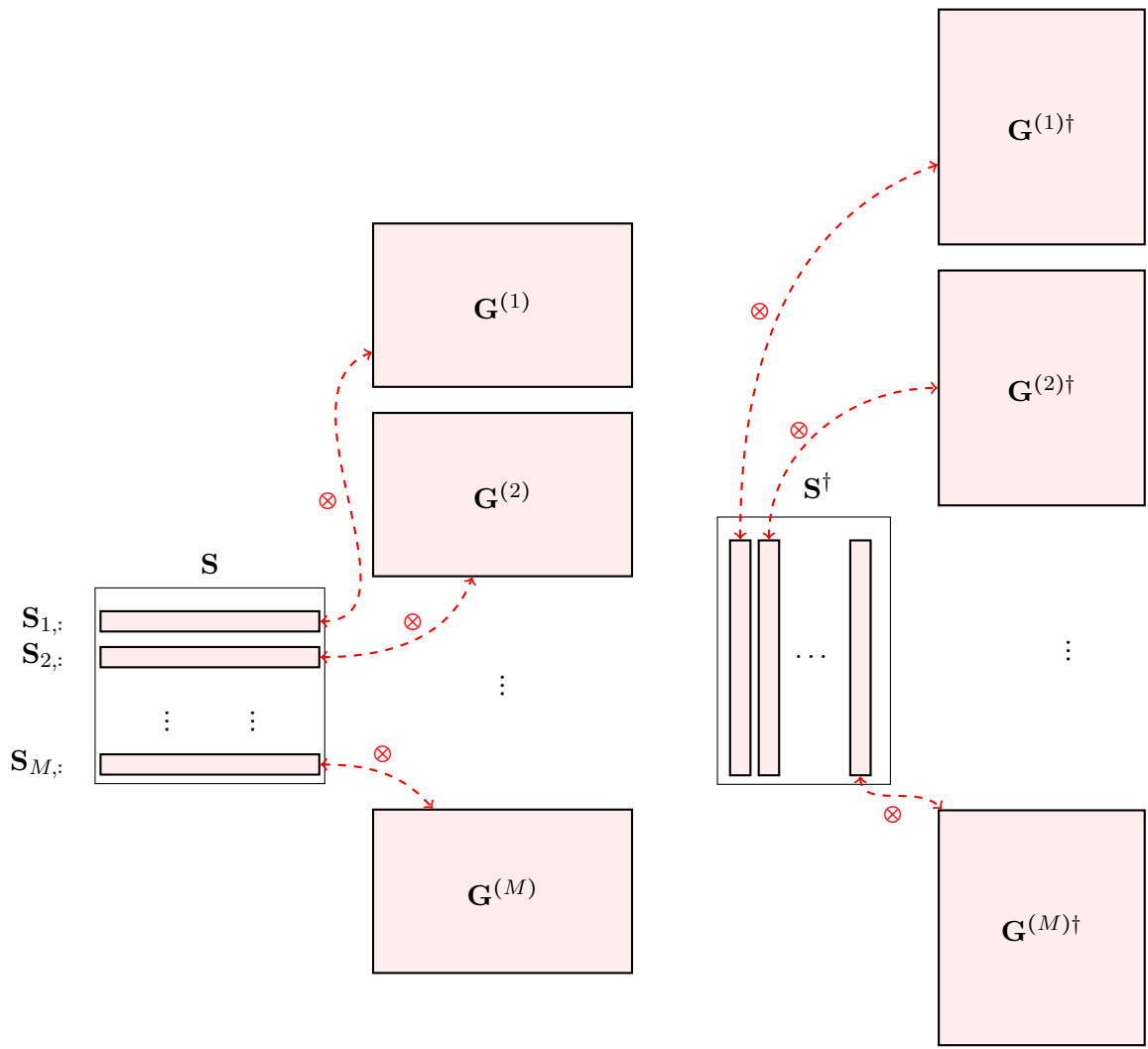

	\begin{multicols}{2}	
		\begin{align*} \\ \\ \\
		& 
		\begin{matrix} \tikzmark{left}{}  & & & & &  & & & &     \\ & & & & &   & & & &   \\  & & & &  & & & & &   \\ & & & & &  & & & & \\  & & & & & & & & &\tikzmark{right}{}  \end{matrix} \Highlight[G1]\\ 
		& \begin{matrix} \tikzmark{left}{}  & & & & &   & & & &     \\ & & & & &  & & & &   \\  & & & &  & & & & &    \\  & & & &  & & & & &   \\  & & & & &   & & & & \tikzmark{right}{}  \end{matrix} \Highlight[G2] \\
		\boxed{ \;
			\begin{matrix} 
			\tikzmark{left}{}& & & & & & & \tikzmark{right}{}  \Highlight[S1] \\ 
			\tikzmark{left}{}& & & & & & & \tikzmark{right}{}  \Highlight[S2] \\ 
			\tikzmark{left}{}& & & & & & & \tikzmark{right}{} \\ 
			& &\vdots & & &  \vdots &  & \\
			\tikzmark{left}{}& & & & & & & \tikzmark{right}{} \Highlight[SN] \end{matrix} 
			\;} \qquad  & \quad \qquad \quad \vdots \\ 
		& \begin{matrix} \tikzmark{left}{}  & & & & &   & & & &     \\ & & & & & & & & &    \\  & & & &  & & & & &    \\  & & & &  & & & & &    \\ & & & & &  & & & &   \tikzmark{right}{}  \end{matrix} \Highlight[GN] 
		\end{align*} 
		
		\tikz[overlay,remember picture] {
			\draw[<->,thick,red,dashed] (S1) to  [out=0, in=200] node [left] {$\kron$} (G1);	
			\node[left=1.8cm] at (S1) {$\bfS_{1,:}$};	
			\node[above=0.5cm] at (S1) {$\bfS$};
			\node[left=1.8cm] at (S2) {$\bfS_{2,:}$};
			\node[left=1.8cm] at (SN) {$\bfS_{M,:}$};
			\node at (G1) {$\bfG^{(1)}$};
		}
		\tikz[overlay,remember picture] {
			\draw[<->,thick,red,dashed] (S2) to  [out=0, in=250] node [above] {$\kron$} (G2);
			\node at (G2) {$\bfG^{(2)}$};
		}
		\tikz[overlay,remember picture] {
			\draw[<->,thick,red,dashed] (SN) to  [out=0, in=130] node [above] {$\kron$} (GN);
			\node at (GN) {$\bfG^{(M)}$};
		}
		\break
		\begin{align*}
		& 
		\begin{matrix} \tikzmark{left}{}   & & & & &   & &     \\  & & & & &   & &    \\   & & & & &   & &   \\  & & & & &   & & \\   & & & & &   & &  \\ & & & & &   & &  \\ & & & & &   & &  \tikzmark{right}{}  \end{matrix} \Highlight[G1p]\\ 
		& \begin{matrix} \tikzmark{left}{}  & & & & &   & &     \\ & & & & &   & &   \\   & & & & &   & &    \\  & & & & &   & &  \\  & & & & &   & &   \\ & & & & &   & &  \\ & & & & &   & & \tikzmark{right}{}  \end{matrix} \Highlight[G2p] \\
		\boxed{ \;\;
			\begin{matrix} \tikzmark{left}{}  
			\\  \\   \\ \\ \\ \\   \tikzmark{right}{} \end{matrix} \Highlight[S1p] \;\;\;\;
			\begin{matrix} \tikzmark{left}{}  
			\\   \\ \\ \\ \\  \\ \tikzmark{right}{} \end{matrix} \Highlight[S2p]\;\;\;\dots \;\;\;\;
			\begin{matrix} \tikzmark{left}{}  
			\\  \\ \\ \\  \\ \\   \tikzmark{right}{} \end{matrix} \Highlight[SNp]\;\;\;
		} \qquad  & \quad \qquad \quad \vdots \\ 
		& \begin{matrix} \tikzmark{left}{} & & & & &   & &    \\ & & & & &   & &    \\ & & & & &   & &   \\  & & & & &   & &    \\ & & & & &   & &   \\ & & & & &   & &  \\ & & & & &   & &   \tikzmark{right}{}  \end{matrix} \Highlight[GNp] 
		\end{align*} 
		
		\tikz[overlay,remember picture] {
			\draw[<->,thick,red,dashed] (S1p) to  [out=90, in=200] node [left] {$\kron$} (G1p);	
			% 	\node[below=1cm] at (S1p) {$(\bfS)_{1,:}^{\dagger}$};	
			\node[above=2.3cm,right=0.7cm] at (S1p) {$\bfS^{\dagger}$};
			% 	\node[below=1cm] at (S2p) {$(\bfS)_{2,:}^{\dagger}$};
			% 	\node[below=1cm] at (SNp) {$(\bfS)_{M,:}^{\dagger}$};
			\node at (G1p) {$\bfG^{(1) \dagger}$};
		}
		\tikz[overlay,remember picture] {
			\draw[<->,thick,red,dashed] (S2p) to  [out=90, in=180] node [left] {$\kron$} (G2p);
			\node at (G2p) {$\bfG^{(2) \dagger}$};
		}
		\tikz[overlay,remember picture] {
			\draw[<->,thick,red,dashed] (SNp) to  [out=270, in=130] node [below] {$\kron$} (GNp);
			\node at (GNp) {$\bfG^{(M) \dagger}$};
		}
	\end{multicols}
	\caption{On the left, the operator $\bfL$; on the right, its generalized inverse, i.e., $\mathfrak{L}$. }
	\label{fig:generalized-inverse}
\end{figure}

As we shall see later, the g-inverse given by \eqref{eq:minnormsol} is particularly suitable for establishing a bound on the reconstruction error as well as for designing a recursive online estimation of the GPTs. Figure \ref{fig:generalized-inverse} schematically shows the computation of $\mathfrak{L}$.

\subsection{Online reconstruction}

In this section we propose very simple formulas to efficiently perform an online reconstruction of the features.

%\[ \vecbadd{\bbM_\theta} = \begin{bmatrix} \mc{R}  (\theta_1) \\ \vdots \\ \mc{R} (\theta_{M}) \end{bmatrix} \vecbadd{\bbM}\]
%where
%\[ \mc{R} (\theta_\ell) = \begin{bmatrix} \bfR_{11} \\ &\bfR_{21} \\ & & \bfR_{12} \\ & & & \ddots \\ & & & & \bfR_{K1}\end{bmatrix} \]
%is orthogonal, with $\bfR_{mn} = \bfR(m \theta_\ell) \kron \bfR(n \theta_\ell)$ $4\times4$ orthogonal matrices.
%...
%
%\[ \bfL  = \begin{bmatrix} \bfL_0 \mathcal{R} (\theta_1) \\ \vdots \\ \bfL_0 \mathcal{R}(\theta_{M}) \end{bmatrix} \]
%
%\[ \bfL  = ( \bfI_M \kron \bfL_0 ) \begin{bmatrix} \mathcal{R} (\theta_1) \\ \vdots \\ \mathcal{R}(\theta_{M}) \end{bmatrix}  = \bfL_0 \bfV^\top \bfD_1 \bfV + ... +   \bfL_0 \bfV^\top \bfD_M \bfV  = ...\]
%
%...
%
%
%
%
%Another method:

\begin{figure}[H]
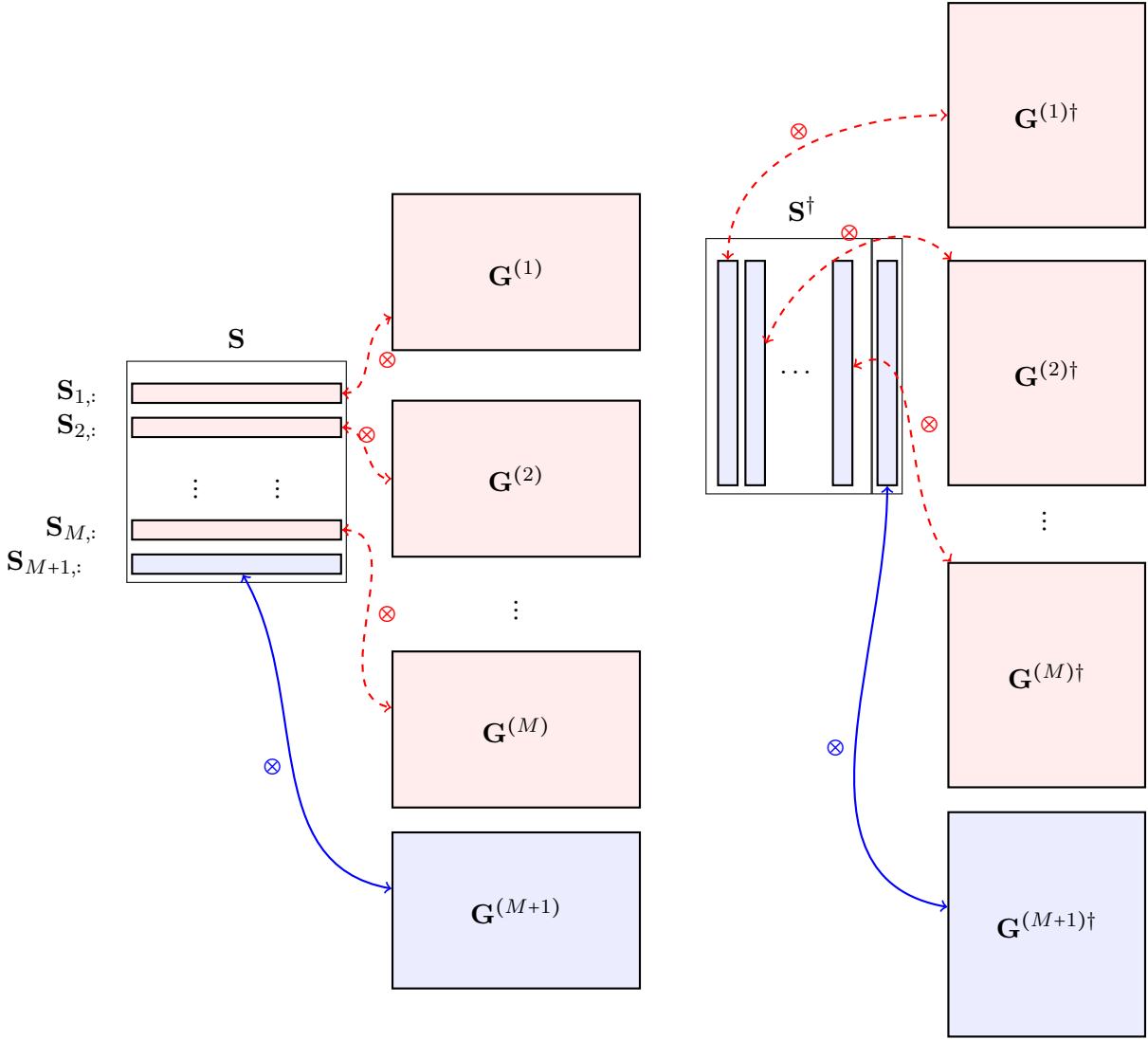

	\begin{multicols}{2}	
		\begin{align*} \\ \\ \\ \\
		& 
		\begin{matrix} \tikzmark{left}{}  & & & & &  & & & &     \\ & & & & &   & & & &   \\  & & & &  & & & & &   \\ & & & & &  & & & & \\  & & & & & & & & &\tikzmark{right}{}  \end{matrix} \Highlight[G1]\\ 
		\boxed{ \;
			\begin{matrix} 
			\tikzmark{left}{}& & & & & & & \tikzmark{right}{}  \Highlight[S1] \\ 
			\tikzmark{left}{}& & & & & & & \tikzmark{right}{}  \Highlight[S2] \\ 
			\tikzmark{left}{}& & & & & & & \tikzmark{right}{} \\ 
			& &\vdots & & &  \vdots &  & \\
			\tikzmark{left}{}& & & & & & & \tikzmark{right}{} \Highlight[SN] \\
			\tikzmark{left}{}& & & & & & & \tikzmark{right}{} \Highlightb[Snew] \end{matrix}  
			\;} \qquad 	& \begin{matrix} \tikzmark{left}{}  & & & & &   & & & &     \\ & & & & &  & & & &   \\  & & & &  & & & & &    \\  & & & &  & & & & &   \\  & & & & &   & & & & \tikzmark{right}{}  \end{matrix} \Highlight[G2] \\
		& \quad \qquad \quad \vdots \\ 
		& \begin{matrix} \tikzmark{left}{}  & & & & &   & & & &     \\ & & & & & & & & &    \\  & & & &  & & & & &    \\  & & & &  & & & & &    \\ & & & & &  & & & &   \tikzmark{right}{}  \end{matrix} \Highlight[GN] \\
		& \begin{matrix} \tikzmark{left}{}  & & & & &   & & & &     \\ & & & & & & & & &    \\  & & & &  & & & & &    \\  & & & &  & & & & &    \\ & & & & &  & & & &   \tikzmark{right}{}  \end{matrix}   \Highlightb[Gnew] 
		\end{align*} 
		
		\tikz[overlay,remember picture] {
			\draw[<->,thick,red,dashed] (S1) to [out=0, in=200] node [right] {$\kron$} (G1);	
			\node[left=1.8cm] at (S1) {$\bfS_{1,:}$};	
			\node[above=0.5cm] at (S1) {$\bfS$};
			\node[left=1.8cm] at (S2) {$\bfS_{2,:}$};
			\node[left=1.8cm] at (SN) {$\bfS_{M,:}$};
			\node at (G1) {$\bfG^{(1)}$};
		}
		\tikz[overlay,remember picture] {
			\draw[<->,thick,red,dashed] (S2) to [out=0, in=180] node [above] {$\kron$} (G2);
			\node at (G2) {$\bfG^{(2)}$};
		}
		\tikz[overlay,remember picture] {
			\draw[<->,thick,red,dashed] (SN)  to [out=0, in=170] node [right] {$\kron$} (GN) ;
			\node at (GN) {$\bfG^{(M)}$};
		}
		
		\tikz[overlay,remember picture] {
			\draw[<->,thick,blue,solid] (Snew)  to [out=300, in=170] node [left] {$\kron$} (Gnew);
			\node at (Gnew) {$\bfG^{(M+1)}$};
			\node[left=2cm] at (Snew) {$\bfS_{M+1 ,:}$};
		}
		\break
		\begin{align*}
		& 
		\begin{matrix} \tikzmark{left}{}   & & & & &   & &     \\  & & & & &   & &    \\   & & & & &   & &   \\  & & & & &   & & \\   & & & & &   & &  \\ & & & & &   & &  \\ & & & & &   & &  \tikzmark{right}{}  \end{matrix} \Highlight[G1p]\\ 
		\boxed{ \;\;
			\begin{matrix} \tikzmark{left}{}  
			\\  \\   \\ \\ \\ \\   \tikzmark{right}{} \end{matrix} \Highlightb[S1p] \;\;\;\;
			\begin{matrix} \tikzmark{left}{}  
			\\   \\ \\ \\ \\  \\ \tikzmark{right}{} \end{matrix} \Highlightb[S2p]\;\;\;\dots \;\;\;\;
			\begin{matrix} \tikzmark{left}{}  
			\\  \\ \\ \\  \\ \\   \tikzmark{right}{} \end{matrix} \Highlightb[SNp]\;\;\;
		} \boxed{\; \begin{matrix} \tikzmark{left}{}  
			\\  \\ \\ \\  \\ \\   \tikzmark{right}{} \end{matrix} \;\Highlightb[Snewp]} \qquad  & \begin{matrix} \tikzmark{left}{}  & & & & &   & &     \\ & & & & &   & &   \\   & & & & &   & &    \\  & & & & &   & &  \\  & & & & &   & &   \\ & & & & &   & &  \\ & & & & &   & & \tikzmark{right}{}  \end{matrix} \Highlight[G2p] \\
		& \quad \quad \quad \vdots \\ 
		& \begin{matrix} \tikzmark{left}{} & & & & &   & &    \\ & & & & &   & &    \\ & & & & &   & &   \\  & & & & &   & &    \\ & & & & &   & &   \\ & & & & &   & &  \\ & & & & &   & &   \tikzmark{right}{}  \end{matrix} \Highlight[GNp] \\
		& \begin{matrix} \tikzmark{left}{} & & & & &   & &    \\ & & & & &   & &    \\ & & & & &   & &   \\  & & & & &   & &    \\ & & & & &   & &   \\ & & & & &   & &  \\ & & & & &   & &   \tikzmark{right}{}  \end{matrix} \Highlightb[Gnewp] 
		\end{align*} 
		
		\tikz[overlay,remember picture] {
			\draw[<->,thick,red,dashed] (S1p) to [out=90, in=180]  node [above] {$\kron$} (G1p);	
			% 	\node[below=1cm] at (S1p) {$(\bfS)_{1,:}^{\dagger}$};	
			\node[above=2.3cm,right=0.7cm] at (S1p) {$\bfS^{\dagger}$};
			% 	\node[below=1cm] at (S2p) {$(\bfS)_{2,:}^{\dagger}$};
			% 	\node[below=1cm] at (SNp) {$(\bfS)_{M,:}^{\dagger}$};
			\node at (G1p) {$\bfG^{(1) \dagger}$};
		}
		\tikz[overlay,remember picture] {
			\draw[<->,thick,red,dashed] (S2p)  to [out=70, in=130] node [above] {$\kron$} (G2p);
			\node at (G2p) {$\bfG^{(2) \dagger}$};
		}
		\tikz[overlay,remember picture] {
			\draw[<->,thick,red,dashed] (SNp)  to [out=30, in=130] node [right] {$\kron$} (GNp);
			\node at (GNp) {$\bfG^{(M) \dagger}$};	
			\node at (Gnewp) {$\bfG^{(M+1) \dagger}$};
		}
		\tikz[overlay,remember picture] {
			\draw[<->,thick,blue,solid] (Snewp)  to [out=270, in=170] node [left] {$\kron$} (Gnewp) ;
		}
	\end{multicols}
	\caption{On the left, the augmented operator $\bfL_{M+1}$; on the right: its generalized inverse, i.e., $\mathfrak{L}_{M+1}$. The parts which change are highlighted in blue.}
	\label{fig:update-gi}
\end{figure}

By inspecting the form of the g-inverse $\mathfrak{L}$ given by \eqref{eq:minnormsol} it is easy to see that, when a new position becomes available, the pseudoinverse of the augmented source matrix $\bfS$ is the only term which needs to be recomputed. As shown in Figure \ref{fig:update-gi}, the pseudoinverses of the matrices $\bfG^{(s)}$ corresponding to different positions intervene in $\mathfrak{L}$ without interfering with each other. Therefore we have the following result.

\begin{lemma} Let us denote $\mathfrak{L}_M$ the generalized inverse given by \eqref{eq:minnormsol} for $M \gg 1$ positions, and let $\bfS_{M+1 , \,:}$ be full column rank. Then
	\begin{equation} \label{eq:minnormsol_update} \mathfrak{L}_{M+1} = \begin{bmatrix}\bfS_{1:M,\,:}^\emph{\dagger}- \mathbf{K}_{M+1} \mathbf{d}_{M+1} & | & \mathbf{K}_{M+1} \end{bmatrix} \kron_C \{\bfG^{(s)\emph{\dagger}}\},\end{equation}
	where 
	\[\mathbf{d}_{M+1} \deff \bfS_{M+1,\,:}\, \bfS_{1:M,\,:}^\emph{\dagger}  \,\, , \]  and \[\mathbf{K}_{M+1} \deff (1 + \mathbf{d}_{M+1} \mathbf{d}_{M+1}^\top)^{-1}  \bfS_{1:M,\,:}^\emph{\dagger} \, \mathbf{d}_{M+1}^\top . \]
\end{lemma}
%\begin{lemma} Let us denote $\mathfrak{L}_M$ the generalized inverse given by \eqref{eq:minnormsol} for $M \gg 1$ positions. Then
%	\begin{equation} \label{eq:minnormsol_update} \mathfrak{L}_{M+1} = \begin{bmatrix}\bfS_{1:M,\,:}^\emph{\dagger}- \mathbf{K}_{M+1} \mathbf{d}_{M+1} & | & \mathbf{K}_{M+1} \end{bmatrix} \kron_C \{\bfG^{(s)\emph{\dagger}}\},\end{equation}
%	where 
%	\[\mathbf{d}_{M+1} = \bfS_{M+1,\,:}\, \bfS_{1:M,\,:}^\emph{\dagger}  \,\, , \] \[\mathbf{c}_{M+1} = \bfS_{M+1,\,:} - \mathbf{d}_{M+1}\,\bfS_{1:M,\,:}  \,\, , \] and \[\mathbf{K}_{M+1} = \begin{cases} \mathbf{c}_{M+1}^\emph{\dagger} & \mbox{if } \mathbf{c}_{M+1} \ne 0 \,. \\ (1 + \mathbf{d}_{M+1} \mathbf{d}_{M+1}^\top)^{-1}  \bfS_{1:M,\,:}^\emph{\dagger} \, \mathbf{d}_{M+1}^\top
%	& \mbox{if } \mathbf{c}_{M+1} = 0 \,.\end{cases} \]
%\end{lemma}
\begin{proof} Appending new positions affects only the factor $\bfS$, which can be updated by means of Greville's recursive formula for the pseudoinverse, see \cite{doi:10.1137/S0895479801388194}.
\end{proof}
\FloatBarrier

\FloatBarrier

\section{Analysis of the design matrix}

\label{sec:rank-analysis}

In this section we want to analyze in detail the form of the acquisition operator. As a result, we provide an in-depth study of the reconstruction.

Minimization problem \eqref{eq:lst_squares} indicates that the null-space of the forward operator $\mc{L}_{\mc{L}}$ we studied so far is related to the capability of uniquely reconstructing the CGPTs, and, in the end, to the classification of a dielectric target.

\subsection{Matrix of receptors}

 The matrix of receptors associated to the $s$-th position is given by 
\[ \bfG^{(s)} = \begin{bmatrix} \dfrac{\cos(\theta_1)}{r_1} &   \dfrac{\sin(\theta_1)}{r_1}  &  \dfrac{\cos(2\theta_{1})}{2 r_{1}^2} &  \dfrac{\sin(2\theta_{1})}{2 r_{1}^2} &   \dots &  \dfrac{\cos(K\theta_{1})}{K r_{1}^K} &  \dfrac{\sin(K\theta_{1})}{K r_{1}^K} \\
\dfrac{\cos(\theta_2)}{r_2} &  \dfrac{\sin(\theta_{2})}{r_{2}}  &   \dfrac{\cos(2\theta_2)}{2 r_2^2} &  \dfrac{\sin(2\theta_2)}{2 r_2^2}  &   \dots & \dfrac{\cos(K\theta_{2})}{K r_{2}^K} &  \dfrac{\sin(K\theta_{2})}{K r_{2}^K} \\
\vdots & & \vdots   & &   \ddots &  & \vdots \\
\dfrac{\cos(\theta_{N_r})}{r_{N_r}} &  \dfrac{\sin(\theta_{N_r})}{r_{N_r}}  &  \dfrac{\cos(2\theta_{N_r})}{2 r_{N_r}^2}  &  \dfrac{\sin(2\theta_{N_r})}{2 r_{N_r}^2}  &   \dots &  \dfrac{\cos(K\theta_{N_r})}{K r_{N_r}^K}  &  \dfrac{\sin(K\theta_{N_r})}{K r_{N_r}^K} \\
\end{bmatrix} .\]

In Appendix \ref{apx:uniqueness}, we show that  $\bfG^{(s)}$ is full column rank as soon as there are $2K \le  N_r$ distinct receptors that are a \emph{general configuration} in the sense of Remark  \ref{rem:abuseofdef}. Furthermore, we have the following Lemma \cite{Ammari2013}. 

\begin{lemma} $2K \leq N_r$ distinct points distributed along a circular arc are a general configuration.
\end{lemma}

It is clear that a single position yields a design matrix $\bfL_{\mathcal{L}}$ which is not full column rank, no matter how many receptors are considered.
However, collecting many electrostatic measurements at different positions ultimately enriches the column space of the matrix $\bfS$. As a matter of fact, if $\bfS$ and $\bfG^{(s)}$ for $s = 1, ... , M$, are full column rank, $\mathfrak{L}$ proves to be a left-inverse and thus $\bfL_{\mathcal{L}}$ is full column rank as well.

\medskip

\medskip

\subsection{Source vector}

The row vector $\bfS_{s,\,:}$, which is referred to the source corresponding to the $s$-th position, is defined by \eqref{eq:q_sr}. For simplicity we consider the case $\xi = 0$, see \eqref{eq:AandB}.

Denote by 
\[
\mathbf{p}_s^\perp =  \mathbf{p}_s \begin{bmatrix}
0 & 1 \\ -1 & 0 
\end{bmatrix}
\] 
the unit vector orthogonal to the dipole moment $\mathbf{p}_s = [\cos \alpha , \sin \alpha ]$, and $\zeta_s = \rho e^{i \overline{\theta}_s}$ be the location of the dipole.  

We can naturally split $\bfS_{s,\,:}$ into the pure dipole term and the distributed source term:

\[ \bfS_{s,\,:} = (\bfS_{dip})_{s,\,:} + (\bfS_{SL})_{s,\,:} \;.\]
Here $\bfS_{dip}$ and $\bfS_{SL}$ are given as follows. %where $(\bfS_{dip})_{s,\,:} :=$
Employing the product given in Definition \ref{eq:kron_vec}, define the block diagonal matrix 
\[ 
\bfP^{(s)}_K \deff 
\bfI_K \kron \left\{(-1)^\ell \begin{bmatrix}
\textbf{p}_s\\
\textbf{p}_s^\perp 
\end{bmatrix}\right\}
%\begin{bmatrix}
%\mathbf{p}_s & \mathbf{0} & \mathbf{0} \\
%\mathbf{p}_s^\perp & \mathbf{0} & \vdots\\
%\mathbf{0} & \ddots &  \vdots \\
%\vdots & \ddots & (-1)^{K+1}\mathbf{p}_s\\
%\textbf{0} & \ldots & (-1)^{K+1}\mathbf{p}_s^\perp  
%\end{bmatrix}
%\bigoplus_{m=1}^K (-1)^{m+1} \begin{bmatrix} \mathbf{p}_s \\  \mathbf{p}_s^\perp \end{bmatrix} ,\]
,\]
the diagonal matrix
\[ 
\bfD_{2,K+1} \deff 
\bfI_K \kron \left\{ \rho^{-(\ell+1)} \bfI_2\right\}
%\begin{bmatrix}
%\rho^{-2} & 0 & \ldots & \ldots & 0\\
%0 & \rho^{-2} & 0 & \ldots & 0\\
%\vdots & \ddots & \ddots & \ddots & \vdots \\
%0 & \ldots & 0 & \rho^{-K-1} & 0\\
%0 & \ldots & \ldots & 0 & \rho^{-K-1}
%\end{bmatrix} 
,\]
%\bfI_K \kron \{ \rho^{-i} \bfI_2 \}_{i=2,...,K+1} \] %\bigoplus_{m=1}^K \rho^{-(m+1)} \bfI_2 ,\]
%\[
%\bfD_{\kappa,\kappa'} = \bfI_{\kappa'-\kappa + 1} \kron \{ \rho^{-i} \bfI_2 \}_{i=\kappa,...,\kappa'},\]

and the row vector
\[ \bfZ_{s,:} \deff \begin{bmatrix}
\cos(2\overline{\theta}_s) & \sin(2\overline{\theta}_s) & \cos(3\overline{\theta}_s) & \sin(3\overline{\theta}_s) & \dots & \cos((K+1)\overline{\theta}_s) & \sin((K+1)\overline{\theta}_s)
\end{bmatrix}  . \]
Then 
\[(\bfS_{dip})_{s,\,:} = \bfZ_{s,\,:\,} \bfP^{(s)\top}_K \bfD_{2,K+1} .\]

%\[ (\bfS_{SL})_{s,\,:} =  \frac{1}{2 \pi} \left [ \int_{\p \Omega_s} \frac{\p u_s}{\p \nu} \biggr |_{+}  \frac{\cos(\theta_y)}{r_y} \text{ d} \sigma_y  \quad  \int_{\p \Omega_s} \frac{\p u_s}{\p \nu} \biggr |_{+}  \frac{\sin(\theta_y)}{r_y} \text{ d} \sigma_y  \right . \]
%\[\int_{\p \Omega_s} \frac{\p u_s}{\p \nu} \biggr |_{+}  \frac{\cos(2\theta_y)}{2 r_y^2} \text{ d} \sigma_y   
%\quad \int_{\p \Omega_s} \frac{\p u_s}{\p \nu} \biggr |_{+}  \frac{\sin(2\theta_y)}{2 r_y^2} \text{ d} \sigma_y  \dots \]
%\[ \left . \dots \int_{\p \Omega_s} \frac{\p u_s}{\p \nu} \biggr |_{+}  \frac{\cos(K\theta_y)}{K r_y^K} \text{ d} \sigma_y   
%\quad \int_{\p \Omega_s} \frac{\p u_s}{\p \nu} \biggr |_{+}  \frac{\sin(K\theta_y)}{K r_y^K} \text{ d} \sigma_y  \right ] \]

On the other hand, given $N$ points $y_i$ uniformly distributed on $\p \Omega_s$, we can discretize the integral defining $(\bfS_{SL})_{s,2k-1}$ and $(\bfS_{SL})_{s,2k}$ as follows
\[ A_{s,k} =   \int_{\p \Omega_s} \frac{\p u^{(s)}}{\p \nu} \biggr |_{+}  \frac{\cos(k\theta_y)}{k r_y^k} \text{ d} \sigma_y \approx \sum_{i=1}^{N}  \frac{\p u^{(s)}}{\p \nu}(y_i) \frac{\cos(k\theta_{y_i})}{k r_i^k}  \Delta y_i,  \]
\[ B_{s,k} =   \int_{\p \Omega_s} \frac{\p u^{(s)}}{\p \nu} \biggr |_{+}  \frac{\sin(k\theta_y)}{k r_y^k} \text{ d} \sigma_y \approx \sum_{i=1}^{N}  \frac{\p u^{(s)}}{\p \nu}(y_i) \frac{\sin(k\theta_{y_i})}{k r_i^k}  \Delta y_i.  \]
Consequently, defining the column vector
\[ \bfU_{s,\,:} \deff   \begin{pmatrix} \displaystyle \frac{\p u^{(s)}}{\p \nu}(y_i)  \Delta y_i \end{pmatrix}_{i = 1}^N \,,\]
we get 
\[(\bfS_{SL})_{s,\,:} = - \bfU_{s,\,:\,} \bfG_*^{(s)}.
\] 
Notice that $\bfG_*^{(s)}$ reduces to $\bfG^{(s)}$ if we choose the receptors as discretization points, i.e., $y_i = x_i^{(s)}$.  %\bfU_{s,:} \bfV_K \widetilde{\bfJ}_K \bfD_K$.

In the end, the source vector can be written in the following form:
\[  \bfS_{s,\,:} = (\bfS_{dip})_{s,\,:} + (\bfS_{SL})_{s,\,:}  = \bfZ_{s, \, :} \bfP^{(s)\top}_K \bfD_{2,\,K+1} - \bfU_{s,\,:} \bfG^{(s)}.  \]

Let $U$ be the background solution, i.e., the potential in the absence of any target. When $\rho \gg \delta$, the dipolar expansion derived in \cite{AMM2013} yields \begin{equation}\bfU_{s,\,:} = \mathbf{u}  + O(\delta^2), \label{eq:dipole_approx}\end{equation}
where $\mathbf{u}  \deff  \begin{pmatrix} \displaystyle \frac{\p U}{\p \nu}(y_i)  \Delta y_i \end{pmatrix}_{i = 1}^{N_r}$. Such first order approximation of $\bfU_{s,\,:}$ depends only on the geometry of the fish and on the position of the receptors.
%\begin{align*} \bfS_{s,\,:} \kron \bfG^{(s)} & = \bfZ_{s,:} \bfP^\top \bfD\kron  \bfV_K \widetilde{\bfJ}_K \bfD_K   - \mathbf{u} \bfV_K \widetilde{\bfJ}_K \bfD_K \kron \bfV_K \widetilde{\bfJ}_K \bfD_K    \\ & =  \bfZ_{s,:} \bfP^\top \overline{\bfD} \kron  \bfV \widetilde{\bfJ} \bfD   - \mathbf{u} \bfV \widetilde{\bfJ} \bfD \kron \bfV \widetilde{\bfJ} \bfD   \\ & =  \bfZ_{s,:} \bfP^\top \overline{\bfD} \kron  \bfV \widetilde{\bfJ} \bfD   - (\mathbf{u} \kron \bfI_{N_r}) ( \bfV \widetilde{\bfJ} \bfD \kron \bfV \widetilde{\bfJ} \bfD ) \\ & =   \bfZ_{s,:} \bfP^\top \overline{\bfD} \kron  \bfV \widetilde{\bfJ} \bfD   - (\mathbf{u} \kron \bfI_{N_r}) (\bfV \widetilde{\bfJ} \bfD )^{\kron 2}. \end{align*}

%\[ \begin{bmatrix} z_1 & \overline{z}_1 & z_1^2 & \overline{z}_1^2 &   \dots & z_1^{K+1} & \overline{z}_1^{K+1} \\
%z_2 & \overline{z}_2 & z_2^2 & \overline{z}_2^2 &   \dots & z_2^{K+1} & \overline{z}_2^{K+1} 
%\\
%\vdots & & \vdots   & &   \ddots &  & \vdots \\
%z_{M} & \overline{z}_M & z_M^2 & \overline{z}_M^2 &   \dots & z_M^{K+1} & \overline{z}_M^{K+1} 
%\end{bmatrix} \]

\subsection{Reconstruction error analysis}

In this section we analyze the relative error in the reconstruction of the CGPTs when the g-inverse $\mathfrak{L}$ given by \eqref{eq:minnormsol} is used, and the MSR data are acquired along a single circular orbit around the target.

Let $\bfW$ be a random matrix $N_r \times M$ with independent and identically distributed $\mc{N}(0,\sigma^2_{\text{noise}})$ entries. Let $\mathbf{E}$ be the matrix $N_r \times M$ of the truncation errors. Recall that the entries of $\bfE$ are of order $\varepsilon^{K+2}$, see \eqref{eq:errorEsr}.

The following multivariate multiple linear regression model for the measurements can be stated:
\begin{equation} \label{eq:measurement_model} 
\mathbf{Q}  = \mc{L}(\mathbb{M}) + \mathbf{E} + \mathbf{W}.
\end{equation}

We restrict ourselves to the situation where the strength of the noise is enough to overpower the truncation error, which we disregard a posteriori for the rest of the analysis. More precisely, we assume that the strength of the noise satisfies
\begin{equation}
\varepsilon^{K+2} \ll \sigma^2_{\text{noise}} \ll \varepsilon^2 .
\end{equation}
We define the signal to noise ratio (SNR) associated to the orbit $\mc{O}$ as
\[ \text{SNR} = \frac{\varepsilon^2}{\sigma^2_{\text{noise}}} .\]

Next, we vectorize equation \eqref{eq:regression_eq}. Define the vectorized error matrix
\begin{equation*} \label{eq:vec_E} 
\vec{\bfW} \sim \gaussd (\mathbf{0}_{M N_r} , \sigma^2_{\text{noise}} \bfI_{M N_r} ),
\end{equation*} 
and response matrix, which is a multivariate normal vector given $\bfS$ and $\bfG^{(s)}$ for $s = 1, ... , M$, namely,
\begin{equation*} \label{eq:vec_response} 
\vec{\bfQ\,|\,\{\bfG^{(s)}\}, \bfS} \sim \gaussd (\vec{\mc{L}(\mathbb{M})} , \sigma^2_{\text{noise}} \bfI_{M N_r} ).
\end{equation*}

Straightforward computations show that the covariance matrix of $\mathfrak{L}\, \vec{\bfW}$ can be written as
\begin{equation} \label{eq:cov} \begin{split} 
\mbox{Cov} ( \mathfrak{L} \,\vec{\bfW} ) &= \sigma_{\text{noise}}^2 \sum_{s = 1}^M \;\bfS^\dagger_{:,s} \bfS^{\dagger \top}_{:,s} \kron  \bfG^{(s)\dagger} \bfG^{(s) \dagger \top} \\ & = \sigma_{\text{noise}}^2  \sum_{s=1}^M \;\begin{bmatrix} (\bfS^\dagger_{1,s})^2   \bfG^{(s)\dagger} \bfG^{(s) \dagger \top} & & * \\ & \ddots & \\  * & & (\bfS^\dagger_{2K,s})^2   \bfG^{(s)\dagger} \bfG^{(s) \dagger \top}  \end{bmatrix}.
\end{split} 
\end{equation}
%Since we have proved that we can approximate $\bfS^\dagger$ with $\bfS_{SL}^\dagger$, we are ready to prove the following theorem.

Hereinafter, we assume that the $N_r$ receptors of $\bfG^{(s)}$ are all distributed along one circular arc of radius $\rho$. This assumption is justified by the fact that if we model the fish skin by two close circular arcs of radius $\rho$ and $\rho + \varepsilon$, with $N_r/2$ receptors on them, and we call $\bfG^{(s)}_\varepsilon$ its matrix of receptors, we can show that $\|\bfG^{(s)}_\varepsilon - \bfG^{(s)}\|_F \to 0$ as $\varepsilon \to 0$.  

\begin{theorem} \label{thm:error_est}
	For $j,k$ so that $(\bfM)_{jk}$ is non-zero, the relative error on the reconstructed CGPT satisfies
	\begin{equation} 
	\sqrt{\frac{\mathbb{E}|(\bfM^{\text{est}})_{jk} - (\bfM)_{jk}|^2}{|(\bfM)_{jk}|^2}} \lesssim \sigma_{\text{noise}} \; \varepsilon^{-\left\lceil j/2 \right\rceil - \left\lceil k/2 \right\rceil } \left\lceil \frac{j}{2} \right\rceil \left\lceil \frac{k}{2} \right\rceil.  
	\label{rel_err} 
	\end{equation} 
	For vanishing $(\bfM)_{jk}$, the error $\sqrt{\mathbb{E}|(\bfM^{\text{est}})_{jk} - (\bfM)_{jk}|^2}$ can be bounded by the right-hand side above with $\varepsilon$ replaced by $\rho^{-1}$.  
\end{theorem}
\begin{proof} We begin by observing that the absolute errors $\sqrt{\mathbb{E}|(\bfM^{\text{est}})_{jk} - (\bfM)_{jk}|^2}$ are the diagonal entries of Cov$(\mathfrak{L} \bfW)$. In particular, the $(j,j)$-th entry of the $k$-th block matrix of Cov$(\mathfrak{L}\, \vec{\bfW})$ given by \eqref{eq:cov}, i.e., $\sum_{s=1}^M(\bfS^\dagger_{k,s})^2 (\bfG^{(s)\dagger} \bfG^{(s) \dagger \top})_{jj}$, corresponds to CGPT $(\bfM)_{jk}$.
	Define $\mathcal{I}_{jk}:=(\bfM^{\text{est}} - \bfM)_{jk}$. By Lemma \ref{lem:estimate_mp_S} we have the inequality
	\[
	|\bfS^\dagger_{k,s}|^2 \leq \|\bfS^\dagger_{k,\,:}\|^2_F \lesssim \frac{\rho^{2 \left\lceil k/2 \right\rceil}}{M}\left\lceil \frac{k}{2} \right\rceil^2.
	\]
	On the other hand it is easy to show that
	\[ | (\bfG^{(s)\dagger} \bfG^{(s) \dagger \top})_{jj} | \lesssim \rho^{2 \left\lceil j/2 \right\rceil}\left\lceil \frac{j}{2} \right\rceil^2. \]

	%	\[
	%	\bfG^{(s)\dagger} \bfG^{(s) \dagger \top}  = \bfC^\dagger  \bfD_{1,K}^\dagger (\widetilde{\bfG}^{(1)}  (\bfI_K \kron \{r_1(i(s-1) \gamma)\}_{i}))^\dagger (\widetilde{\bfG}^{(1)}  (\bfI_K \kron \{r_1(i(s-1) \gamma)\}_{i}))^{\dagger \top} \bfC^\dagger  \bfD_{1,K}^\dagger.
	%	\]
	%	Since
	%	\[
	%	(\bfG^{(s)\dagger} \bfG^{(s) \dagger \top})_{jj} = \rho^{\left\lceil j/2 \right\rceil}\left\lceil \frac{j}{2} \right\rceil \left((\widetilde{\bfG}^{(1)}  (\bfI_K \kron \{r_1(i (s-1) \gamma)\}_{i}))^\dagger (\widetilde{\bfG}^{(1)}  (\bfI_K \kron \{r_1(i(s-1) \gamma)\}_{i}))^{\dagger \top}\right)_{jj} \rho^{\left\lceil j/2 \right\rceil} \left\lceil \frac{j}{2} \right\rceil,
	%	\]
	Therefore, we obtain the following estimate
	%	\[
	%	\begin{aligned}
	%	\mathbb{E}(\mathcal{I}_{jk})^2 = \sum_{s=1}^M (\bfS^\dagger_{k,s})^2   (\bfG^{(s)\dagger} \bfG^{(s) \dagger \top})_{jj} \lesssim \rho^{2 (\left\lceil j/2 \right\rceil + \left\lceil k/2 \right\rceil ) } \left\lceil \frac{j}{2} \right\rceil^2\left\lceil \frac{k}{2} \right\rceil^2 .
	%	\end{aligned}
	%	\]
	\[
	\begin{aligned}
	\mathbb{E}(\mathcal{I}_{jk})^2 = \sum_{s=1}^M (\bfS^\dagger_{k,s})^2   (\bfG^{(s)\dagger} \bfG^{(s) \dagger \top})_{jj} \lesssim \rho^{2 (\left\lceil j/2 \right\rceil + \left\lceil k/2 \right\rceil ) } \left\lceil \frac{j}{2} \right\rceil^2\left\lceil \frac{k}{2} \right\rceil^2 .
	\end{aligned}
	\]
	The scaling property
	\[ 
	(\bfM)_{jk}(\delta B) = \delta^{\left\lceil j/2 \right\rceil + \left\lceil k/2 \right\rceil}(\bfM)_{jk} (B),
	\] 
	together with the above control on $\mathbb{E}(\mathcal{I}_{jk})^2$ show that   the relative error satisfies \eqref{rel_err}.
	
\end{proof}

\begin{rem}	In the proof we used the following inequality:
	\[ | (\bfG^{(s)\dagger} \bfG^{(s) \dagger \top})_{jj} | \lesssim \rho^{2 \left\lceil j/2 \right\rceil}\left\lceil \frac{j}{2} \right\rceil^2 , \]
	where the unspecified constant depends on the number of receptors $N_r$ and the angle of view $\gamma$. In the next section we analyze such dependency in the limit as $N_r \to \infty$,  observing that the upper bound in \eqref{rel_err} decays like $1/N_r$.
\end{rem}

Following \cite{Ammari2014}, given the SNR and a tolerance level $\tau_0$, the resolving order is defined as

\begin{equation*} K^* = \min \left \{ 1 \le k \le K : \sqrt{\frac{\mathbb{E}|(\bfM^{\text{est}})_{kk} - (\bfM)_{kk}|^2}{|(\bfM)_{kk}|^2}} \le \tau_0 \right \}. \end{equation*}
It is readily seen that the resolving order $K^*$ satisfies
\begin{equation} \label{eq:res_order} (K^* \varepsilon^{1-K^*})^2 \simeq \tau_0 \text{SNR}. \end{equation}

\subsection{Angular resolution}

In this section we want to discuss the issues related with the limited-view configuration, which is an intrinsic feature of the fish geometry. In order to get a grasp of how much the angle of view has an impact on the error estimate provided by Theorem \ref{thm:error_est}, we restrict ourselves to a special configuration of receptors. In particular, we assume that there are $N_r$ receptors evenly distributed on an arc of the unit circle, with aperture angle $\gamma \in (0, 2 \pi)$, and we let $N_r$ go to $\infty$.

Instead of studying the spectrum of  $(\bfG^{(s) \mathsf{H}} \bfG^{(s)})^\dagger$  we shall equivalently consider that of  $(\bfV_K^{\mathsf{H}} \bfV_K)^{\dagger}$, where $\bfV_K$ is defined in Appendix \ref{apx:uniqueness}, and $\bfV_K^\mathsf{H}$ denotes the matrix $\overline{\bfV}_K^\top$.
For the sake of notation we refer to $\bfV_{K}$ as the block matrix $\begin{bmatrix}
\bfW_K & \overline{\bfW}_K 
\end{bmatrix},$ obtained from $\bfV_{K}$ by permuting the columns as in Appendix \ref{apx:uniqueness}. We are interested in the asymptotic expansion of $\bfV_K^\mathsf{H} \bfV_K$ as $N_r \to \infty$. Hereinafter, we denote $\lim\limits_{N_r\to\infty} \frac{1}{N_r} \textbf{A}$  by $(\textbf{A})_\infty$. 

With this particular geometry of receptors the limit matrix can be analytically computed:
\[
(\bfV_K^\mathsf{H} \bfV_K)_\infty = 
\begin{bmatrix} 
(\overline{\bfW}_K^\top \bfW_K)_\infty & (\overline{\bfW}_K^\top \overline{\bfW}_K)_\infty \\ \\ 
(\bfW_K^\top \bfW_K)_\infty & (\bfW_K^\top \overline{\bfW}_K)_\infty
\end{bmatrix},
\]
where
\[
(\bfW_K^\top \bfW_K)_\infty =
\begin{bmatrix}
-\frac{1-e^{i2\gamma}}{2i\gamma}& -\frac{1-e^{i3\gamma}}{3i\gamma}& -\frac{1-e^{i4\gamma}}{4i\gamma} & \ldots & -\frac{1-e^{i(K+1)\gamma}}{(K+1)i\gamma} \\ \\
-\frac{1-e^{i3\gamma}}{3i\gamma} & -\frac{1-e^{i4\gamma}}{4i\gamma} & -\frac{1-e^{i5\gamma}}{5i\gamma} & \ldots & -\frac{1-e^{i(K+2)\gamma}}{(K+2)i\gamma} \\\\
\vdots & \vdots & \vdots & &  \vdots \\\\
-\frac{1-e^{i(K+1)\gamma}}{(K+1)i\gamma} & -\frac{1-e^{i(K+2)\gamma}}{(K+2)i\gamma} & -\frac{1-e^{i(K+3)\gamma}}{(K+3)i\gamma} & \ldots & -\frac{1-e^{2Ki\gamma}}{2Ki\gamma}
\end{bmatrix},
\]
%\[
%(\overline{\bfW}_K^\top \overline{\bfW}_K )_\infty =
%\begin{bmatrix}
%\frac{1-e^{-i2\gamma}}{2i\gamma}& \frac{1-e^{-i3\gamma}}{3i\gamma}& \frac{1-e^{-i4\gamma}}{4i\gamma} & \ldots & \frac{1-e^{-i(K+1)\gamma}}{(K+1)i\gamma} \\ \\
%\frac{1-e^{-i3\gamma}}{3i\gamma} & \frac{1-e^{-i4\gamma}}{4i\gamma} & \frac{1-e^{-i5\gamma}}{5i\gamma} & \ldots & \frac{1-e^{-i(K+2)\gamma}}{(K+2)i\gamma} \\\\
%\vdots & \vdots & \vdots & &  \vdots \\\\
%\frac{1-e^{-i(K+1)\gamma}}{(K+1)i\gamma} & \frac{1-e^{-i(K+2)\gamma}}{(K+2)i\gamma} & \frac{1-e^{-i(K+3)\gamma}}{(K+3)i\gamma} & \ldots & \frac{1-e^{-2Ki\gamma}}{2Ki\gamma}
%\end{bmatrix},
%\]
\[
(\overline{\bfW}_K^\top {\bfW}_K)_\infty = \begin{bmatrix}
1 & -\frac{1-e^{i\gamma}}{i\gamma}& -\frac{1-e^{i2\gamma}}{2i\gamma} & \ldots & -\frac{1-e^{i(K-1)\gamma}}{(K-1)i\gamma} \\\\
\frac{1-e^{-i\gamma}}{i\gamma} & 1 & -\frac{1-e^{i\gamma}}{i\gamma} & \ldots & -\frac{1-e^{i(K-2)\gamma}}{(K-2)i\gamma} \\ \\
\frac{1-e^{-2i\gamma}}{2i\gamma} & \frac{1-e^{-i\gamma}}{i\gamma} & 1 & \ldots & -\frac{1-e^{i(K-3)\gamma}}{(K-3)i\gamma}\\ \\
\vdots & \vdots & \vdots &  &  \vdots \\ \\
\frac{1-e^{-i(K-1)\gamma}}{i(K-1)\gamma} & \frac{1-e^{-i(K-2)\gamma}}{i(K-2)\gamma} & \ldots & \frac{1-e^{-i\gamma}}{i\gamma} & 1
\end{bmatrix},
\]
%\[
%({\bfW}_K^\top \overline{\bfW}_K)_\infty = 
%\begin{bmatrix}
%1 & \frac{1-e^{-i\gamma}}{i\gamma}& \frac{1-e^{-i2\gamma}}{2i\gamma} & \ldots & \frac{1-e^{-i(K-1)\gamma}}{(K-1)i\gamma} \\\\
%-\frac{1-e^{i\gamma}}{i\gamma} & 1 & \frac{1-e^{-i\gamma}}{i\gamma} & \ldots & \frac{1-e^{-i(K-2)\gamma}}{(K-2)i\gamma} \\ \\
%-\frac{1-e^{2i\gamma}}{2i\gamma} & -\frac{1-e^{i\gamma}}{i\gamma} & 1 & \ldots & \frac{1-e^{-i(K-3)\gamma}}{(K-3)i\gamma} \\ \\
%\vdots & \vdots & \vdots &  &  \vdots \\ \\
%-\frac{1-e^{i(K-1)\gamma}}{i(K-1)\gamma} & -\frac{1-e^{i(K-2)\gamma}}{i(K-2)\gamma} & \ldots & -\frac{1-e^{i\gamma}}{i\gamma} & 1
%\end{bmatrix}.
%\]
and $({\bfW}_K^\top \overline{\bfW}_K)_\infty = \overline{(\overline{\bfW}_K^\top {\bfW}_K)_\infty}$, $(\overline{\bfW}_K^\top \overline{\bfW}_K)_\infty = \overline{({\bfW}_K^\top {\bfW}_K)_\infty}$.

Hence,
\[
\bfV_K^\mathsf{H} \bfV_K =   N_r  \left ( (\bfV_K^\mathsf{H} \bfV_K)_\infty + 
O \left (\frac{1}{N_r} \right ) \right ), \mbox{ as } N_r \to \infty.
\]

Applying the results contained in \cite{deif} we obtain the following result.

\begin{lemma} \label{lem:singular-val}
	For $N_r$ large,
	\begin{equation}\label{eq:kernel-expansion} 
	(\bfV_K^\mathsf{H} \bfV_K)^{-1} = \frac{1}{N_r} \left((\bfV_K^\mathsf{H} \bfV_K)_\infty + 
	O \left ( \frac{1}{N_r}\right ) \right)^{-1},
	\end{equation}
	where the entries of $(\bfV_K^\mathsf{H} \bfV_K)_\infty$ depend only on the angle of view $\gamma$ and the truncation order $K$. Moreover, we have
	\begin{equation}\label{eq:singular-val} 
	\frac{1}{\widetilde{\sigma}_k} = \frac{1}{N_r} \left (\sigma_k + O \left ( \frac{1}{N_r}\right ) \right)^{-1},
	\end{equation} 
	where $\widetilde{\sigma}_k$ (resp. $\sigma_k$) are the eigenvalues of the matrix $\bfV_K^\mathsf{H} \bfV_K$ (resp. $(\bfV_K^\mathsf{H} \bfV_K)_\infty$).
\end{lemma}

\begin{rem} Lemma \ref{lem:singular-val} highlights the following facts. On one hand, as far as all the eigenvalues are away from zero, the reconstructed CGPTs have an upper bound on the relative error in \eqref{rel_err}  which decays like $1/N_r$. 
	%This follows from the fact that if $(\bfV_K^\top \bfV_K)_\infty$ is invertible, then we are able to expand \eqref{eq:kernel-expansion} by \cite{ONERRORBOUND}, with its first order term decaying exactly like $1/N_r$. 
	More precisely, $(\bfG^{(s)\dagger} \bfG^{(s) \dagger \mathsf{H}})_{jj} \lesssim 1/N_r \sum 1/|\sigma_k|$. 
	On the other hand, when some eigenvalues occur to be very small (e.g. $\sigma_l \le 10^{-8}$), inequality \eqref{rel_err} becomes uninformative, making us unable to predict the behavior of the relative error. 
\end{rem}

%The previous remark suggests that there is an angular resolving order $K_{\#} =  K_{\#}(\gamma, K)$, depending on the angle-view, which is the maximum order we can resolve given the 

\begin{figure}[H]
	\centering
	\hspace*{-5mm}
	\begin{tabular}{cc}
		\includegraphics[scale=0.5]{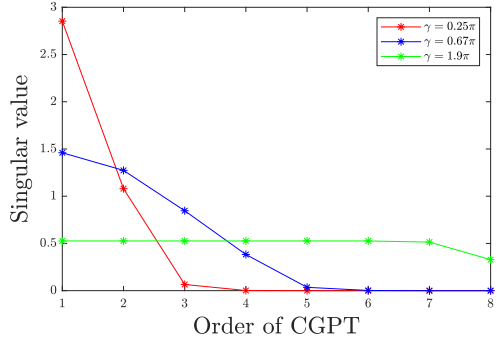} & \raisebox{33ex - \height}{\includegraphics[scale=0.5]{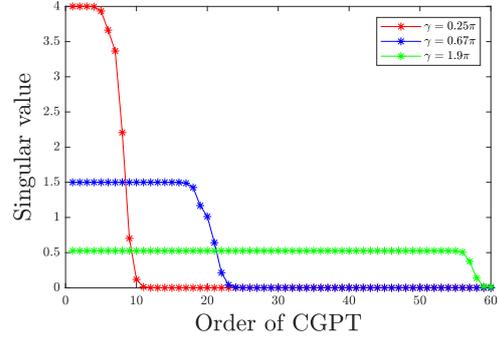}} \\
		\raisebox{2.5ex - \height}{(a) $K=4$} & \raisebox{2.5ex - \height}{(b) $K=30$} \\
	\end{tabular}
	\caption{The distribution of eigenvalues at different angles of view $\gamma$. The solid line corresponds to $\sigma(\bfV_K^\mathsf{H} \bfV_K / N_r)$, the starred line to $\sigma((\bfV_K^\mathsf{H} \bfV_K)_\infty)$. We use $N_r = 10^4$ receptors.}
	\label{fig:eig}
\end{figure}

Figure \ref{fig:eig} provides the distribution of $\widetilde{\sigma}_k$, ${\sigma}_k/N_r$ at different values of the reconstruction order $K$ and angles of view $\gamma$, as $N_r \to \infty$.
Firstly, we clearly observe the asymptotic behavior of the spectrum stated by \eqref{eq:singular-val}. Secondly, we notice that the effect of the limited-view configuration is reflected by the decaying of the eigenvalues of the matrix of receptors. As expected, the closer the angle of  view is to $2\pi$, the more informative estimate \eqref{rel_err} becomes. Furthermore, if $K$ is small, the angle of view has a minor impact on the reconstructed CGPTs.

%	we are able to lead a numerical study to illustrate the dependence of our reconstruction error on $\gamma$ and $K$ by studying the eigenvalues of $(\bfV_K^\top \bfV_K)_\infty$, where the dependence on $N_r$ is neglected. Figures X1 and X2 show the distribution of eigenvalues $\widetilde{\sigma}_k$ (resp. ${\sigma}_k$) of the matrix $\bfV_K^\top \bfV_K$ (resp. $(\bfV_K^\top \bfV_K)_\infty$) for different orders $K$ and large $N_r$. Each eigenvalue is related to the reconstruction error affecting the correspondent GPT. Observe that the largest eigenvalues correspond to the first GPTs. From these results, we see clearly the effect of the limited-view aspect [CITE TRACKING]. One can see also that when the angle view is closer to $2\pi$, we are able to reconstruct correctly a larger number of GPTs. Lemma \ref{lem:singular-val} shows that for the eigenvalues of $\bfV_K^\top \bfV_K$ corresponding to the first GPTs, i.e., $\widetilde{\sigma}_k$ is not too close to zero, we can make the dependence on $N_r$ explicit. Their formula is given by \eqref{eq:singular-val}. For $N_r$ large enough and $K<N_r$, this means that the information on the dependence of the eigenvalues $\widetilde{\sigma}_k$ on $\gamma$ and $K$ is all encoded in the eigenvalues of $(\bfV_K^\top \bfV_K)_\infty$ when these are not too close to zero. Figure X1 and X2 show that such eigenvalues are the majority of the spectrum of $\bfV_K^\top \bfV_K$. 
%
%
%
%Figures... clearly show that the estimate \eqref{rel_err} is uninformative for some higher-orders CGPTs.

\section{Recognition}
\label{sec:recognition}

In the previous section, we established an upper bound on the reconstruction error which essentially depends on the length-scale of the acquisition orbit. As a consequence, the closer the fish gets to the target, the higher is the order of the features it can retrieve from the noisy measurements.
This result suggests that, when it comes to classification, it is of preeminent importance to design recognition algorithms that exploit the information contained in measurements collected at multiple scales.

This section aims at presenting a novel multi-scale algorithm for target classification.

\subsection{Complex CGPTs and shape descriptors}

Let us briefly recall the definition of the shape descriptors.

We introduce convenient complex combinations of CGPTs. For any pair of indices $m, n = 1 , 2 , ... $ ,  we introduce the following quantities
\begin{equation*} \label{eq:CCGPTs} \begin{matrix} \mathbf{N}_{mn}^{(1)} (\gamma , D) = (M_{mn}^{cc} - M_{mn}^{ss}) + i  (M_{mn}^{cs} + M_{mn}^{sc}),  \\ \mathbf{N}_{mn}^{(2)} (\gamma , D) = (M_{mn}^{cc} + M_{mn}^{ss}) + i  (M_{mn}^{cs} - M_{mn}^{sc}) . \end{matrix}  \end{equation*}
We define the complex CGPT matrices by
\[ \mathbf{N}^{(1)}: = ( \mathbf{N}^{(1)}_{mn} )_{m,n}, \qquad \mathbf{N}^{(2)}: = ( \mathbf{N}^{(2)}_{mn} )_{m,n} .  \]

We call a dictionary $\mathcal{D}$ a collection of standard shapes $B$, centered at the origin, with characteristic size of order $1$.
We assume that a reference dictionary $\mathcal{D}$ is initially given. Furthermore, suppose to consider a shape $D$, which is unknown, that is obtained from an element $B \in \mathcal{D}$ by applying some unknown rotation $\theta$, scaling $s$ and translation $z$, i.e., $D = T_z s R_\theta \,B$.

Following \cite{Ammari2014}, let $\eta = \frac{\mathbf{N}_{12}^{(2)} (D)}{2 \mathbf{N}_{11}^{(2)} (D)}$.
We define the following quantities
\begin{equation*} \mathcal{J}^{(1)} (D) = \mathbf{N}^{(1)} ( T_{-\eta} D ) = \mathbf{C}^{-\eta} \mathbf{N}^{(1)} (D) ( \mathbf{C}^{-\eta})^T , \end{equation*}
\begin{equation*} \mathcal{J}^{(2)} (D) = \mathbf{N}^{(2)} ( T_{-\eta} D ) = \overline{\mathbf{C}^{-\eta}} \mathbf{N}^{(2)} (D) ( \mathbf{C}^{-\eta})^T , \end{equation*}
where the matrix $\mathbf{C}^{-\eta}$ is a lower triangular matrix with the $m,n$-th entry given by
\[ \mathbf{C}^{-\eta}_{mn} = \binom{m}{n} \, (-\eta)^{m-n} .\]

These quantities are translation invariant.

\medskip

From $\mathcal{J}^{(1)}(D) = (\mathcal{J}_{mn}^{(1)}(D) )_{m,n}$, $\mathcal{J}^{(2)} (D)  = (\mathcal{J}_{mn}^{(2)} (D) )_{m,n}$, for each pair of indices $m,n$, we define  the scaling invariant quantities:
\begin{equation*} \mathcal{S}^{(1)}_{mn}(D)  = \frac{ \mathcal{J}_{mn}^{(1)}(D)}{(\mathcal{J}_{mm}^{(2)}(D)  \mathcal{J}_{nn}^{(2)}(D) )^{1/2}} , \qquad \mathcal{S}^{(2)}_{mn} (D) = \frac{ \mathcal{J}_{mn}^{(2)}(D) }{(\mathcal{J}_{mm}^{(2)}(D)  \mathcal{J}_{nn}^{(2)}(D) )^{1/2}} . \end{equation*}

The CGPT-based shape descriptors $\mathcal{I}^{(1)} = ( \mathcal{I}^{(1)}_{mn})_{m,n}$ and $\mathcal{I}^{(2)} = ( \mathcal{I}^{(2)}_{mn})_{m,n}$ are defined as:
\[\mathcal{I}^{(1)}_{mn} = | \mathcal{S}^{(1)}_{mn} (\gamma,D)  | , \qquad \mathcal{I}^{(2)}_{mn} = | \mathcal{S}^{(2)}_{mn} (\gamma,D)  | , \]
where $| \, \cdot \, |$ denotes the modulus of a complex number. Recall that $\mathcal{I}^{(1)}$ and $\mathcal{I}^{(2)}$ are invariant under translation, rotation, and scaling.

The details of this construction can be found in \cite{Ammari2014}.

\subsection{Multi-scale acquisition setting} \label{sbsec:setting_orbits}

Suppose that the scanning movement consists of $\mc{M}$ concentric circular orbits 
$\mc{O}_1 , \mc{O}_2 ,  ... , \mc{O}_{\mc{M}}$, with radii $\rho_1 > \rho_2 > ... > \rho_{\mc{M}}$ respectively (ordered from the farthest to the nearest), the target being located at the common center. On each orbit only a discrete number of positions accounts for the data acquisition process, as described in Section  \ref{sbsec:setting}. Precisely, $M_j$ different positions are sampled along the orbit $\mc{O}_j$, and for each position $s$ the corresponding electric signal $u_{j}^{(s)} - H_{j}^{(s)}$ is measured by $N_r$ receptors on the skin, $\{ x_{j,r}^{(s)} \}_{r = 1}^{N_r}$.

%Precisely, along orbit $\p B_{r_j}$ $M_j$ different positions are sampled, and the signal is statically recorded at the $N_r$ receptors $x_r^{(s)}$ corresponding to the position $s$. 

%\begin{center}
%	\begin{tabular}{@{} *5l @{}}    \toprule
%		\emph{Symbol} & \emph{Meaning}   \\\midrule
%		$\Omega_{s}$   & Fish body \\ 
%		$\mathbf{p}_{s}$ & dipole moment \\ 
%		$\zeta_{s}$ & electric organ \\ 
%		$x_r^{(s)}$  & $r$-th receptor \\$u_{j,s}$  & electric potential solution to \eqref{eq:model_u} \\
%		$H_{j,s}$  & $\dots$ \\\midrule
%		$N_r$  & number of receptors \\
%		$M$  & number of positions \\
%		$\mc{O}$ & orbit\\\bottomrule
%		\hline
%	\end{tabular}
%	\captionof{table}{Notation referred to position $s \in \{ 1 , ... , M \}$ on the orbit $\mc{O}$. \label{Tab:Tcr2}}
%\end{center}

Therefore, for any orbit $\mc{O}_j$ we get the $N_r \times M_j$ MSR matrix, whose $(r,s)$-entry is defined as
\begin{equation}(\bfQ_{\mc{O}_j})_{r,s} = u_{j}^{(s)}(x_{j,r}^{(s)}) - H_{j}^{(s)}(x_{j,r}^{(s)}), \qquad j \in \{ 1 , ... , \mc{M} \} . \label{eq:MSR-multiscale} \end{equation}

Notice that so far the setting described above is very general, as no restriction has been given on the radii of the orbits yet.

%It will become clear soon that not every configuration is admissible in some sense.

Since the orbits are at different length-scales, it is clear that the closer the orbit is to the center, the stronger the MSR signal is. However, resorting to the error estimate on the reconstruction order \eqref{rel_err}, we are able to choose the orbits in such a way that the resolving order is enhanced.

In the multi-scale setting described above, formula \eqref{eq:res_order} reads
\begin{equation} \label{eq:res_order_multi} (K_{j}^* \varepsilon_j^{1-K^*_j})^2 \simeq \tau_0 \text{SNR}_j, \end{equation}
where $\text{SNR}_j = \frac{\varepsilon^2_j}{\sigma^2_{\text{noise}}}$ is the signal-to-noise ratio associated to the $j$-th orbit.

Thus, a length-scale dependent resolving order is introduced, and obviously $K_{j+1}^* \ge K_{j}^* \ge 2$.

The noisy MSR matrix is given by the following formula
\begin{equation} \label{eq:regression-multiscale}
\bfQ_{\mc{O}_j} = \mathcal{L}_{j}(\bbM^{(K_j)}) + \mathbf{E}_{\mc{O}_j} + \mathbf{W} .
\end{equation}

On each orbit $\mc{O}_j$, the CGPTs can be retrieved from the data \eqref{eq:regression-multiscale}, for instance, by either using the classical Moore-Penrose inverse or the generalized inverse $\mathfrak{L}$ given by \eqref{eq:minnormsol}.
Moreover, we denote by $(\mathcal{I}^{(1)}(D; \mc{O}_j) ,  \mathcal{I}^{(2)}(D; \mc{O}_j) )$ the measured descriptors associated to the small target $D$, which are computed from the reconstructed CGPTs $\bbM^{(K_j)}$.

\begin{figure}[H]
	\begin{tikzpicture}[scale=0.65]
	\begin{axis}[hide axis, axis lines=middle,xmin=-1.8,xmax=1.8, ymin=-1.5,ymax=1.5, enlargelimits,xtick=\empty, ytick=\empty, ylabel={$$}, xlabel={$$}, title={\LARGE Multi-scale aquisition setting$$},at={(-0.5\linewidth,0.2\linewidth)},width=0.9\textwidth]
	\node[anchor=south west,inner sep=0] (image) at (250,25) {\includegraphics[width=0.2\textwidth]{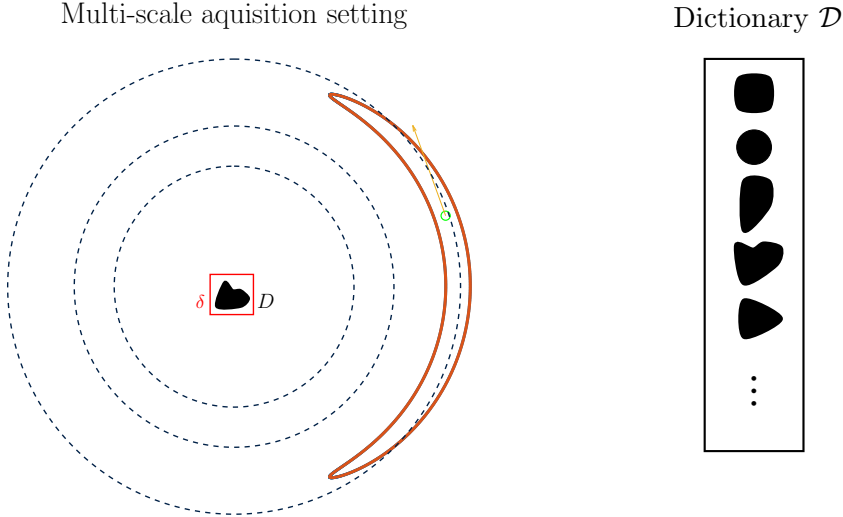}};
	%		\node[anchor=south west,inner sep=0] (image) at (250,20) {\includegraphics[width=0.2\textwidth]{sss2}};	
	%	\node[anchor=south west,inner sep=0] (image) at (100,20) {\includegraphics[width=0.5\textwidth]{sss1}};
	%		\node[anchor=south west,inner sep=0] (image) at (250,20) {\includegraphics[width=0.5\textwidth]{ssss2}};
	%	\node[circle,inner sep=1pt,left] at (axis cs:-1.5,0.4) {\textcolor{arsenic}{\small{$K_1$}}}; 
	%	\node[circle,inner sep=1pt,left] at (axis cs:-0.9,0.7) {\textcolor{arsenic}{\small{$K_2$}}};

	%\draw plot [smooth cycle] coordinates {(0,3) (6,0.35) (6.5, 0.2) (7,0.5) (7,1.65) (6.5,2.75) (.80\textwidth,0.75\textwidth) (0.3\textwidth,.45\textwidth) (1,0.85) } node at (6,1.7) {$D'$};
	\draw[ultra thick, draw=black,fill=black] plot [smooth cycle, xshift=4.5cm,yshift=4.4cm, scale=6] coordinates {(0,-1) (1,2) (2,1) (3,1) (4.,0) (3,-1)}  node[xshift=4.5cm,yshift=4cm, above=0.3cm, right=0.5cm] at (6,1.7) {\large $D$};
	\draw[red, thick,xshift=1.2cm,yshift=1cm,scale=250] (0.47,0.45) node[left=0.2cm,above=0cm] {$\mathsf{\delta}$} rectangle (0.6,0.57);

	\addplot[domain=0:2*pi, samples=400,smooth, dashed, thick,oxfordblue]({1.7*sin(deg(x))},{1.7*cos(deg(x))});

	\addplot[domain=0:2*pi, samples=400,smooth,dashed, thick,oxfordblue]({1.2*sin(deg(x))},{1.2*cos(deg(x))});

	\addplot[domain=0:2*pi, samples=400,smooth,dashed, thick,oxfordblue]({0.9*sin(deg(x))},{0.9*cos(deg(x))});

	%		
	%		\addplot[domain=0.9:1.1, samples=400,smooth, thick,  dotted,oxfordblue]({x*sin(deg(0))},{x*cos(deg(0))});
	%		\addplot[domain=0.9:1.1, samples=400,smooth, thick,  dotted,oxfordblue]({x*sin(deg(pi+pi/4))},{x*cos(deg(pi+pi/4))});
	%		\addplot[domain=0.9:1.1, samples=400,smooth, thick,  dotted,oxfordblue]({x*sin(deg(pi/2+pi/4))},{x*cos(deg(pi/2+pi/4))});
	%		\addplot[domain=0.9:1.1, samples=400,smooth, thick,  dotted,oxfordblue]({x*sin(deg(pi/2+pi+pi/4))},{x*cos(deg(pi/2+pi+pi/4))});
	%		\addplot[domain=0.9:1.1, samples=400,smooth, thick,  dotted,oxfordblue]({x*sin(deg(pi/2-pi))},{x*cos(deg(pi/2-pi))});
	%		\addplot[domain=0.9:1.1, samples=400,smooth, thick,  dotted,oxfordblue]({x*sin(deg(pi/2-pi/4))},{x*cos(deg(pi/2-pi/4))});
	%		\addplot[domain=0.9:1.1, samples=400,smooth, thick,  dotted,oxfordblue]({x*sin(deg(pi))},{x*cos(deg(pi))});
	%		\addplot[domain=0.9:1.1, samples=400,smooth, thick,  dotted,oxfordblue]({x*sin(deg(pi/2))},{x*cos(deg(pi/2))});
	
	%\addplot[domain=0:2*pi, samples=400,smooth,thick,oxfordblue]({2*sin(deg(x))},{2*cos(deg(x))});

	%	\draw[<->,thick,red,dashed] (0,5)  node [left] {$K_3$}  to  [out=0, in=120] node [left] {$$} (2.6,4);	
	%	\draw[<->,thick,red,dashed] (0.5,5.5)  node [left] {$K_2$}  to  [out=0, in=120] node [left] {$$} (2.9,4.5);	
	%	\draw[<->,thick,red,dashed] (1,6)  node [left] {$K_1$}  to  [out=0, in=140] node [left] {$$} (3.1,5.4);	

	\end{axis}
	\begin{scope}[shift={(5,2.5)}]
	\draw[ultra thick, draw=black,fill=black, rotate=30] plot [smooth cycle, xshift=6cm,yshift=3.05cm, scale=0.25] coordinates {(0,-1) (1,2) (2,1) (3,1) (4.,0) (3,-1)};
	
	\draw[ultra thick, draw=black,fill=black] plot [smooth cycle, xshift=3.8cm,yshift=7cm, scale=0.5] coordinates {(0,-1) (0,1) (1,1) (1,0)};
	
	\filldraw[xshift=4cm,yshift=8.2cm] (0,0) circle (10pt);
	
	\draw[ultra thick, draw=black,fill=black] plot [smooth cycle, xshift=3.7cm,yshift=9cm, scale=0.6] coordinates {(0,0) (0,1) (1,1) (1,0)} ;

	\draw[ultra thick, draw=black,fill=black] plot [smooth cycle, xshift=3.8cm,yshift=4.7cm, scale=0.37] coordinates {(0,-1) (0,1) (2,0)} ;
	\draw[black, thick,xshift=3cm,yshift=7cm,scale=1] (0,3) node[left=-0.7cm,above=0.2cm] {Dictionary $\mathcal{D}$} rectangle (2,-5);
	
	\node[xshift=2.6cm,yshift=2.1cm, scale=1.5] (0,-1) {$\vdots$};
	
	\end{scope}		
	\end{tikzpicture}	
	\caption{MSR data are collected by swimming along multiple concentric orbits as shown on the left. The classification problem is to use the features extracted from the data in order to classify, up to rotation and scaling, the small dielectric target $D$ among the elements of a dictionary $\mathcal{D}$.}
	\label{fig:multiscale-dictionary}
\end{figure}

Given a dictionary $\mc{D}$ of $N$ standard shapes, which are denoted by $B_1 , B_2 , ... , B_N$, we want to design a matching algorithm, which generalizes the one proposed in \cite{Ammari11652} to a multi-scale configuration, see Figure \ref{fig:multiscale-dictionary}. Firstly, a matching procedure as in \cite{Ammari11652} is required on each orbit, which consists of a comparison between the theoretical shape descriptors $((\mathcal{I}^{(1)}(B_\kappa) ,  \mathcal{I}^{(2)}(B_\kappa) ))_{\kappa = 1 , ... , N}$ and the measured ones $(\mathcal{I}^{(1)}(D; \mc{O}_j) ,  \mathcal{I}^{(2)}(D; \mc{O}_j) )$, up to a properly chosen length-scale dependent order $K_j$.

Let us define the following scores:
\begin{equation} \label{eq:scores} \Delta_j (B_\kappa , D) =  \left ( \| \mathcal{I}^{(1)}(B_\kappa) - \mathcal{I}^{(1)}(D; \mc{O}_j) \|_F^2 + \| \mathcal{I}^{(2)}(B_\kappa) - \mathcal{I}^{(2)}(D; \mc{O}_j) \|_F^2 \right )^{1/2}, \end{equation}
%\begin{equation} \Delta_{j\,n}(D) =  \left ( \| \mathcal{I}^{(1)}(B_n) - \mathcal{I}^{(1)}(D; \mc{O}_j) \|_F^2 + \| \mathcal{I}^{(2)}(B_n) - \mathcal{I}^{(2)}(D; \mc{O}_j) \|_F^2 \right )^{1/2}, \end{equation}
where $\| \cdot \|_F$ denotes the Frobenius norm of matrices, $j = 1 , ... , \mc{M}$ and $\kappa = 1 , ... , N$.

For every $j$, the scores \eqref{eq:scores} are used to perform the (local) comparison. Precisely, let $\phi_j(D) \in \mc{D}^N$ be the vector formed by the elements of the dictionary, rearranged in ascending order according to $\Delta_j$, i.e., 
\begin{equation*} \phi_j(D) = (B_{\sigma_j(1)} , ... , B_{\sigma_j(N)} ), \end{equation*}
where $\sigma_j$ is a permutation such that $\Delta_j (B_{\sigma_j(l')} , D) \le \Delta_j (B_{\sigma_j(l)} , D)$ , for each $l < l'$.

Notice that, for efficiency reasons, it is convenient to cut the vector $\phi_j(D)$ retaining the first $n \le N$ components only, which are the elements of $\mc{D}$ that produce the lowest scores.

Thus far we have just sorted the elements of the dictionary by matching the descriptors on each orbit separately.  Instead of simply returning $B_{\sigma_j(1)}$ for each $j$, that is in fact the algorithm in \cite{Ammari2014} applied on each orbit for a fixed reconstructing order $K_j$, we aim at fusing the descriptors at the score level. Of course, the scores corresponding to descriptors which have different orders are not directly comparable. The proposed approach is inspired by \cite{MONDEJARGUERRA201557}.

We consider, for every orbit $\mc{O}_j$, the evidence distribution
\begin{equation*} \pi_j := \pi_j(D) = (\eta_{j1} , ... , \eta_{jn}) , \end{equation*}
where
\begin{equation*} \eta_{j \kappa} = \begin{cases} \varphi_j \left ( \frac{\Delta_j(B_{\sigma_j(1)},D)}{\Delta_j(B_{\kappa},D)} \right )^\beta, & \mbox{ if } \sigma_j(\kappa) \le n, \\ \, 0  \;, &  \mbox{ otherwise}. \end{cases}  \end{equation*}
Here, $\varphi_j$ is a normalization constant such that the integral of the evidence distribution $\pi_j$ is $1$, and $\beta$ is a smoothing parameter.

Besides $\pi_j$, we may consider the Shannon's entropy as a confidence factor. Tracing \cite{Huynh2010AdaptivelyEW}, we define
\begin{equation*} c_j := c_j(D) = 1 - \frac{\mc{H}(\pi_j)}{\log(n)}  , \end{equation*}
where $\mc{H}(\pi_j) =  \sum_{\kappa = 1}^n - \eta_{j\kappa} \log(\eta_{j\kappa})$ is the Shannon entropy of the distribution $\pi_j$.

Then, for each $j$, we define the basic belief assignment (BBA)
%\begin{align*} & \mathsf{m}_j(B_{\kappa})  = \eta_{j \kappa} ,\\ &\mathsf{m}_j(\mc{D}) = c_j. \end{align*}
\begin{align*}
\mathsf{m}_j :=(\mathsf{m}_j(B_1),...,\mathsf{m}_j(B_N),\mathsf{m}_j(\mc{D}))\propto (\eta_{j 1},...,\eta_{j N}, c_j).
\end{align*}
$\mathsf{m}_j$ quantifies the evidence given to each element of $\mc{D}$  by the comparison of the descriptors on the $j$-th orbit.

There exist many ways to fuse the evidence which are expressed as BBAs. One of the simplest formulas is the TBM conjunctive rule \eqref{eq:TBM-conjunctive-rule}, which is associative.
Therefore, we start with blending $\mathsf{m}_1$ in with $\mathsf{m}_2$, obtaining $\mathsf{m}_{12} := \mathsf{m}_{1} \myointersection \mathsf{m}_2$. Then we combine $\mathsf{m}_{12}$ with $\mathsf{m}_{3}$, obtaining
\[\mathsf{m}_{123} := (\mathsf{m}_{1} \myointersection \mathsf{m}_2) \myointersection \mathsf{m}_{3} = \mathsf{m}_{1} \myointersection \mathsf{m}_2 \myointersection \mathsf{m}_{3}, \]
and so on, until we compute $\mathsf{m}_{12...\,\mathcal{M}}$.

\medskip

Finally, from the fused BBAs we define the pignistic probability 
\eqref{eq:pignisticprob}, which is used to select the best candidate among the elements of $\mc{D}$. 

\medskip

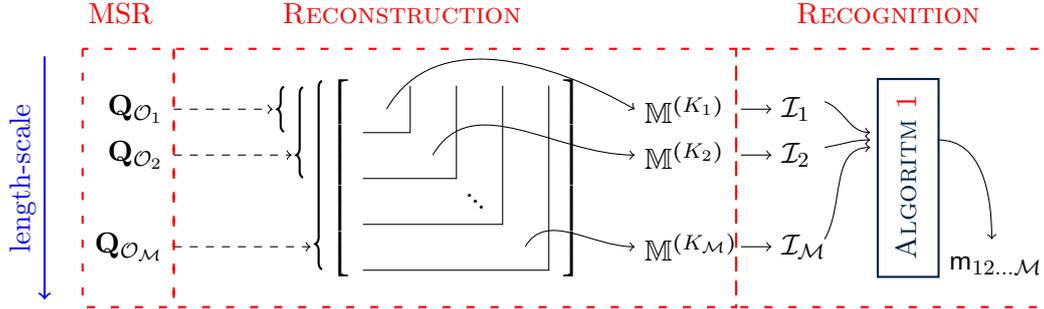
\begin{figure}[H]
	\centering	
	\begin{tikzpicture} 
	\matrix[
	matrix of math nodes,
	nodes in empty cells,
	left delimiter={[},
	right delimiter={]},
	nodes={text height=8pt,text depth=2pt,text width=10pt}
	] (mat)
	{
		& & &  \\
		& & &  \\
		& &\ddots &  \\
		& & &  \\
	};
	\foreach \Valor in {1,...,4}
	\draw (mat-\Valor-1.south west) -| (mat-1-\Valor.north east);
	\draw[->] 
	(mat-1-1.center) 
	to[out=60,in=150] 
	([xshift=1cm]mat.east|-mat-1-1) 
	node[anchor=west] {$\mathbb{M}^{(K_1)}$}
	;  
	\draw[->] 
	(mat-2-2.center) 
	to[out=60,in=180] 
	([xshift=1cm]mat.east|-mat-2-2) 
	node[anchor=west] {$\mathbb{M}^{(K_2)}$} 
	;
	\draw[->] 
	(mat-4-4.center) 
	to[out=60,in=180] 
	([xshift=1cm]mat.east|-mat-4-4) 
	node[anchor=west] {$\mathbb{M}^{(K_\mc{M})}$}
	;  
	
	%\draw[->] 
	%([xshift=2.3cm]mat.east|-mat-2-2) 
	%to[out=0,in=180]  (5,0.2)
	%;  

	\draw[->] 
	([xshift=2.4cm]mat.east|-mat-1-1) 
	to[out=0,in=180] 
	([xshift=2.8cm]mat.east|-mat-1-1)  node[anchor=west] %{$(\mathcal{I}^{(1)}(D; \mc{O}_1) ,  \mathcal{I}^{(2)}(D; \mc{O}_1) )$}
	{$\mathcal{I}_1$}
	;
	
	\draw[->] 
	([xshift=2.4cm]mat.east|-mat-2-2) 
	to[out=0,in=180] 
	([xshift=2.8cm]mat.east|-mat-2-2) node[anchor=west] {$\mathcal{I}_2$}
	;
	
	\draw[->] 
	([xshift=2.4cm]mat.east|-mat-4-4) 
	to[out=0,in=180] 
	([xshift=2.8cm]mat.east|-mat-4-4) node[anchor=west] 
	{$\mathcal{I}_{\mc{M}}$};

	\draw[->] 
	([xshift=3.5cm,yshift=0.1cm]mat.east|-mat-4-4)
	to[out=60,in=180] 
	([xshift=4.1cm,yshift=0.1cm]mat.east|-mat-2-4);
	
	\draw[->] 
	([xshift=3.5cm,yshift=0.1cm]mat.east|-mat-2-2)
	to[out=10,in=180] 
	([xshift=4.1cm,yshift=0.2cm]mat.east|-mat-2-4);
	
	\draw[->] 
	([xshift=3.5cm,yshift=0.1cm]mat.east|-mat-1-1)
	to[out=340,in=180] 
	([xshift=4.1cm,yshift=0.3cm]mat.east|-mat-2-4);
	
	% AFTER COMBINATION
	\draw[->] 
	([xshift=5cm,yshift=0.2cm]mat.east|-mat-2-4)
	to[out=0,in=90] 
	([xshift=5.7cm,yshift=0.1cm]mat.east|-mat-4-4) node[right=0.05cm,below=0.1cm] {$\mathsf{m}_{12 ... \mc{M}}$};
	%{$(\mathcal{I}^{(1)}(D; \mc{O}_\mc{M}) ,  \mathcal{I}^{(2)}(D; \mc{O}_\mc{M}) )$}
	;
	%	\draw[->] 
	%	(mat-3-3.center) 
	%	to[out=60,in=180] 
	%	([xshift=1cm]mat.east|-mat-3-3) 
	%	node[anchor=west] {some text}
	%	;  

%	\node[draw,circle,inner sep=-0.2pt, scale=1.4] at ([xshift=4.5cm,yshift=0.2cm]mat.east|-mat-2-4)	{$\bigcap$};

	\draw[oxfordblue, thick,xshift=1.2cm,yshift=1cm,scale=1.] ([xshift=4.2cm,yshift=0.4cm]mat.east|-mat-1-4) node[sloped,left=-0.1cm,below=1.3cm, rotate=90] {\textsc{Algoritm} \ref{matching_algorithm}} rectangle ([xshift=5cm,yshift=-0.4cm]mat.east|-mat-4-4);

	%
	%\draw[->,dashed] 
	%([xshift=-6cm] mat-1-1.center) node[anchor=east] 
	%{MSR $\mc{O}_1$}
	%to[out=0,in=180]  
	%([xshift=-1cm]mat-1-1.center) 
	%;  
	%\draw[->,dashed] 
	%([xshift=-6cm] mat-2-1.center) node[anchor=east] 
	%{MSR $\mc{O}_2$}
	%to[out=0,in=180]  
	%([xshift=-1cm]mat-1-1.center) 
	%;  
	%\draw[->,dashed] 
	%([xshift=-6cm] mat-4-1.center) node[anchor=east] 
	%{MSR $\mc{O}_{\mc{M}}$}
	%to[out=0,in=180]  
	%([xshift=-1cm]mat-1-1.center) 
	%;  

	\draw[->,dashed] 
	([xshift=-2.8cm] mat-1-1.center) node[anchor=east] 
	{$\bfQ_{\mc{O}_1}$}
	to[out=0,in=180]  
	([xshift=-1.5cm]mat-1-1.center) 
	;  
	
	\draw[mybrace=0.5, thick] ([xshift=-1.35cm,yshift=-0.3cm]mat-1-1.center)  
	to[out=90,in=270]   ([xshift=-1.35cm,yshift=0.3cm]mat-1-1.center) ;
	
	\draw[->,dashed] 
	([xshift=-2.8cm] mat-2-1.center) node[anchor=east] 
	{$\bfQ_{\mc{O}_2}$}
	to[out=0,in=180]  
	([xshift=-1.25cm]mat-2-1.center) 
	;  
	
	\draw[mybrace=0.25, thick] ([xshift=-1.1cm,yshift=-0.3cm]mat-2-1.center)  
	to[out=90,in=270]   ([xshift=-1.1cm,yshift=0.9cm]mat-2-1.center) ;
	
	\draw[->,dashed] 
	([xshift=-2.8cm] mat-4-1.center) node[anchor=east] 
	{$\bfQ_{\mc{O}_{\mc{M}}}$}
	to[out=0,in=180]  
	([xshift=-1cm]mat-4-1.center) 
	;

	\draw[mybrace=0.135, thick] ([xshift=-0.85cm,yshift=-0.3cm]mat-4-1.center)  
	to[out=90,in=270]   ([xshift=-0.85cm,yshift=2.2cm]mat-4-1.center) ;
	
	%% time line %%
	\draw[->,thick, blue] 
	([xshift=-4.5cm,yshift=0.7cm] mat-1-1.center) 
	to[out=270,in=90]  node[sloped,midway,below=0.6cm, rotate=180]{length-scale}
	([xshift=-4.5cm,yshift=-0.7cm] mat-4-1.center)
	;  
	
	\draw[red, thick, loosely dashed] ([xshift=-4cm,yshift=0.8cm] mat-1-1.center)  node[right=0.5cm,above=0.2cm] {MSR} rectangle ([xshift=-2.8cm,yshift=-0.8cm]mat-4-1.center) ;
	
	\draw[red, thick, loosely dashed] ([xshift=-2.8cm,yshift=0.8cm] mat-1-1.center)  node[right=3cm,above=0.2cm] {\textsc{Reconstruction}} rectangle ([xshift=4.6cm,yshift=-0.8cm]mat-4-1.center) ;

	\draw[red, thick, loosely dashed] ([xshift=4.6cm,yshift=0.8cm]mat-1-1.center)  node[right=2cm,above=0.2cm] {\textsc{Recognition}} rectangle ([xshift=8.7cm,yshift=-0.8cm]mat-4-1.center) ;

	%%% radius line %%
	%\draw[<-,thick, blue] 
	%([xshift=-8cm,yshift=0.7cm] mat-1-1.center) 
	%to[out=270,in=90]  node[sloped,midway,above=-0.5cm]{distance}
	%([xshift=-8cm,yshift=-0.7cm] mat-4-1.center)
	%;  
	
	\end{tikzpicture}
	\caption{Overview of the three relevant stages involved in our multi-scale approach. For conciseness, we denote $(\mathcal{I}^{(1)}(D; \mc{O}_j) ,  \mathcal{I}^{(2)}(D; \mc{O}_j) )$ by $\mc{I}_j$.}
	\label{fig:overview}
\end{figure}

The procedure described hereabove is summarised in Algorithm \ref{matching_algorithm}. See also Figure \ref{fig:overview}.

\bigskip
\begin{algorithm}[H]
	\Input{On each orbit $\mc{O}_j \in \{\mc{O}_1 , ... , \mc{O}_{\mc{M}} \}$, the first $k$-th order shape descriptors $\mathcal{I}^{(1)}(D, \mc{O}_j)$, $\mathcal{I}^{(2)}(D, \mc{O}_j)$ of an unknown target $D$.}
	\nl \For{$\mc{O}_j \in \{\mc{O}_1 , ... , \mc{O}_{\mc{M}} \}$}{
		%        \nl \;
		\nl \For{$B_\kappa \in \{ B_1, ... , B_N\}$}{
			\nl  $\Delta_j (B_\kappa,D) \leftarrow \left ( \| \mathcal{I}^{(1)}(B_\kappa) - \mathcal{I}^{(1)}(D; \mc{O}_j) \|_F^2 + \| \mathcal{I}^{(2)}(B_\kappa) - \mathcal{I}^{(2)}(D; \mc{O}_j) \|_F^2 \right )^{1/2}$ \;
			
		}
		\nl $\sigma_j(1) \leftarrow \text{argmin}_\kappa \Delta_j(B_\kappa,D)$\;
		\nl \For{$B_\kappa \in \{B_1, ... , B_N\}$}{
			\nl $\eta_{j\kappa} \leftarrow  \varphi_j  \left( \frac{\Delta_j(B_{\sigma_j(1)},D)}{\Delta_j(B_\kappa,D)} \right)^\beta$ \;
			\nl $\mathsf{m}_j(B_\kappa) \leftarrow \eta_{j\kappa}$\;
		}
		\nl $c_j \leftarrow 1 - \frac{1}{\log(N)}\sum_{\kappa = 1}^N- \eta_{j\kappa} \log(\eta_{j\kappa})$ \;
		\nl $\mathsf{m}_j(\mc{D}) \leftarrow c_{j}$\;
		%        \nl $d (\mc{O}_j) \leftarrow \frac{e (\mc{O}_j)}{\max_n e_n(\mc{O}_j)}$\;
	}
	\nl $\mathsf{m}_{1 ...\mc{M}} \leftarrow \mathsf{m}_{1} \myointersection ... \myointersection \mathsf{m}_{\mc{M}}$\;
	\Output{the best matching element of the dictionary $\kappa^* \leftarrow \text{argmax}_{\kappa} \text{BetP}(B_{\kappa})$.}
	\caption{Shape identification for a multi-scale setting based on transform invariant descriptors}
	\label{matching_algorithm}
\end{algorithm}

\section{Numerical results}

\label{sec:numerical-exp}
In this section, we show some numerical results which illustrate how Algorithm \ref{matching_algorithm} can significantly improve the robustness of the recognition procedure.

\subsection{Setting} \label{subsec-setting}
\begin{figure}[H]
	\centering
	\begin{tabular}{ccccc}
		%		& 1 & 2  & 3  \\
		%		& & &  \\
		\quad &  \raisebox{7ex - \height}{\includegraphics[scale=0.33]{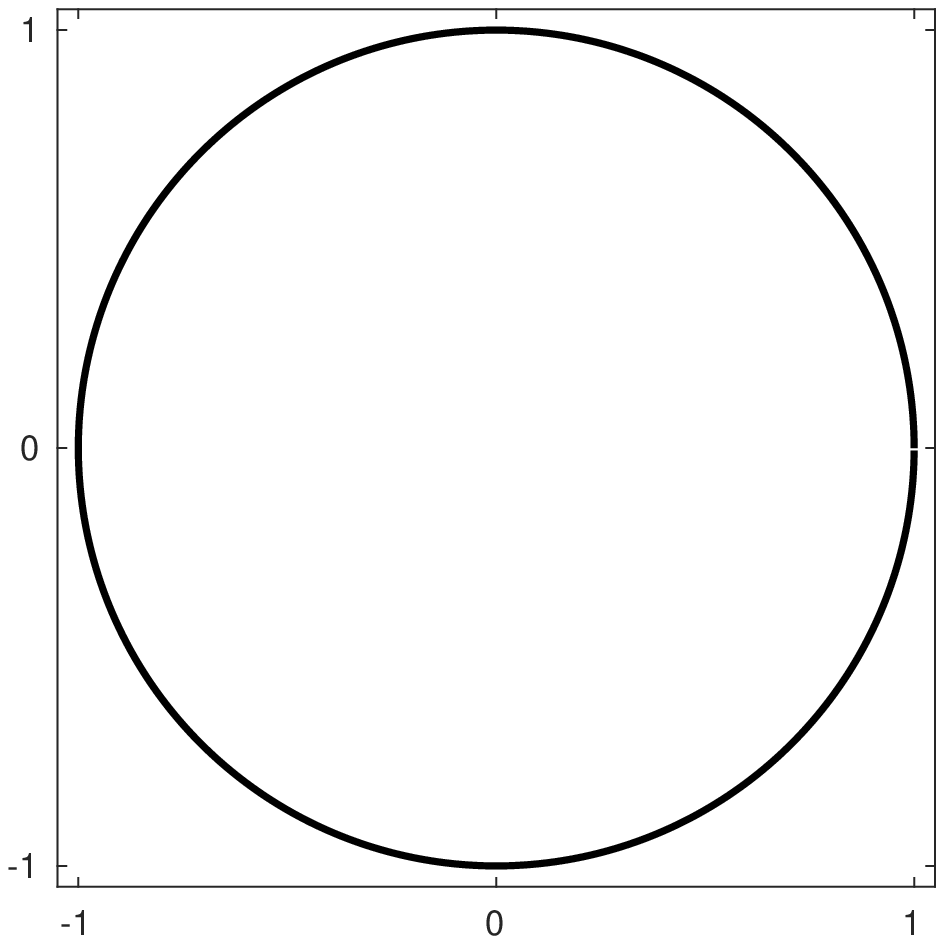}}   &
		\raisebox{7ex - \height}{\includegraphics[scale=0.33]{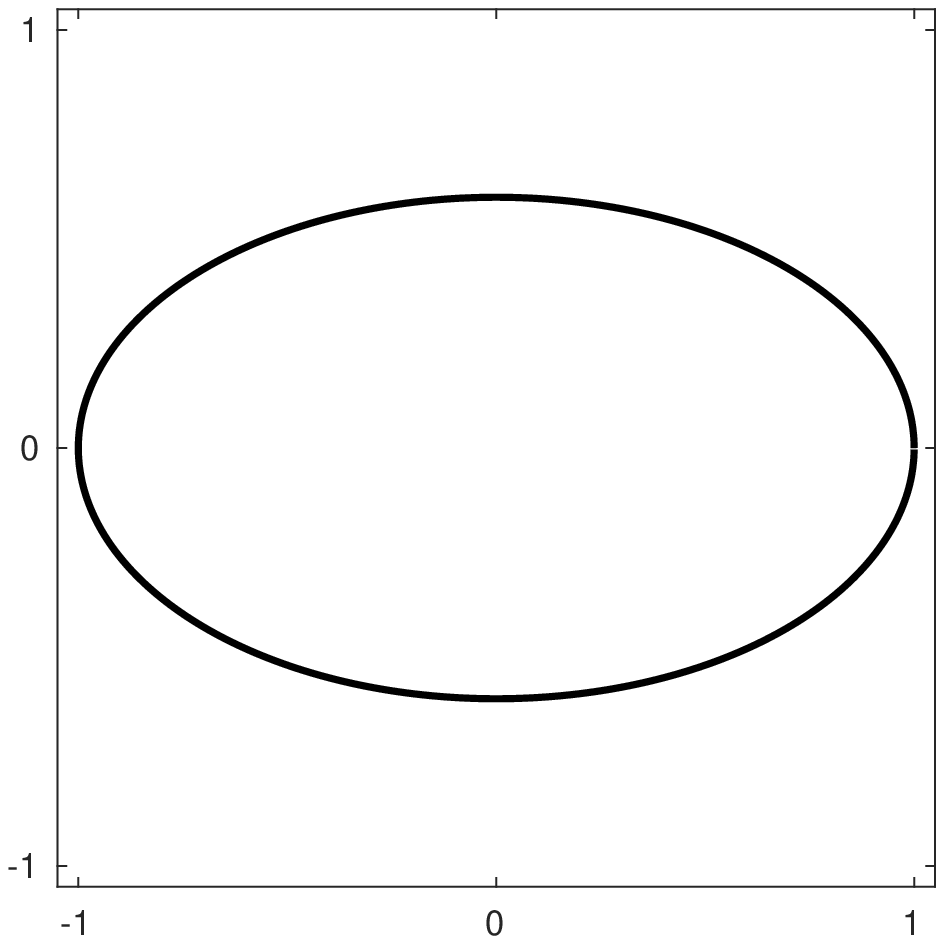}} & 
		\raisebox{7ex - \height}{\includegraphics[scale=0.33]{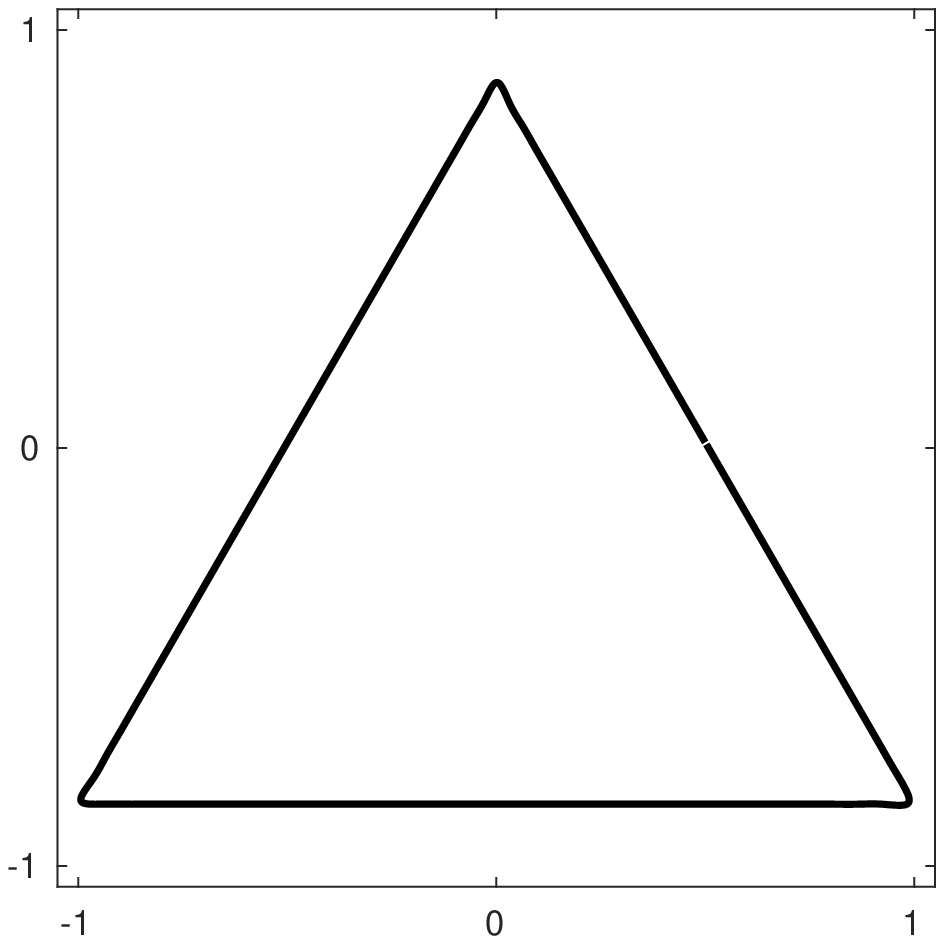}} & 
		\raisebox{7ex - \height}{\includegraphics[scale=0.33]{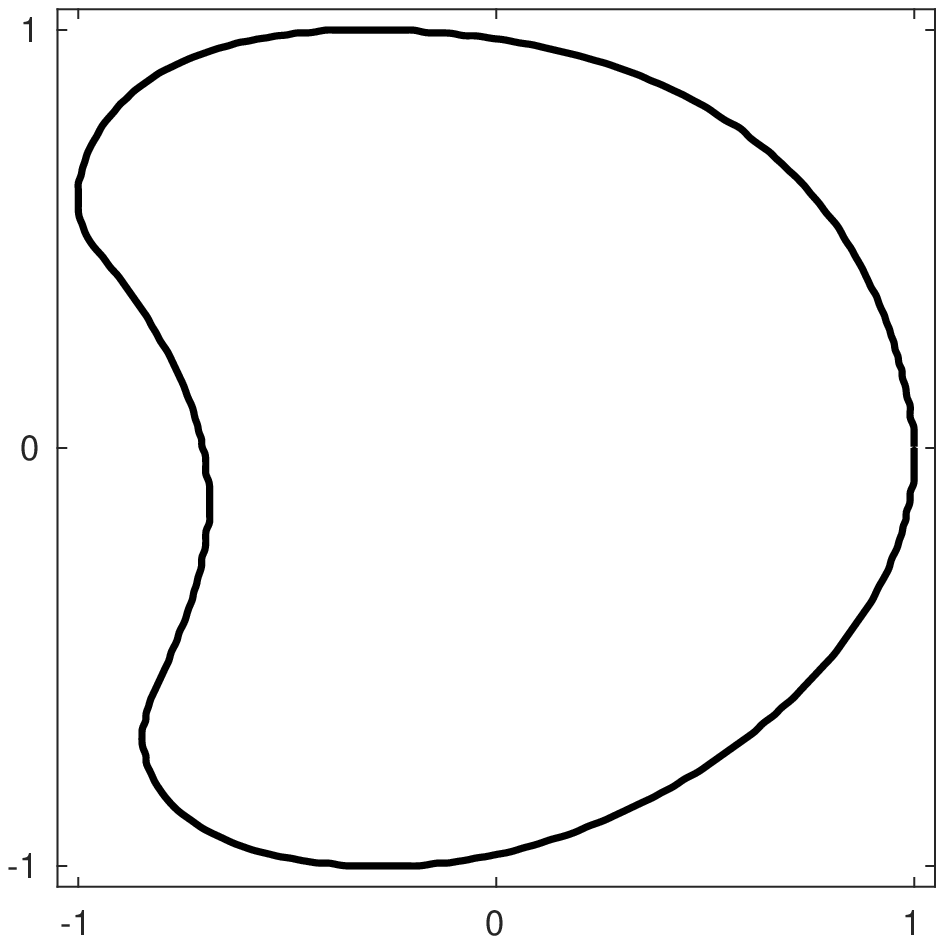}} \\
		& & & &\\
		\quad & \raisebox{7ex - \height}{\includegraphics[scale=0.33]{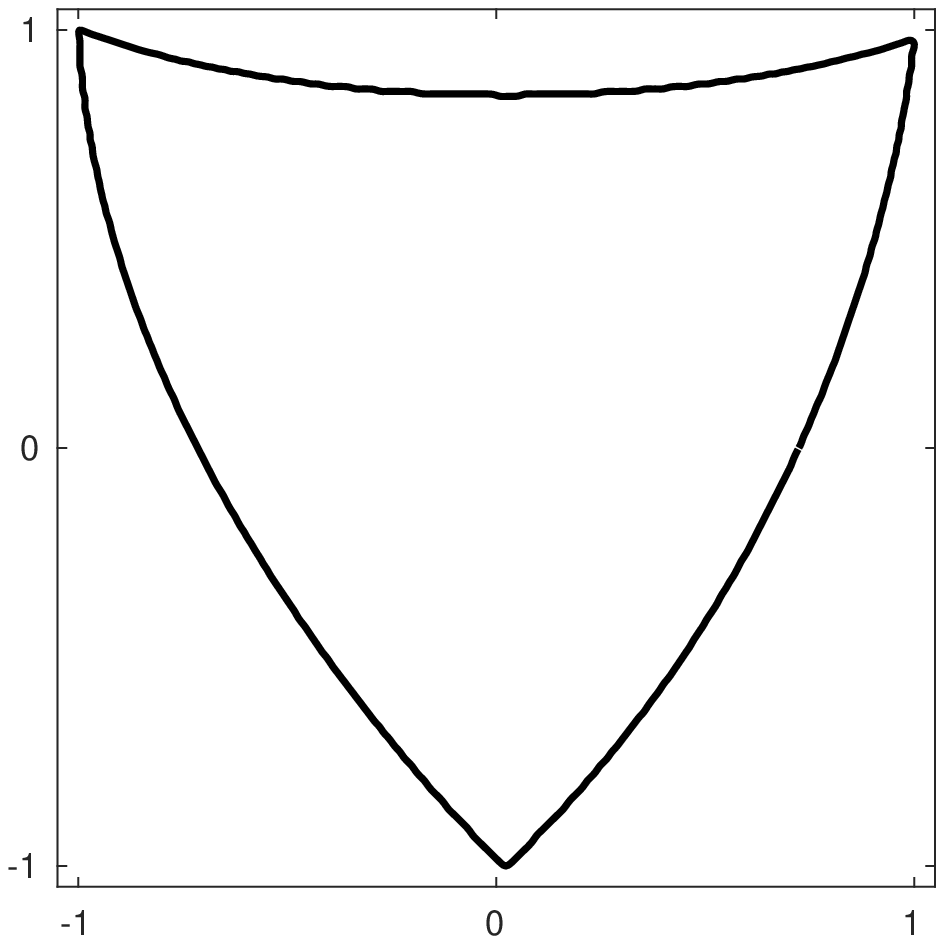}} &
		\raisebox{7ex - \height}{\includegraphics[scale=0.33]{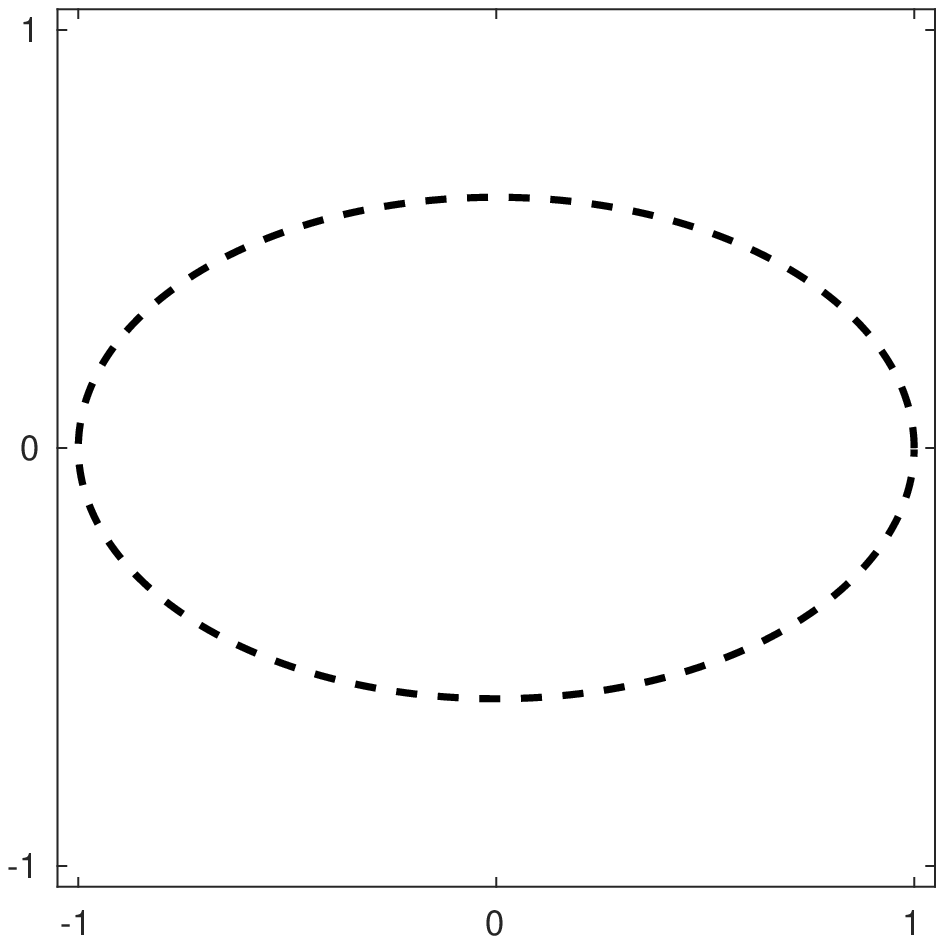}}  &
		\raisebox{7ex - \height}{\includegraphics[scale=0.33]{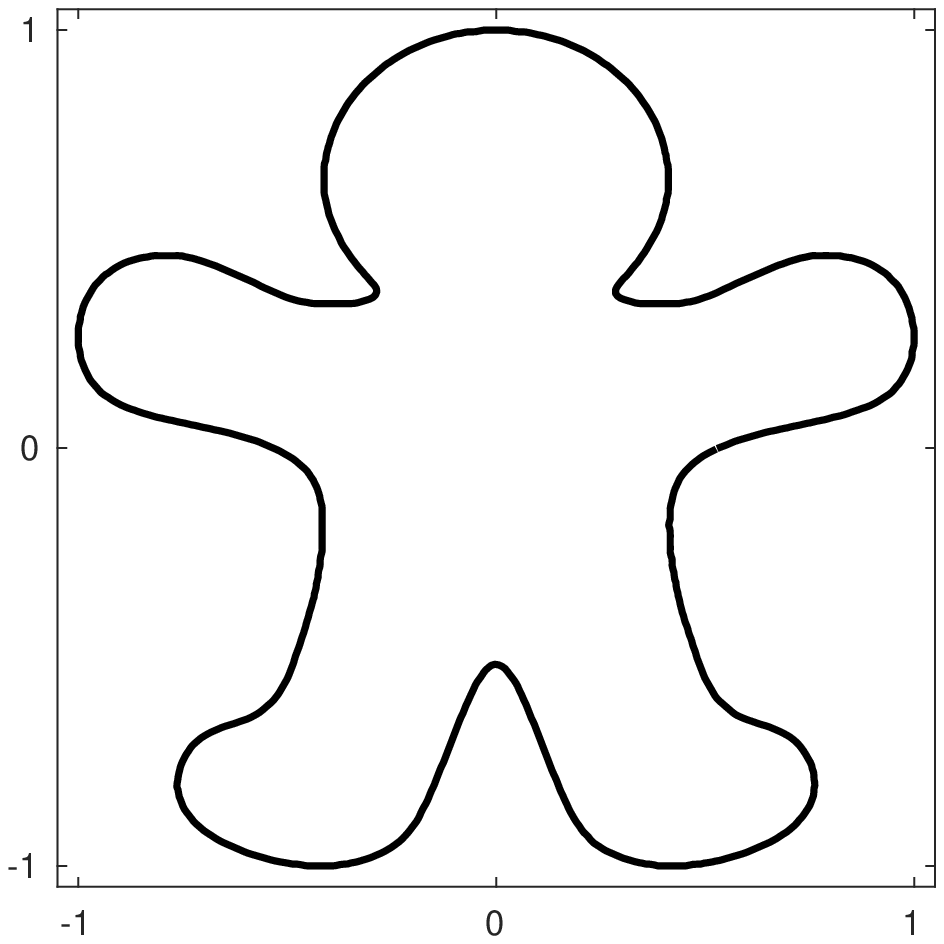}} & 
		\raisebox{7ex - \height}{\includegraphics[scale=0.33]{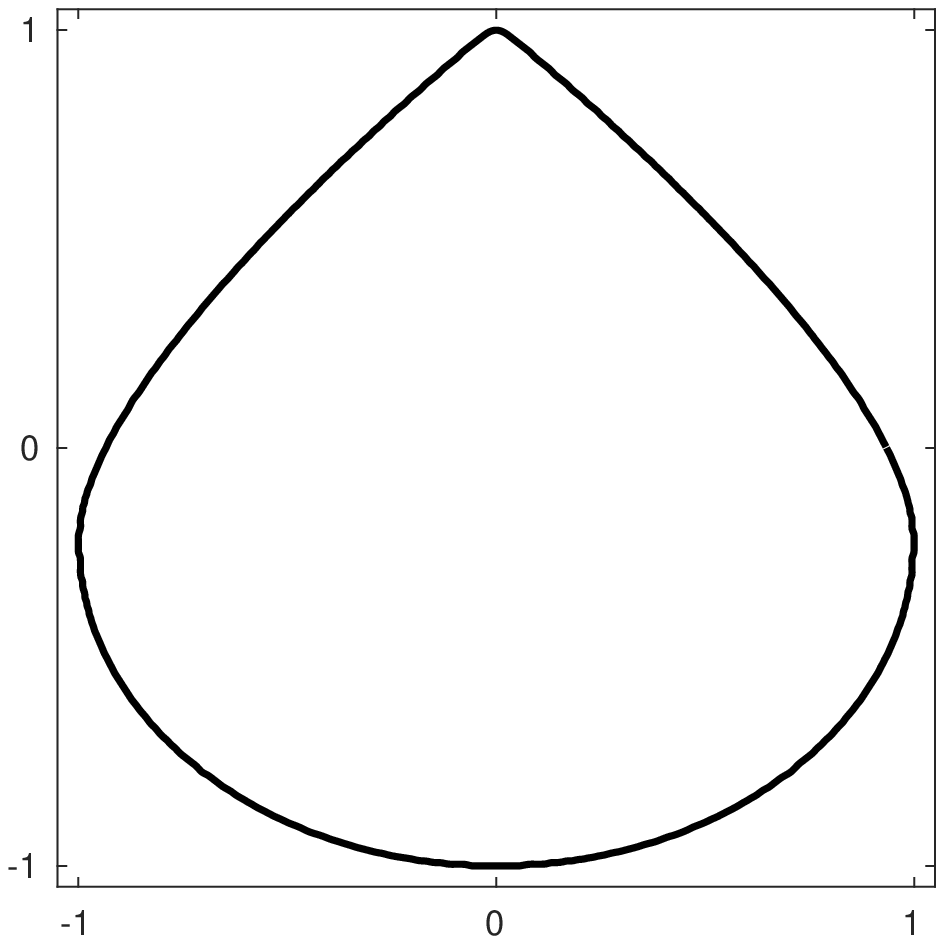}}   	\\	
	\end{tabular}
	\caption{Dictionary $\mathcal{D}$.}
	\label{smalldico}
\end{figure}

Let $\mathcal{D}$ be a dictionary containing 8 standard shapes, as illustrated in \figref{smalldico}. Each solid shape is equipped with homogeneous conductivity having parameter $k = 3$ (Circle, Ellipse, Triangle, Bent Ellipse, Curved Triangle, Gingerbread Man, Drop) whereas the dashed one (Ellipse) has conductivity $k = 10$. All the shapes have the same characteristic size, which is of order one. 

% Describe D

The targets $D$ we are considering for the experiments are located at the origin as the standard shapes, and are obtained by scaling and rotating the elements of $\mc{D}$, with scaling coefficient and  rotation angle chosen as $\delta = 0.2$ and $\theta = \pi/3$,
respectively. 

\FloatBarrier

\subsection{Experiment} \label{subsec-experiment}

The data are acquired as described in Section \ref{sbsec:setting_orbits}. Precisely, given $\mc{O}_1 , \mc{O}_2$ and $\mc{O}_3$ three circular orbits around $D$, of radii $\rho_1 = 1.6$, $\rho_2 = 1.1$, $\rho_3 = 0.9$, respectively, we sample $M_1 = M_2 = M_3 = 200$ positions on each trajectory, and build the corresponding MSR matrices. We consider $N_r = 2^{10}$ receptors evenly distributed on the the body of the fish.

In the numerical experiments the MSR data are simulated using the code developed in \cite{WangRep}.

\begin{figure}[H]
	\centering
	\hspace*{-5mm}
	\includegraphics[scale=0.7]{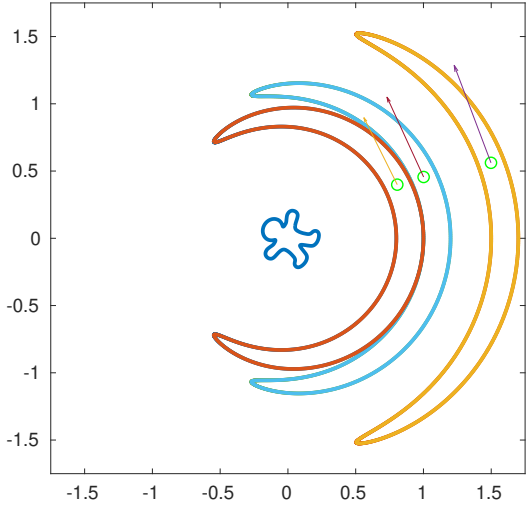}
	\caption{A single position of the electric fish per orbit. The target is located at $(0,0)$.}
	\label{fig:settingOrbits}
\end{figure}

\begin{rem} In principle, the acquisition operator $\mc{L}$ depends on the measurements, which is not a desirable property. In order to overcome this difficulty, we perform the numerical simulations using the surrogate acquisition operator obtained from the dipolar approximation \eqref{eq:dipole_approx}.	
\end{rem}

The CGPTs are reconstructed on each orbit from the MSR matrix by exploiting either the g-inverse $\mathfrak{L}$ or the Moore-Penrose inverse $\bfL^\dagger$. The reconstruction orders are set $K_1 = 2$ on $\mc{O}_1$, $K_2 = 3$ on $\mc{O}_2$, and $K_3 = 4$ on $\mc{O}_3$. The corresponding shape descriptors are then used as a rationale for building the BBAs $\mathsf{m}_1 , \mathsf{m}_2$ and $\mathsf{m}_3$, with parameter $\beta = 2$, as described in Section \ref{sbsec:setting_orbits}. For efficiency reasons,  positive mass is given only to the first $n = 3$ best matching elements of $\mc{D}$.

We study the robustness of the fused descriptors given by Algorithm \ref{matching_algorithm} with moderate noise in the measurements. Precisely, given a target $D$ and $\sigma_{\text{noise}} \in  [0.0025,0.050]$, we test the recognition algorithm by considering $10^4$ experiments, and computing the frequencies.

The results arising from the beliefs produced on each orbit, i.e., $\mathsf{m}_1 , \mathsf{m}_2$ and $\mathsf{m}_3$, are compared with the ones obtained from the fused beliefs synthetized by the TBM conjunctive rule, i.e., $\mathsf{m}_{12} = \mathsf{m}_1 \myointersection \mathsf{m}_2$ and $\mathsf{m}_{123} = \mathsf{m}_{12} \myointersection \mathsf{m}_3$.

\FloatBarrier
\subsubsection{Reconstruction by the generalized inverse $\mathfrak{L}$} \label{subsubsec-experiment_gi}

The results of this part are obtained by employing the reflexive minimum norm g-inverse $\mathfrak{L}$ of the acquisition operator given by \eqref{eq:minnormsol} for the reconstruction of the CGPTs from the MSR data.

\begin{figure}[H]
	\centering
	\hspace*{-5mm}
	\begin{tabular}{cc}
		\includegraphics[scale=0.6]{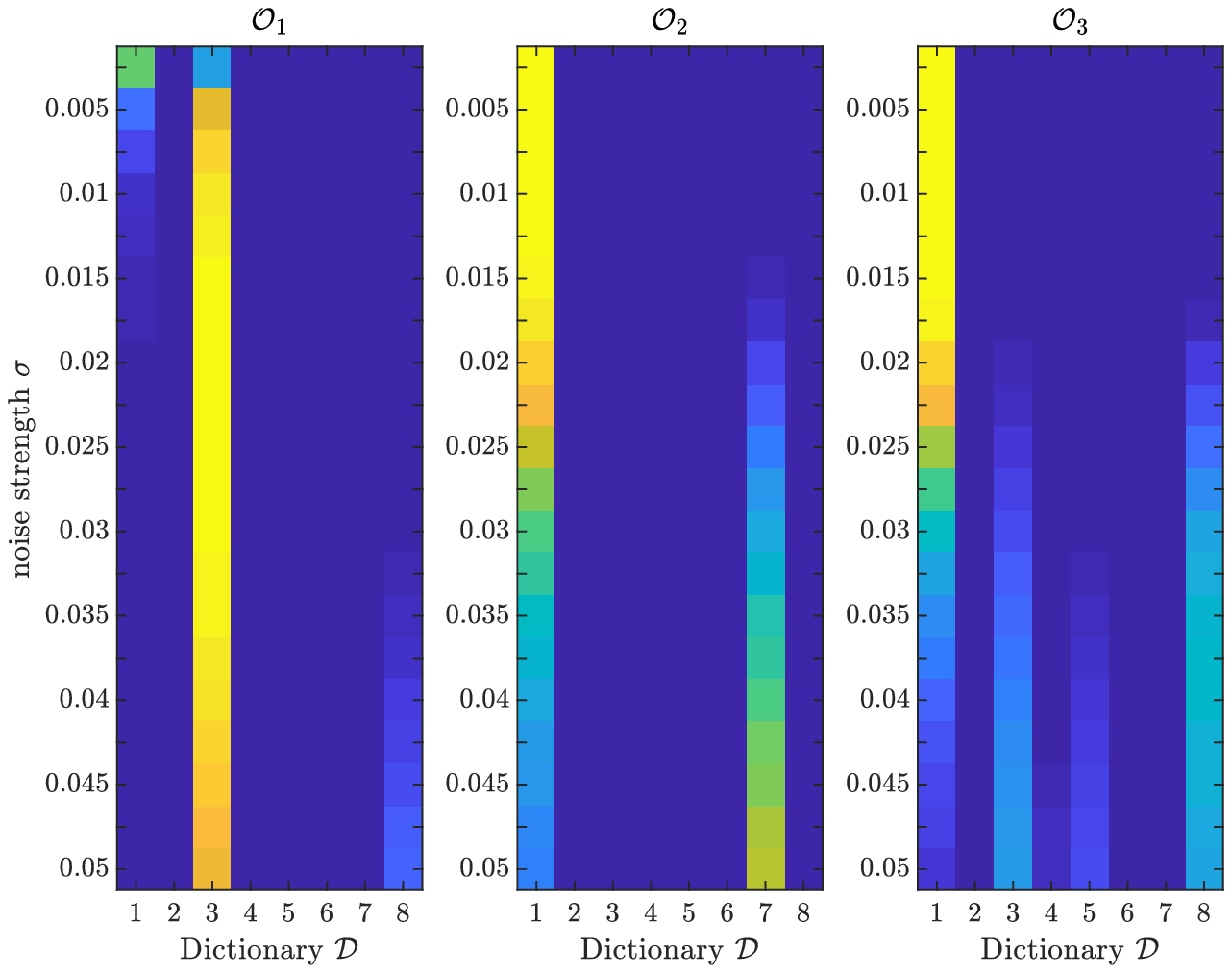} & 
		\fbox{\includegraphics[scale=0.61]{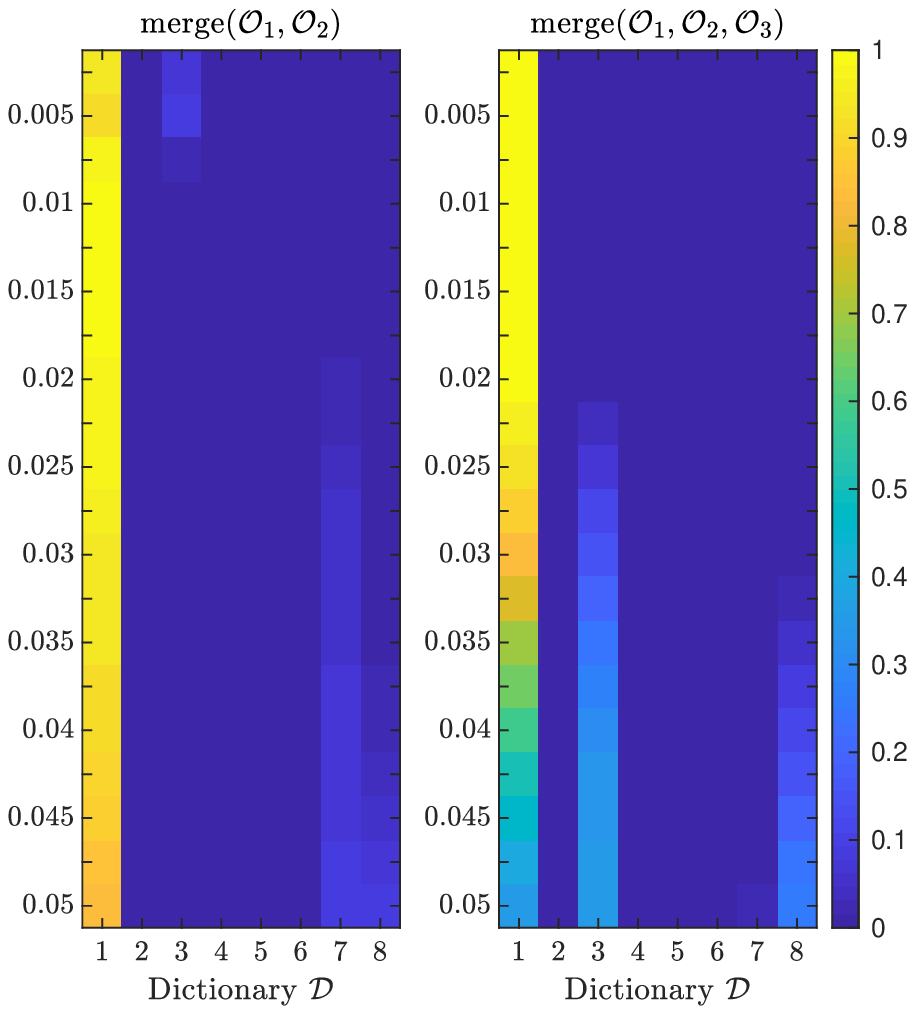}}
	\end{tabular}
	\caption{Circle (1)}
	\label{fig:recoCircle}
\end{figure}

\begin{figure}[H]
	\centering
	\hspace*{-5mm}
	\begin{tabular}{cc}
		\includegraphics[scale=0.6]{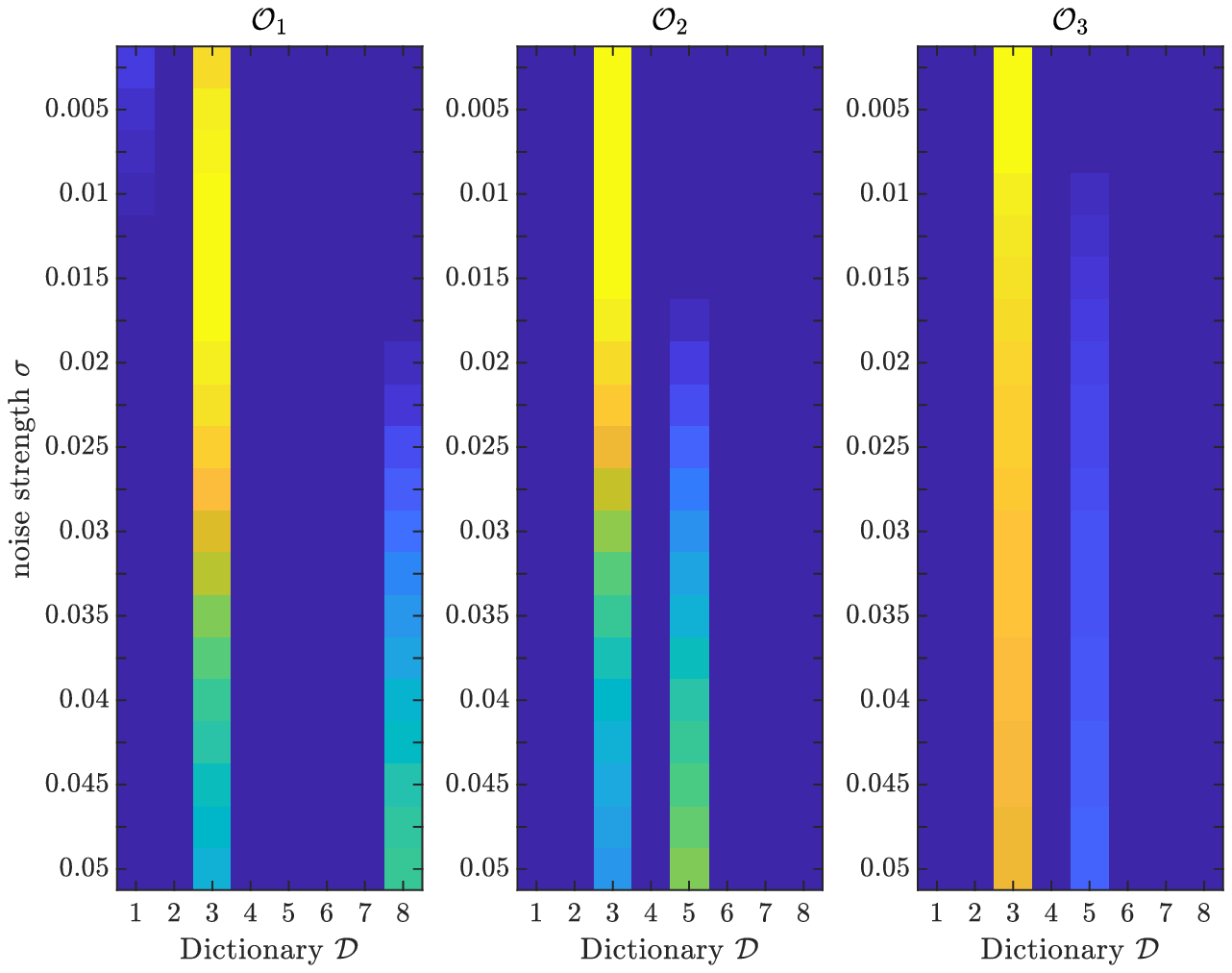} & 
		\fbox{\includegraphics[scale=0.61]{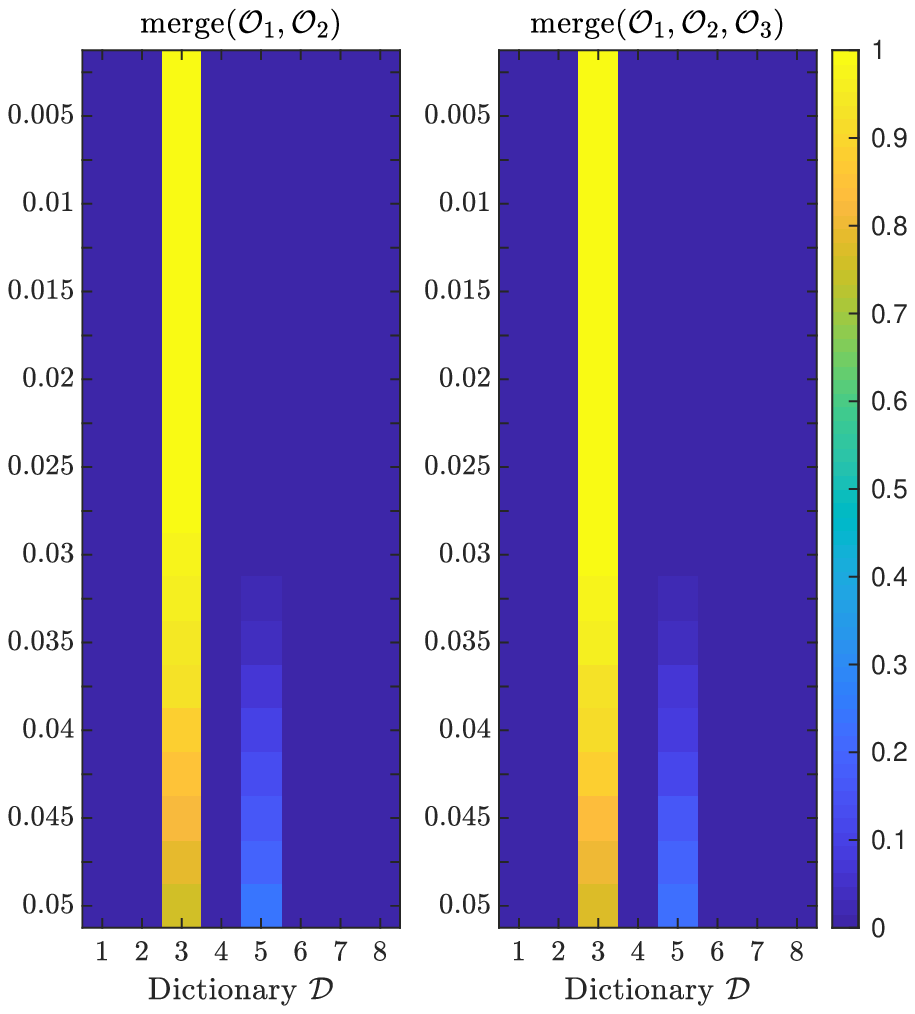}}
	\end{tabular}
	\caption{Triangle (3)}
	\label{fig:recoTriangle}
\end{figure}

\begin{figure}[H]
	\centering
	\hspace*{-5mm}
	\begin{tabular}{cc}
		\includegraphics[scale=0.6]{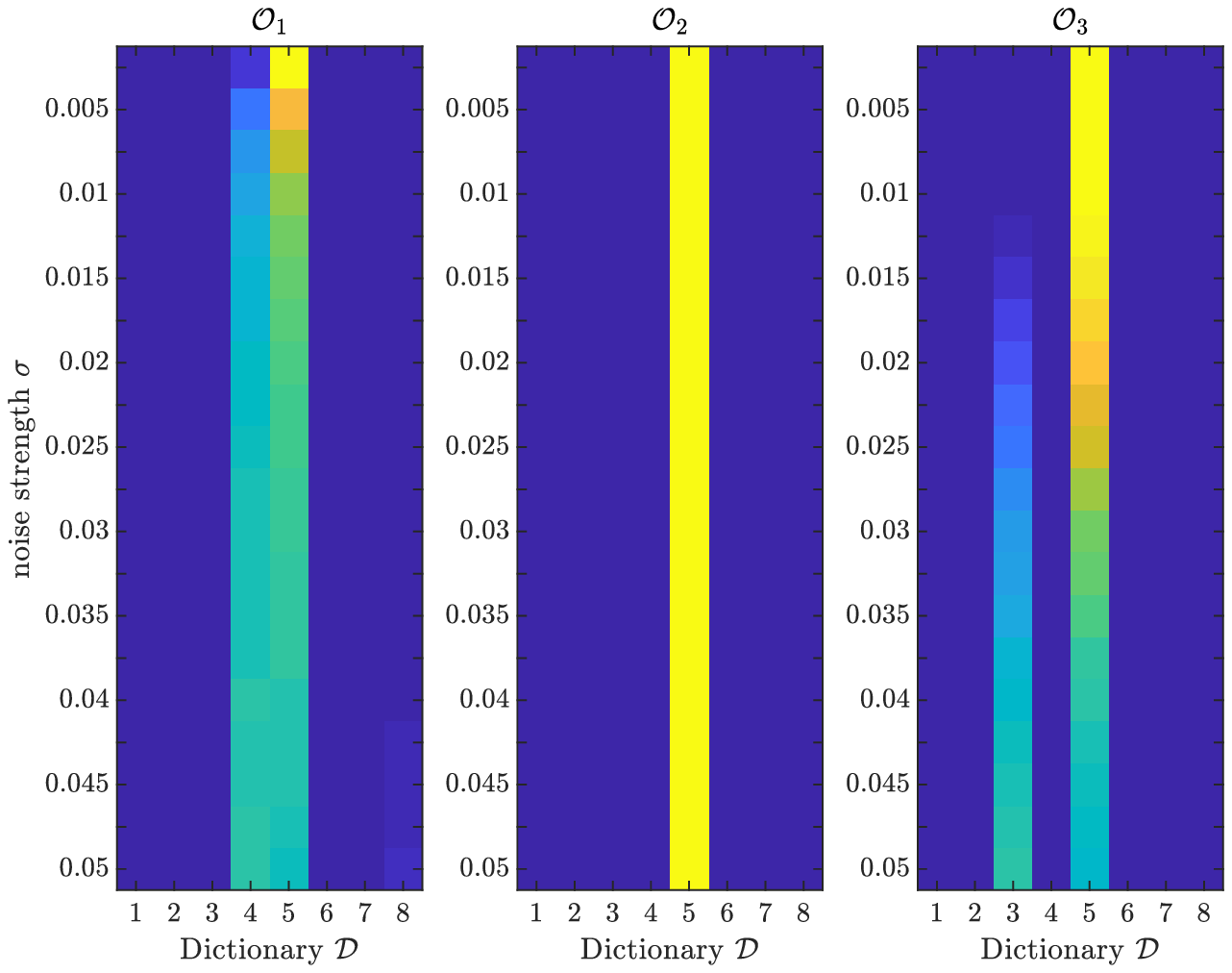} & 
		\fbox{\includegraphics[scale=0.61]{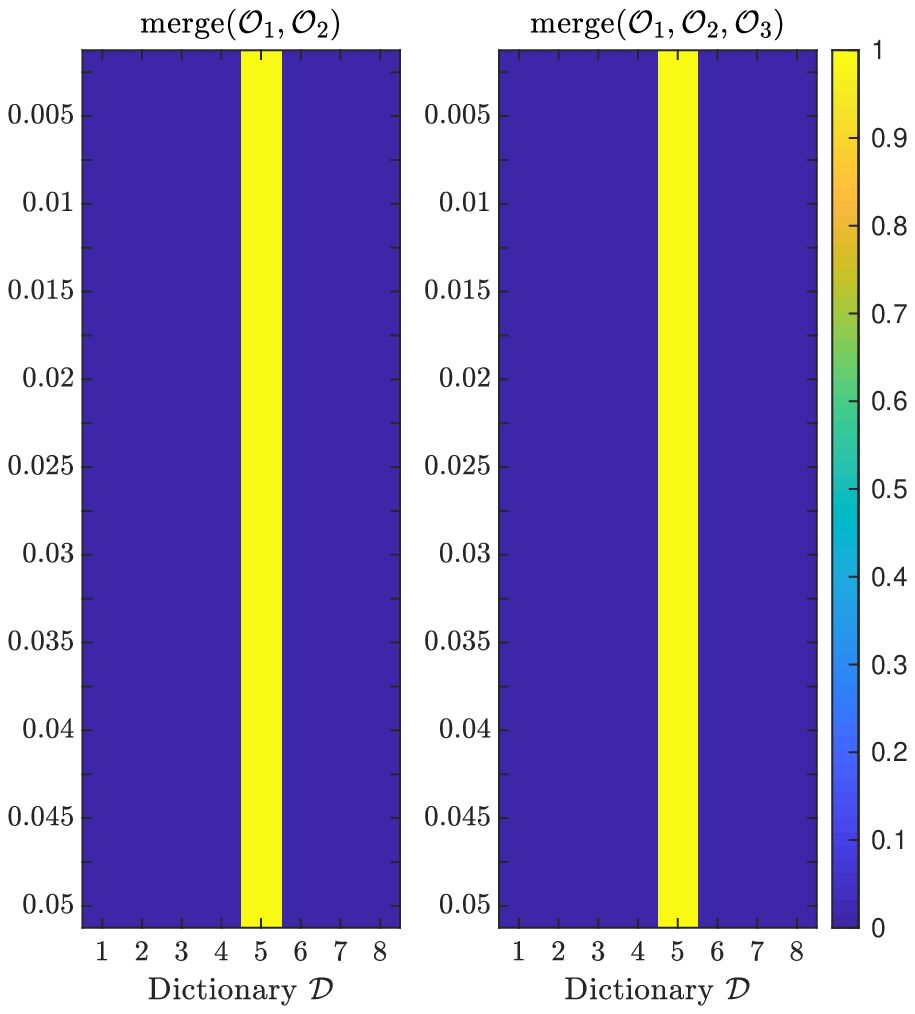}}
	\end{tabular}
	\caption{Curved Triangle (5)}
	\label{fig:recoTriaShield}
\end{figure}

The frequencies are reported in Figures \ref{fig:recoCircle}, \ref{fig:recoTriangle} and \ref{fig:recoTriaShield}.

\begin{table}[H]
	\centering
	\begin{tabular}{l|>{\columncolor[gray]{0.9}}cccccccc}
	\toprule		
	&$\boxed{\text{1}}$&2&3&4&5&6&7&8\\
	\midrule
	$\mathsf{m}_{1}$ &0.0021& \rule{7pt}{0pt}0& \rule{7pt}{0pt}0.8029& \rule{7pt}{0pt}0& \rule{7pt}{0pt}0& \rule{7pt}{0pt}0& \rule{7pt}{0pt}0& \rule{7pt}{0pt}0.1950\\
	$\mathsf{m}_{2}$&0.2805& \rule{7pt}{0pt}0& \rule{7pt}{0pt}0		& \rule{7pt}{0pt}0& \rule{7pt}{0pt}0& \rule{7pt}{0pt}0& \rule{7pt}{0pt}0.7195& \rule{7pt}{0pt}0\\
	$\mathsf{m}_{3}$ &0.0728 & \rule{7pt}{0pt}0& \rule{7pt}{0pt}0.3582& \rule{7pt}{0pt}0.0381& \rule{7pt}{0pt}0.1400& \rule{7pt}{0pt}0& \rule{7pt}{0pt}0.0088& \rule{7pt}{0pt}0.3821\\ 
	\hline
	$\mathsf{m}_{12}$ &0.8288& \rule{7pt}{0pt}0& \rule{7pt}{0pt}0 & \rule{7pt}{0pt}0& \rule{7pt}{0pt}0& \rule{7pt}{0pt}0& \rule{7pt}{0pt}0.0834& \rule{7pt}{0pt}0.0878\\ 
	$\mathsf{m}_{123}$&0.3575& \rule{7pt}{0pt}0& \rule{7pt}{0pt}0.3585& \rule{7pt}{0pt}0.0001& \rule{7pt}{0pt}0.0013& \rule{7pt}{0pt}0& \rule{7pt}{0pt}0.0174& \rule{7pt}{0pt}0.2652\\ 
	\bottomrule
\end{tabular}
	
	\medskip
	
	\caption{Frequency table for the identification of the Circle (1) with the strongest noise, i.e., $\sigma = 0.05$. Each row displays the relative frequencies for all the elements of the dictionary corresponding to different BBAs.}
	\label{frequency-table-reco1}
\end{table}

Looking at Table \ref{frequency-table-reco1} we can clearly see that the combination of classifiers outperforms each classifier. In particular we can see that combining the classifiers on $\mc{O}_1$ and $\mc{O}_2$ yields a great improvement in the recognition already.

%it appears clear that the recognition on each orbit is

%\myworries{TODO}

Tables \ref{frequency-table-reco3} and \ref{frequency-table-reco5} refer to the Triangle (3) and the Curved Triangle (5). Because of their similar silhouette, they are troublesome from a dictionary approach point of view \cite{SCAPIN20191872}, and thus it is interesting to have a close look at what happens in both cases.

\begin{table}[H]
	\centering
	\begin{tabular}{l|cc>{\columncolor[gray]{0.9}}cccccc}
	\toprule
	&1&2&$\boxed{\text{3}}$&4&5&6&7&8\\
	\midrule
	$\mathsf{m}_{1}$ &0.0006 & \rule{7pt}{0pt}0& \rule{7pt}{0pt}0.4229& \rule{7pt}{0pt}0& \rule{7pt}{0pt}0.0115& \rule{7pt}{0pt}0& \rule{7pt}{0pt}0& \rule{7pt}{0pt}0.5650\\
	$\mathsf{m}_{2}$&0& \rule{7pt}{0pt}0& \rule{7pt}{0pt}0.3371 & \rule{7pt}{0pt}0& \rule{7pt}{0pt}0.6629& \rule{7pt}{0pt}0& \rule{7pt}{0pt}0& \rule{7pt}{0pt}0\\
	$\mathsf{m}_{3}$ &0 & \rule{7pt}{0pt}0& \rule{7pt}{0pt}0.8038 & \rule{7pt}{0pt}0& \rule{7pt}{0pt}0.1962 & \rule{7pt}{0pt}0& \rule{7pt}{0pt}0& \rule{7pt}{0pt}0\\ 
	\hline
	$\mathsf{m}_{12}$ &0& \rule{7pt}{0pt}0& \rule{7pt}{0pt}0.7519& \rule{7pt}{0pt}0& \rule{7pt}{0pt}0.2382& \rule{7pt}{0pt}0& \rule{7pt}{0pt}0& \rule{7pt}{0pt}0.0099\\ 
	$\mathsf{m}_{123}$&0& \rule{7pt}{0pt}0& \rule{7pt}{0pt}0.7702& \rule{7pt}{0pt}0& \rule{7pt}{0pt}0.2297& \rule{7pt}{0pt}0& \rule{7pt}{0pt}0& \rule{7pt}{0pt}0.0001\\ 
	\bottomrule
\end{tabular}
	
	\medskip
	
	\caption{Frequency table for the identification of the Triangle (3) with the strongest noise, i.e., $\sigma = 0.05$.}
	\label{frequency-table-reco3}
\end{table}

\begin{table}[H]
	\centering
	\begin{tabular}{l|cccc>{\columncolor[gray]{0.9}}cccc}
	\toprule
	&1&2&3&4&$\boxed{\text{5}}$&6&7&8\\
	\midrule
	$\mathsf{m}_{1}$ &0& \rule{7pt}{0pt}0& \rule{7pt}{0pt}0& \rule{7pt}{0pt}0.5015& \rule{7pt}{0pt}0.4650& \rule{7pt}{0pt}0& \rule{7pt}{0pt}0.0002& \rule{7pt}{0pt}0.0333\\
	$\mathsf{m}_{2}$&0& \rule{7pt}{0pt}0& \rule{7pt}{0pt}0.0064& \rule{7pt}{0pt}0& \rule{7pt}{0pt}0.9936& \rule{7pt}{0pt}0& \rule{7pt}{0pt}0& \rule{7pt}{0pt}0\\
	$\mathsf{m}_{3}$ &0 & \rule{7pt}{0pt}0& \rule{7pt}{0pt}0.5395& \rule{7pt}{0pt}0& \rule{7pt}{0pt}0.4605 & \rule{7pt}{0pt}0& \rule{7pt}{0pt}0& \rule{7pt}{0pt}0\\ 
	\hline
	$\mathsf{m}_{12}$ &0& \rule{7pt}{0pt}0& \rule{7pt}{0pt}0.0011 & \rule{7pt}{0pt}0.0016& \rule{7pt}{0pt}0.9967& \rule{7pt}{0pt}0& \rule{7pt}{0pt}0& \rule{7pt}{0pt}0.0006\\ 
	$\mathsf{m}_{123}$&0& \rule{7pt}{0pt}0& \rule{7pt}{0pt}0.0017& \rule{7pt}{0pt}0& \rule{7pt}{0pt}0.9983 & \rule{7pt}{0pt}0& \rule{7pt}{0pt}0& \rule{7pt}{0pt}0\\ 
	\bottomrule
\end{tabular}
	
	\medskip
	
	\caption{Frequency table for the identification of the Curved Triangle (5) with the strongest noise, i.e., $\sigma = 0.05$. }
	\label{frequency-table-reco5}
\end{table}

It is worth noticing that in Table \ref{frequency-table-reco3}  the highest frequency is not attained by the fusion of the descriptors. Instead, the third orbit alone produces the best matching. However, merging all the three BBAs enhances considerably the classification success rate in the worst case scenario, which is strikingly lower than the rate in the best case scenario.
%It is expected that merging conflicting beliefs .....

Clearly, this is not a drawback of our method. As a matter of fact, since in advance we don't know  which classifier performs the best, the above results indicate that using their combination is a valid -as well as natural- trade-off.

\FloatBarrier

\subsubsection{Reconstruction by the Moore-Penrose inverse} \label{subsubsec-experiment-MP}

In this part we make use of the Moore-Penrose inverse to reconstruct the CGPTs from the MSR simulated data, as shown in \eqref{eq:MPinverse}.

\begin{figure}[H]
	\centering
	\hspace*{-5mm}
	\begin{tabular}{cc}
		\includegraphics[scale=0.6]{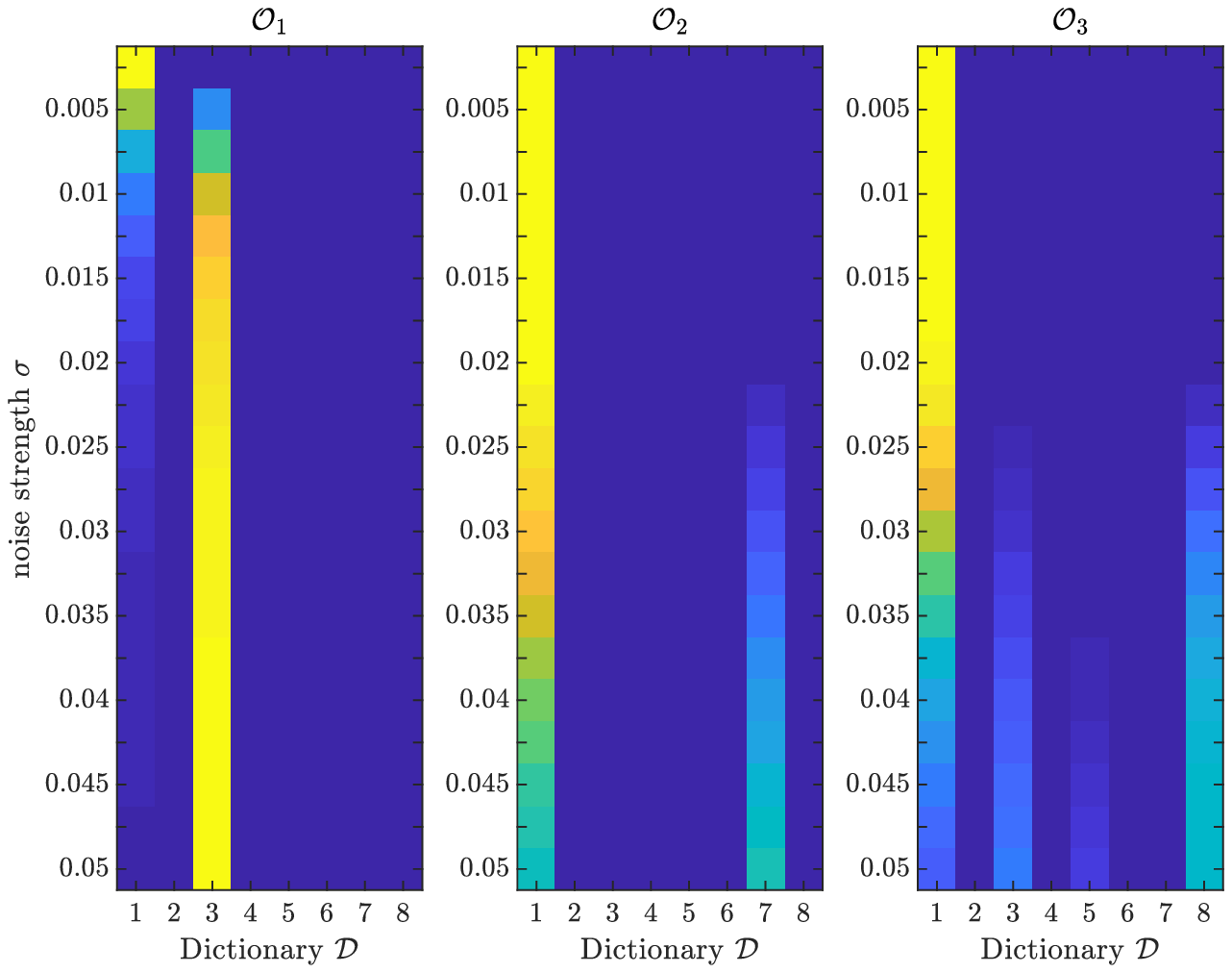} & 
		\fbox{\includegraphics[scale=0.61]{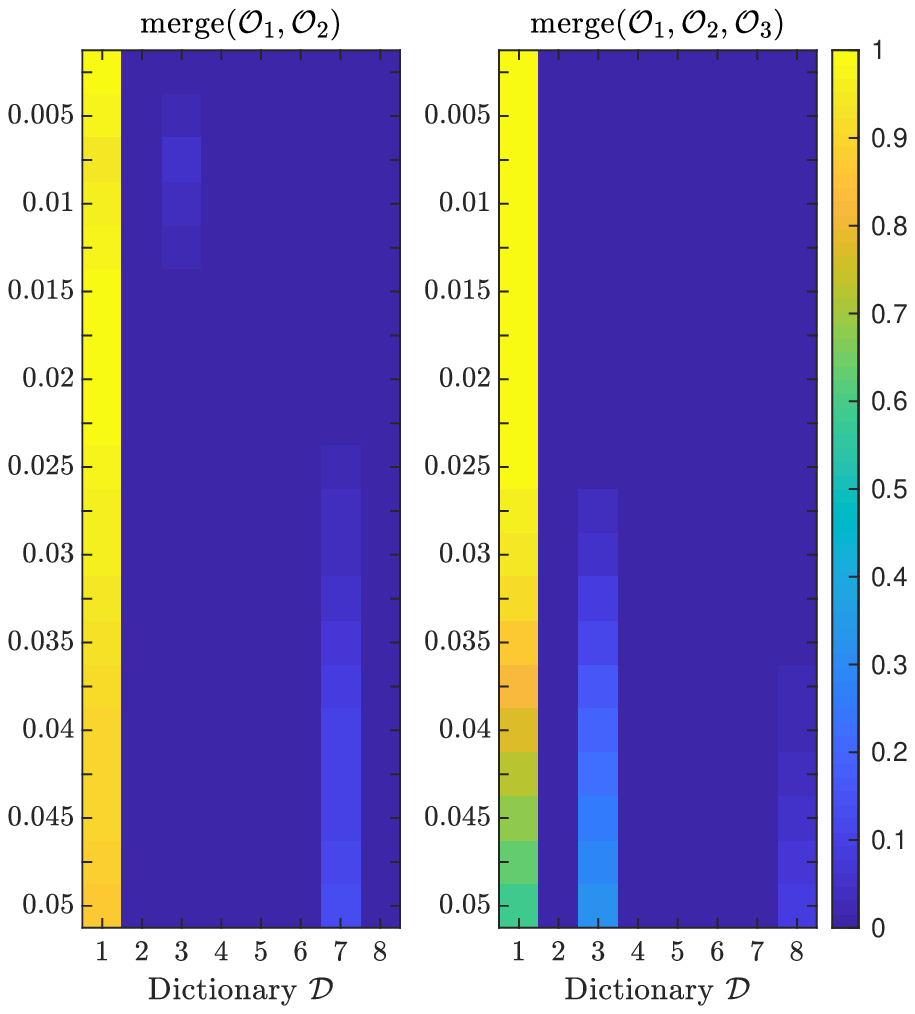}}
	\end{tabular}
	\caption{Circle (1)}
	\label{fig:recoCircleMP}
\end{figure}

\begin{figure}[H]
	\centering
	\hspace*{-5mm}
	\begin{tabular}{cc}
		\includegraphics[scale=0.6]{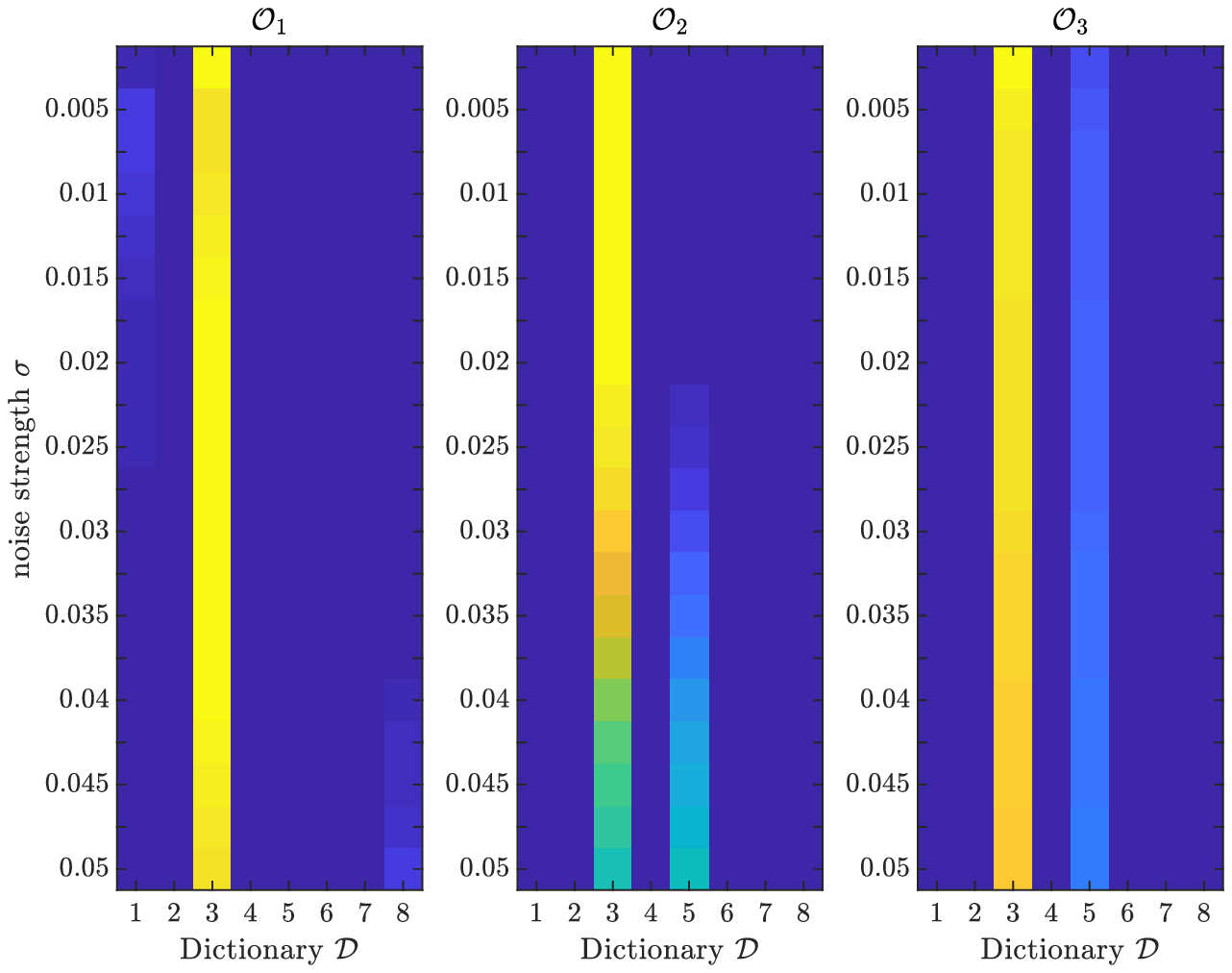} & 
		\fbox{\includegraphics[scale=0.61]{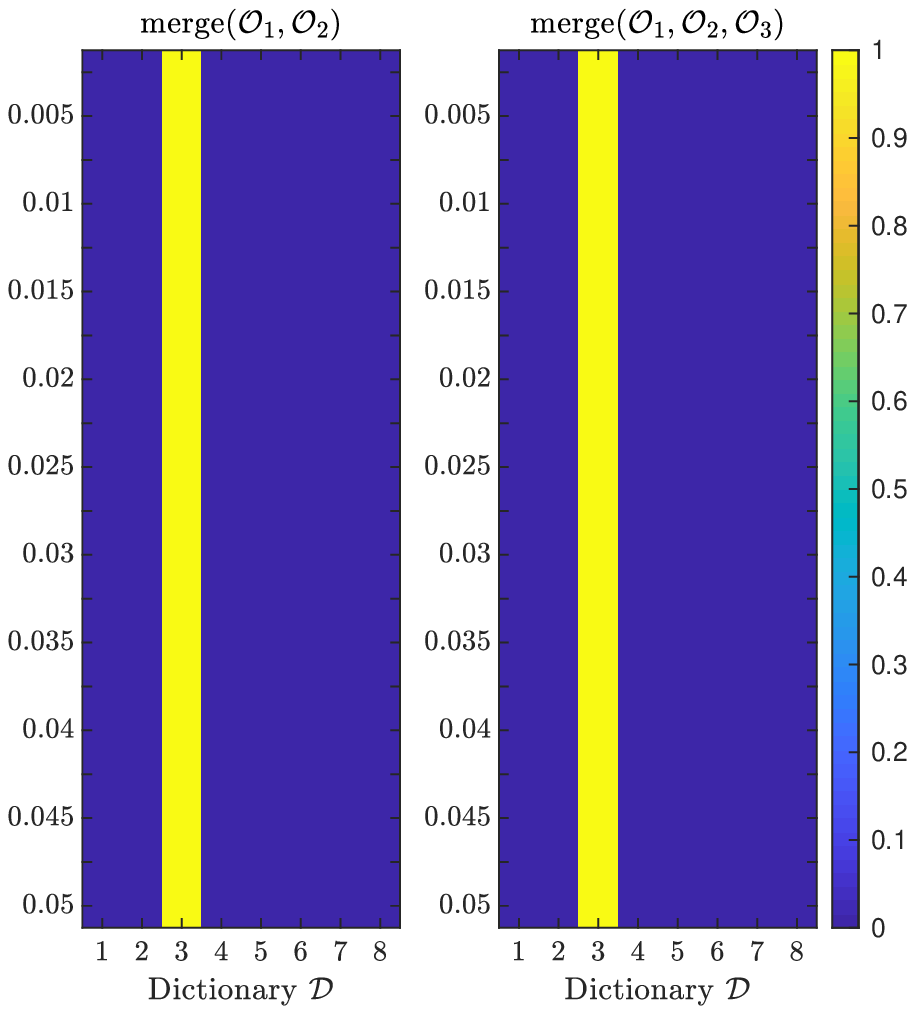}}
	\end{tabular}
	\caption{Triangle (3)}
	\label{fig:recoTriangleMP}
\end{figure}

\begin{figure}[H]
	\centering
	\hspace*{-5mm}
	\begin{tabular}{cc}
		\includegraphics[scale=0.6]{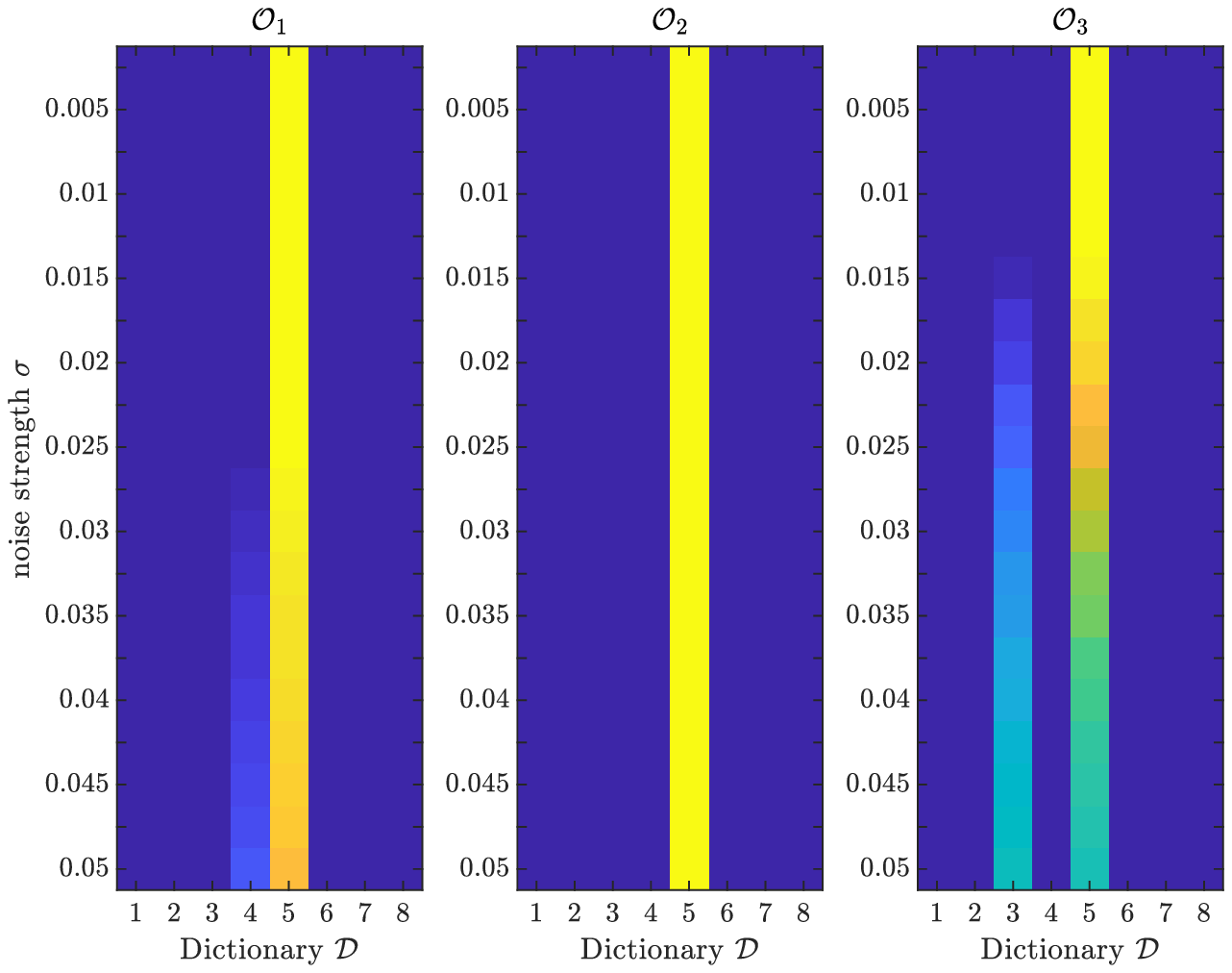} & 
		\fbox{\includegraphics[scale=0.61]{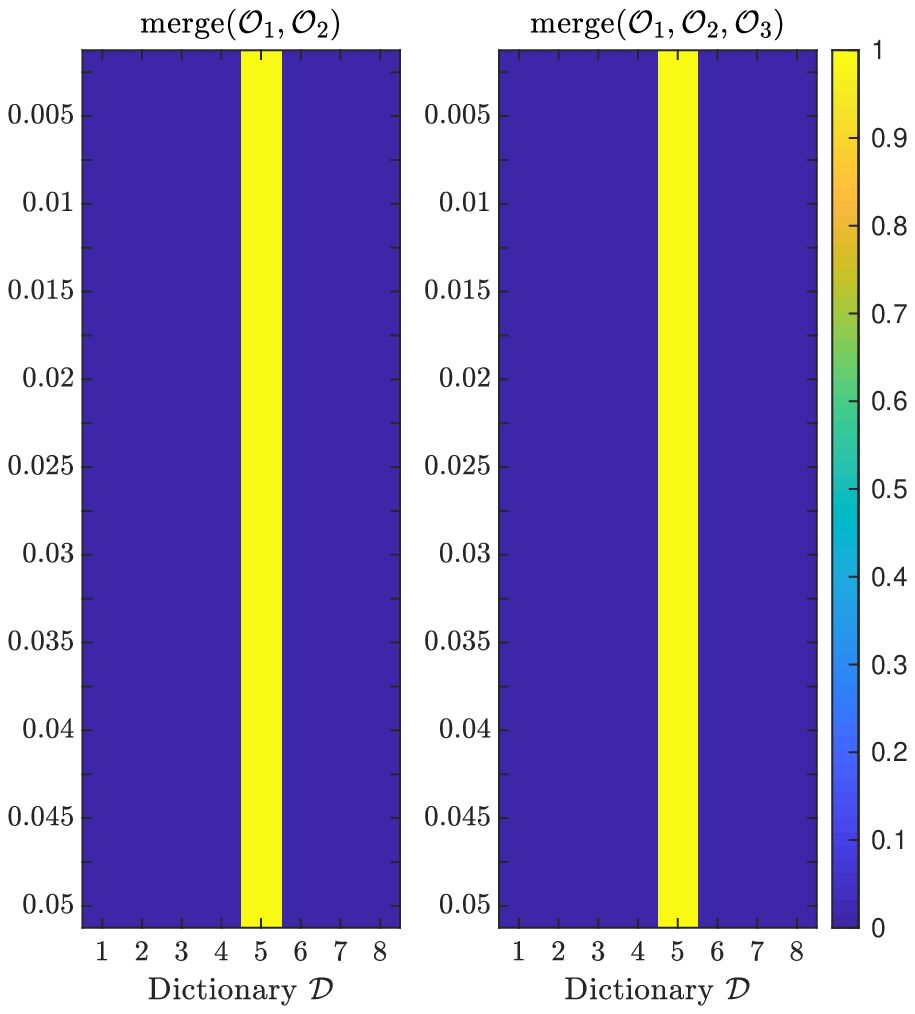}}
	\end{tabular}
	\caption{Curved Triangle (5)}
	\label{fig:recoTriaShieldMP}
\end{figure}

\begin{table}[H]
	\centering
%	\begin{tabular}{l|>{\columncolor[gray]{0.9}}cccccccc}
%		\toprule
%		&$\boxed{\text{1}}$&2&3&4&5&6&7&8\\
%		\midrule
%		$\mathsf{m}_{1}$ &\rule{7pt}{0pt}0& \rule{7pt}{0pt}0& \rule{7pt}{0pt}0& \rule{7pt}{0pt}0& \rule{7pt}{0pt}0& \rule{7pt}{0pt}0& \rule{7pt}{0pt}0& \rule{7pt}{0pt}0\\
%		$\mathsf{m}_{2}$ &\rule{7pt}{0pt}0& \rule{7pt}{0pt}0& \rule{7pt}{0pt}0 & \rule{7pt}{0pt}0& \rule{7pt}{0pt}0& \rule{7pt}{0pt}0& \rule{7pt}{0pt}0& \rule{7pt}{0pt}0\\
%		$\mathsf{m}_{3}$ &\rule{7pt}{0pt}0 & \rule{7pt}{0pt}0& \rule{7pt}{0pt}0& \rule{7pt}{0pt}0& \rule{7pt}{0pt}0 & \rule{7pt}{0pt}0& \rule{7pt}{0pt}0& \rule{7pt}{0pt}0\\ 
%		\hline
%		$\mathsf{m}_{12}$ &\rule{7pt}{0pt}0& \rule{7pt}{0pt}0& \rule{7pt}{0pt}0 & \rule{7pt}{0pt}0& \rule{7pt}{0pt}0& \rule{7pt}{0pt}0& \rule{7pt}{0pt}0& \rule{7pt}{0pt}0\\ 
%		$\mathsf{m}_{123}$&\rule{7pt}{0pt}0& \rule{7pt}{0pt}0& \rule{7pt}{0pt}0& \rule{7pt}{0pt}0& \rule{7pt}{0pt}0 & \rule{7pt}{0pt}0& \rule{7pt}{0pt}0& \rule{7pt}{0pt}0\\ 
%		\bottomrule
%	\end{tabular}

	\begin{tabular}{l|>{\columncolor[gray]{0.9}}cccccccc}
	\toprule
	&$\boxed{\text{1}}$&2&3&4&5&6&7&8\\
	\midrule
	$\mathsf{m}_{1}$ &\rule{7pt}{0pt}0.0106& \rule{7pt}{0pt}0& \rule{7pt}{0pt}0.9885& \rule{7pt}{0pt}0& \rule{7pt}{0pt}0& \rule{7pt}{0pt}0& \rule{7pt}{0pt}0& \rule{7pt}{0pt}0.0009\\
	$\mathsf{m}_{2}$ &\rule{7pt}{0pt}0.4889& \rule{7pt}{0pt}0& \rule{7pt}{0pt}0 & \rule{7pt}{0pt}0& \rule{7pt}{0pt}0& \rule{7pt}{0pt}0& \rule{7pt}{0pt}0.5111& \rule{7pt}{0pt}0\\
	$\mathsf{m}_{3}$ &\rule{7pt}{0pt}0.1784 & \rule{7pt}{0pt}0& \rule{7pt}{0pt}0.2565& \rule{7pt}{0pt}0.0138& \rule{7pt}{0pt}0.0817 & \rule{7pt}{0pt}0& \rule{7pt}{0pt}0.0021& \rule{7pt}{0pt}0.4675\\ 
	\hline
	$\mathsf{m}_{12}$ &\rule{7pt}{0pt}0.8729& \rule{7pt}{0pt}0& \rule{7pt}{0pt}0 & \rule{7pt}{0pt}0& \rule{7pt}{0pt}0& \rule{7pt}{0pt}0& \rule{7pt}{0pt}0.1270& \rule{7pt}{0pt}0.0001\\ 
	$\mathsf{m}_{123}$&\rule{7pt}{0pt}0.5914& \rule{7pt}{0pt}0& \rule{7pt}{0pt}0.3185& \rule{7pt}{0pt}0.0001& \rule{7pt}{0pt}0.0016& \rule{7pt}{0pt}0& \rule{7pt}{0pt}0.0067& \rule{7pt}{0pt}0.0817\\ 
	\bottomrule
\end{tabular}
	
	\medskip
	
	\caption{Frequency table for the identification of the Circle (1) with the strongest noise, i.e., $\sigma = 0.05$. Each row displays the relative frequencies for all the elements of the dictionary corresponding to different BBAs.}
	\label{frequency-table-reco-MP-1}
\end{table}

\begin{table}[H]
	\centering
	\begin{tabular}{l|cc>{\columncolor[gray]{0.9}}cccccc}
	\toprule
	&1&2&$\boxed{\text{3}}$&4&5&6&7&8\\
	\midrule
	$\mathsf{m}_{1}$ &    0.0033 & \rule{7pt}{0pt}0& \rule{7pt}{0pt}0.9121& \rule{7pt}{0pt}0& \rule{7pt}{0pt}0& \rule{7pt}{0pt}0& \rule{7pt}{0pt}0& \rule{7pt}{0pt}0.0846\\
	$\mathsf{m}_{2}$&0& \rule{7pt}{0pt}0& \rule{7pt}{0pt}0.5027 & \rule{7pt}{0pt}0& \rule{7pt}{0pt}0.4973& \rule{7pt}{0pt}0& \rule{7pt}{0pt}0& \rule{7pt}{0pt}0\\
	$\mathsf{m}_{3}$ &0 & \rule{7pt}{0pt}0& \rule{7pt}{0pt}0.7652 & \rule{7pt}{0pt}0& \rule{7pt}{0pt}0.2346 & \rule{7pt}{0pt}0& \rule{7pt}{0pt} 0.0002 & \rule{7pt}{0pt}0\\ 
	\hline
	$\mathsf{m}_{12}$ & 0& \rule{7pt}{0pt}0& \rule{7pt}{0pt}0.9933     & \rule{7pt}{0pt}0& \rule{7pt}{0pt}0.0046         & \rule{7pt}{0pt}0& \rule{7pt}{0pt}0& \rule{7pt}{0pt}0.0021\\ 
	$\mathsf{m}_{123}$&0& \rule{7pt}{0pt}0& \rule{7pt}{0pt}0.9910 & \rule{7pt}{0pt}0& \rule{7pt}{0pt}0.0090& \rule{7pt}{0pt}0& \rule{7pt}{0pt}0& \rule{7pt}{0pt}0\\ 
	\bottomrule
\end{tabular}
	
	\medskip
	
	\caption{Frequency table for the identification of the Triangle (3) with the strongest noise, i.e., $\sigma = 0.05$.}
	\label{frequency-table-reco-MP-3}
\end{table}

\begin{table}[H]
	\centering
	\begin{tabular}{l|cccc>{\columncolor[gray]{0.9}}cccc}
	\toprule
	&1&2&3&4&$\boxed{\text{5}}$&6&7&8\\
	\midrule
	$\mathsf{m}_{1}$ &\rule{7pt}{0pt}0& \rule{7pt}{0pt}0& \rule{7pt}{0pt}0& \rule{7pt}{0pt}0.1476& \rule{7pt}{0pt}0.8512& \rule{7pt}{0pt}0& \rule{7pt}{0pt}0& \rule{7pt}{0pt}0.0012\\
	$\mathsf{m}_{2}$ & \rule{7pt}{0pt}0& \rule{7pt}{0pt}0& \rule{7pt}{0pt}0.0013 & \rule{7pt}{0pt}0& \rule{7pt}{0pt}0.9987& \rule{7pt}{0pt}0& \rule{7pt}{0pt}0& \rule{7pt}{0pt}0\\
	$\mathsf{m}_{3}$ & \rule{7pt}{0pt}0 & \rule{7pt}{0pt}0& \rule{7pt}{0pt}0.4855& \rule{7pt}{0pt}0& \rule{7pt}{0pt}0.5145 & \rule{7pt}{0pt}0& \rule{7pt}{0pt}0& \rule{7pt}{0pt}0\\ 
	\hline
	$\mathsf{m}_{12}$ & \rule{7pt}{0pt}0& \rule{7pt}{0pt}0& \rule{7pt}{0pt}0 & \rule{7pt}{0pt}0.0008& \rule{7pt}{0pt}0.9992& \rule{7pt}{0pt}0& \rule{7pt}{0pt}0& \rule{7pt}{0pt}0\\ 
	$\mathsf{m}_{123}$& \rule{7pt}{0pt}0& \rule{7pt}{0pt}0& \rule{7pt}{0pt}0.0001& \rule{7pt}{0pt}0& \rule{7pt}{0pt}0.9999 & \rule{7pt}{0pt}0& \rule{7pt}{0pt}0& \rule{7pt}{0pt}0\\ 
	\bottomrule
\end{tabular}
	
	\medskip
	
	\caption{Frequency table for the identification of the Curved Triangle (5) with the strongest noise, i.e., $\sigma = 0.05$.}
	\label{frequency-table-reco-MP-5}
\end{table}

\FloatBarrier

While the g-inverse $\mathfrak{L}$ lacks of the property of being a least-square g-inverse, the Moore-Penrose $\bfL^\dagger$ provides the solution to the minimization problem \eqref{eq:lst_squares_vec}. Therefore, it is not surprising that the classification rates obtained by using the latter are generally better than the ones resulting from using $\mathfrak{L}$.

\section{Concluding remarks} 
In this paper, we have presented a dictionary-matching approach for classification in
electro-sensing that takes advantage of measurements at different length-scales.
We have performed a careful analysis of the acquisition operator that was not available yet. In particular, by exploiting its peculiar block Kronecker form, we have studied its rank and established a length-scale dependent estimate on the reconstruction error. We have also discussed to which extent the limited-view configuration impacts on predicting the committed error.

\vspace{1mm}	

\section*{Acknowledgments}
The authors gratefully acknowledge Prof. H. Ammari for his guidance. During the preparation of this work, the authors were financially supported by a Swiss National Science Foundation grant (number 200021-172483).  

\vspace{1mm}

\appendix

\section{Generalized polarization tensors and boundary layer potentials}
\label{apx:GPTs}

We briefly summarise some fundamental concepts that are essential for understanding the problem. For further references see \cite{AMM2013,Ammari2013,Ammari11652,SCAPIN20191872}.

Let $\Omega$ be a simply-connected bounded domain. We assume
$\Omega \in C^{2,\alpha}$. 
%Given an arbitrary function $u$ defined on $\mathbb{R}^2\setminus \partial \Omega$ and $x \in \partial \Omega$, we define
%\[u(x)|_{\pm} := \lim_{t\to 0} u(x\pm t\nu(x)), \]
%\[ \left. \frac{\partial u}{\partial \nu}(x) \right|_{\pm} := \lim_{t \to 0} \nabla u (x \pm t \nu(x) ) \cdot \nu(x). \]

\begin{definition}
	Denote by $\Gamma$ the fundamental solution of the Laplacian in $\mathbb{R}^2$, i.e.,
	\[
	\Gamma(x-y) := \frac{1}{2 \pi} \log |x-y|, \qquad x,y \in \mathbb{R}^2.
	\]
\end{definition}

\begin{definition}
	For any $\phi \in L^2(\partial \Omega)$, the single- and double-layer potentials on $\Omega$ are given by the following formulas:
	\[\mathcal{S}_\Omega [\phi] (x) := \int_{\partial \Omega} \Gamma(x,y) \phi(y)  \text{ d} \sigma_y , \qquad x \in \mathbb{R}^2, \]
	\[\mathcal{D}_\Omega [\phi] (x) := \int_{\partial \Omega} \frac{\partial \Gamma}{\partial \nu_y}(x,y) \phi(y)  \text{ d} \sigma_y , \qquad x \in \mathbb{R}^2\setminus \partial \Omega. \]
\end{definition}

Recall that for $\phi \in L^2(\Omega)$, the functions $\mathcal{S}_\Omega$ and $\mathcal{D}_\Omega$ are harmonic functions in $\mathbb{R}^2 \setminus \partial \Omega$.

\begin{definition}
	The operator $\mathcal{K}_\Omega$ and its $L^2$-adjoint $\mathcal{K}^*_\Omega$ are given by the following formulas:
	\[ \mathcal{K}_\Omega [\phi] (x) := \frac{1}{2\pi} \text{p.v.} \int_{\partial \Omega} \frac{(y-x) \cdot \nu(y)}{|x-y|^2} \phi(y)  \text{ d} \sigma_y , \qquad x \in \partial \Omega, \]
	\[ \mathcal{K}^*_\Omega [\phi] (x) := \frac{1}{2\pi} \text{p.v.} \int_{\partial \Omega} \frac{(x-y) \cdot \nu(x)}{|x-y|^2} \phi(y)  \text{ d} \sigma_y , \qquad x \in \partial \Omega\] where p.v. stands for the Cauchy principal value.
\end{definition}

$\mathcal{K}^*_\Omega$ is also known as the Neumann-Poincar\'e operator.

We introduce the generalized polarization tensor (GPT).

\begin{definition} Let $\alpha,\beta \in \mathbb{N}^2$ be multi-indices. We define the generalized polarization tensor associated with the domain $\Omega$ and the contrast $\lambda$ by
	\[M_{\alpha \beta} (\lambda, \Omega):= \int_{\partial \Omega} (\lambda {I} - \mathcal{K}_\Omega^*)^{-1} \left[\frac{\partial y^\alpha}{\partial \nu}\right] y^\beta \text{ d}\sigma_y,\]
	where ${I}$ is the identity operator.
\end{definition}
We can also define the contracted generalized polarization tensor (CGPT) as follows.

\begin{definition} \label{def:CGPTs}
	Let $m,n \in \mathbb{N}$. We define the contracted generalized polarization tensors by
	\[M^{cc}_{mn} = \sum_{|\alpha|=m} \sum_{|\beta|=n} a^m_\alpha a^n_\beta M_{\alpha \beta} ,\]
	\[M^{cs}_{mn} = \sum_{|\alpha|=m} \sum_{|\beta|=n} a^m_\alpha b^n_\beta M_{\alpha \beta} ,\]
	\[M^{sc}_{mn} = \sum_{|\alpha|=m} \sum_{|\beta|=n} b^m_\alpha a^n_\beta M_{\alpha \beta} ,\]
	\[ M^{ss}_{mn} = \sum_{|\alpha|=m} \sum_{|\beta|=n} b^m_\alpha b^n_\beta M_{\alpha \beta} ,\]
	where the real numbers $a^m_\alpha$ and $b^m_\beta$ are defined by the following relation
	\[
	(x_1+ix_2)^m = \sum_{|\alpha|=m} a^m_\alpha x^\alpha + \sum_{|\beta|=m} b^m_\beta x^\beta.
	\]
\end{definition}

\section{Kronecker products and generalized inverses}

Let us denote by $\Mat{m}{n}$ the space of $m \times n$ matrices.

\begin{definition}[vec operator] Given a matrix $\bfX = \begin{bmatrix}
	\mathbf{x}^1 & \mathbf{x}^2 & \dots & \mathbf{x}^m \end{bmatrix} \in \Mat{k}{m}$, define the vectorization operator $\vec{\cdot} : \Mat{k}{m} \rightarrow \Mat{km}{1}$ as follows
	\begin{equation} \label{eq:vec} 
	\vec{\bfX} = \begin{bmatrix} \mathbf{x}^1 \\ \vdots \\ \mathbf{x}^m
	\end{bmatrix} \,\, \in \Mat{km}{1} .
	\end{equation}
\end{definition}

%The operator Vec
%converts a k × m matrix A into a vector Vec(A) ∈ R km by stacking the columns
%one underneath the other.

\begin{definition}[Kronecker product] Given $\bfX = (x_{ij}) \in \Mat{m}{n}$ and $\bfY = (y_{ij}) \in \Mat{p}{q}$, define
	\begin{equation} \label{eq:kron} 
	\bfX \kron \bfY = \begin{bmatrix}
	x_{11} \bfY & x_{12} \bfY &  \dots & x_{1n} \bfY \\
	x_{21} \bfY & x_{22} \bfY &  \dots & x_{2n} \bfY \\
	\vdots & \vdots  &  \ddots & \vdots \\
	x_{m1} \bfY & x_{m2} \bfY &  \dots & x_{mn} \bfY \\
	\end{bmatrix} \,\, \in \Mat{mp}{nq} .
	\end{equation}
\end{definition}

%\begin{definition}[Hadamard product] Given $\bfX = (x_{ij}) , \bfY = (y_{ij}) \in \Mat{m}{n}$, define
%	\begin{equation} \label{eq:kron} 
%	\bfX \hadam \bfY = \begin{pmatrix}
%	x_{11} y_{11} & x_{12} y_{12} &  \dots & x_{1n} y_{1n} \\
%	x_{21} y_{21} & x_{22} y_{22} &  \dots & x_{2n} y_{2n} \\
%	\vdots & \vdots  &  \ddots & \vdots \\
%	x_{m1} y_{m1} & x_{m2} y_{m2} &  \dots & x_{mn} y_{mn} \\
%	\end{pmatrix} \,\, \in \Mat{m}{n},
%	\end{equation}
%	i.e., the element-wise product of $\bfX$ and $\bfY$.
%\end{definition}

For all matrices $\bfA$, $\bfB$ and $\bfC$ such that the product $\bfA \bfB \bfC$ is well defined we have
\begin{equation} \label{eq:kron_vec} 
\vec{\bfA \bfB \bfC} = (\bfC^\top \kron \bfA ) \vec{\bfB}.
\end{equation}

We introduce a generalized Kronecker product.
\begin{definition}[Generalized Kronecker product \cite{REGALIA1989,Marco2018}] \label{def:kronecker_mitra} Given a matrix $\bfX$ ($M \times l$) and a set of $M$ matrices $\bfY$ ($N\times r$)  we define the  matrix $\bfX \kron \{\bfY_\ell\}$ ($M N \times lr$ )  as
	\begin{equation}
	\bfX \kron \{\bfY_\ell\} \deff \begin{bmatrix} \bfX_{1,:} \kron \bfY_1 \\ \vdots \\ \bfX_{M,:} \kron \bfY_M  \end{bmatrix} .
	\end{equation}
	
\end{definition}

Notice that in \cite{REGALIA1989} the factors are swapped compared to this definition. The reason is that, unlike us, they considered a \emph{left} Kronecker product.

\medskip

\begin{definition}[Generalized column-wise Kronecker product] \label{def:col_gkp} Given a matrix $\bfX$ ($M \times l$) and a set of $M$ matrices $\bfY_\ell$ ($N\times r$)  we define the  matrix $\bfX \kron_C \{\bfY_\ell\}$ ($M N \times lr$ )  as
	\begin{equation}
	\bfX \kron_C \{\bfY_\ell\} \deff \begin{bmatrix} \bfX_{:,1} \kron \bfY_1 & \dots & \bfX_{:,l} \kron \bfY_M  \end{bmatrix} .
	\end{equation}	
\end{definition}

\begin{definition}[Moore-Penrose inverse] The Moore-Penrose inverse of the matrix $\bfM$ ($m \times p$) is the unique matrix $\bfG$ ($p \times m$) satisfying the four \emph{Penrose conditions}:	
	\begin{enumerate}
		\item $\bfM \bfG \bfM = \bfM$,
		\item $ \bfG \bfM \bfG = \bfG$,
		\item $(\bfM \bfG)^\top = \bfM \bfG$,
		\item $(\bfG \bfM)^\top = \bfG \bfM$.		
	\end{enumerate}	
\end{definition}
A matrix $\bfG$ satisfying (1) is called \emph{a generalized inverse} or \emph{g-inverse}. 

A matrix $\bfG$ satisfying (1) and (2) is called \emph{reflexive g-inverse}. 

A matrix $\bfG$ satisfying (1) and (3) is called \emph{least-square g-inverse}. 

A matrix $\bfG$ satisfying (1) and (4) is called \emph{minimum norm  g-inverse}.

\section{Reminder in the asymptotic expansion \eqref{eq:q_sr}} \label{apx:reminder}

Let $\varepsilon = \delta/\rho$ be the length-scale. By definition, the truncation error at the receptor $x_r^{(s)}$ is given by
\[ \bfE_{rs} \deff u^{(s)}(x_r^{(s)}) - H^{(s)}(x_r^{(s)}) - \mc{L}^{(s)}(\bbM^{(K)}).\]
The asymptotic behavior of $\bfE_{sr}$ in the far-field regime is assessed by the following proposition. 

\begin{prop}  For $\delta , \varepsilon \ll 1$ we have the following asymptotic behavior 
	\begin{equation}\label{eq:errorEsr} |\bfE_{rs}| = O (\varepsilon^{K+2}). \end{equation}
\end{prop}

\begin{proof} 
	Hereinafter, we simplify our notation by fixing the index position $s$, i.e., $x^{(s)}_r = x_r$, $H^{(s)}=H$, $\bfE_{rs}=E_r$. 
	
	Let us define $H_K(x) = \displaystyle \sum_{|\alpha| = 0}^{K} \frac{1}{\alpha!} \p^{\alpha} H(z) (x - z)^\alpha $.
	The truncation error can be expressed as
	\begin{equation}\label{eq:split} 
	E_r = \int_{\p D} \Gamma_K(x_r-y) (\lambda I - \mc{K}_D^* )^{-1} [ \p_\nu H - \p_\nu H_K ] \mbox{ d} s_y + \int_{\p D} ( \Gamma - \Gamma_K )(x_r-y)  (\lambda I - \mc{K}_D^* )^{-1} [ \p_\nu H ] \mbox{ d} s_y .
	\end{equation}
	%	\begin{equation}
	%	\begin{split} u(x_r) - H(x_r) &= \int_{\p D} \Gamma(x_r-y) (\lambda I - \mc{K}_D^* )^{-1} [ \p_\nu H ] \mbox{ d} s_y
	%	\\ & = \int_{\p D} \Gamma(x_r-y) (\lambda I - \mc{K}_D^* )^{-1} [ \p_\nu H - \p_\nu H_K ] \mbox{ d} s_y  
	%	\\ & + \int_{\p D} \Gamma(x_r-y) (\lambda I - \mc{K}_D^* )^{-1} [ \p_\nu H_K ] \mbox{ d} s_y 
	%	\\ & = \int_{\p D} \Gamma(x_r-y) (\lambda I - \mc{K}_D^* )^{-1} [ \p_\nu H - \p_\nu H_K ] \mbox{ d} s_y \\ & + \int_{\p D} ( \Gamma - \Gamma_K )(x_r-y)  (\lambda I - \mc{K}_D^* )^{-1} [ \p_\nu H_K ] \mbox{ d} s_y 
	%	\\ & +  \int_{\p D} \Gamma_K(x_r-y) (\lambda I - \mc{K}_D^* )^{-1} [ \p_\nu H_K ] \mbox{ d} s_y .\end{split}
	%	\end{equation}
	We want to estimate each term separately.
	
	\medskip
	
	Denote the first term:
	
	\begin{equation}E_r^{(1)} :=  \int_{\p D}  \Gamma_K(y;x_r,z) (\lambda {I} - \mc{K}_D^* )^{-1} \left [ \frac{\p H}{\p \nu} -  \frac{\p H_K}{\p \nu} \right ](y) \mbox{ d} s_y  .
	\end{equation}
	
	We have
	
	\begin{equation*} \begin{split} \left  |E_r^{(1)} \right | & \le \left | \int_{\p D}  \Gamma_K(y;x_r,z) (\lambda {I} - \mc{K}_D^* )^{-1} \left [ \frac{\p H}{\p \nu} -  \frac{\p H_K}{\p \nu} \right ](y) \mbox{ d} s_y \right | 
	\\& \le \sup_{y \in \p D} |\Gamma_K(y;x_r,z)|  \int_{\p D} \left |(\lambda {I} - \mc{K}_D^* )^{-1} \left [ \frac{\p H}{\p \nu} -  \frac{\p H_K}{\p \nu} \right ](y)\right| \mbox{ d} s_y 
	\\& \le \| \Gamma_K( \, \cdot\, ; x_r,z) \|_{L^{\infty}(\p D)}  \left \| (\lambda {I} - \mc{K}_D^* )^{-1} \left [ \frac{\p H}{\p \nu} -  \frac{\p H_K}{\p \nu} \right ]\right \|_{L^2(\p D)}  |\p D|^{1/2} 
	\\&\le C  \| \Gamma_K( \, \cdot\, ; x_r,z) \|_{L^{\infty}(\p D)} \left \| \frac{\p H}{\p \nu} -  \frac{\p H_K}{\p \nu} \right \|_{L^2(\p D)} |\p D|^{1/2}  
	\\& \le C' \frac{\delta}{\rho}  \left \| \frac{\p H}{\p \nu} -  \frac{\p H_K}{\p \nu} \right \|_{L^2(\p D)} |\p D|^{1/2}  
	\\& \le C' \frac{\delta^2}{\rho}  \left \| \frac{\p H}{\p \nu} -  \frac{\p H_K}{\p \nu} \right \|_{L^{\infty}(\p D)} |\p B| .  
	\end{split}
	\end{equation*}

	Now we estimate the term 
	\begin{equation}
	\left \| \frac{\p H}{\p \nu} -  \frac{\p H_K}{\p \nu} \right \|_{L^{\infty}(\p D)}.	\end{equation}
	Recall the integral form of the reminder of Taylor's formula:
	\begin{equation*}
	(H - H_K ) (y) =  \sum_{|\alpha| = K+1}^{} \frac{1}{\alpha!} \int_0^1 (1-t)^K \p^{\alpha} H(ty) \,\mbox{d}t\; y^\alpha ,
	\end{equation*}
	then
	\begin{equation*}
	\left ( \frac{\p H}{\p \nu} -  \frac{\p H_K}{\p \nu} \right ) (y) =  \sum_{|\alpha| = K+1}^{} \frac{1}{\alpha!}  \left [ \int_0^1 (1-t)^K \p_\nu \p^{\alpha} H(ty) \,\mbox{d}t\; y^\alpha + \int_0^1 (1-t)^K \p^{\alpha} H(ty) \,\mbox{d}t\; \p_\nu y^\alpha  \right ].
	\end{equation*}
	It is immediate to show that
	\begin{equation*} \begin{split}
	\left | \left ( \frac{\p H}{\p \nu} -  \frac{\p H_K}{\p \nu} \right ) \biggr |_{\p D} \right | &
	\le  C \left ( \delta^{K+1} \sum_{|\alpha| = K+1}^{} \frac{1}{\alpha!} \| \nabla \p^{\alpha} H \|_{C^0(\overline{D})} \; +  \delta^{K} \sum_{|\alpha| = K+1}^{} \frac{1}{\alpha!} \| \p^{\alpha} H \|_{C^0(\overline{D})}  \right ) .
	\end{split}
	\end{equation*}
	Assume $\rho \ge \mbox{dist} (\p \Omega , 0) > 1$, $\delta \ll 1$. By using formulas for the derivatives of $H$ we get
	\begin{equation} \label{eq:dnu_dnuH} \left  \|  \frac{\p H}{\p \nu} -  \frac{\p H_K}{\p \nu}   \right \|_{L^\infty(\p D)} \le \left ( \frac{C_1'}{\rho^{K+3}} + \frac{C_2'}{\rho^{K+2}}\right ) \delta^{K+1}  + \left ( \frac{C_1}{\rho^{K+2}} + \frac{C_2}{\rho^{K+1}}\right ) \delta^K, \end{equation}
	hence
	\[ \left  \|  \frac{\p H}{\p \nu} -  \frac{\p H_K}{\p \nu}   \right \|_{L^\infty(\p D)} = O \left ( \frac{\delta^{K}}{\rho^{K+1}} \right ) .\]
	Therefore
	\begin{equation*}
	\left  | E_r^{(1)} \right |  \le  C' \frac{\delta}{\rho}  \left \| \frac{\p H}{\p \nu} -  \frac{\p H_K}{\p \nu} \right \|_{L^2(\p D)} |\p D|^{1/2} \le C'' \frac{\delta^{K+2}}{\rho^{K+2}} = O ( \varepsilon^{K+2} ) .
	\end{equation*}

	Denote the second term in \eqref{eq:split}:
	\begin{equation}E_r^{(2)} :=  \int_{\p D} ( \Gamma - \Gamma_K )(x_r-y)  (\lambda I - \mc{K}_D^* )^{-1} \left [ \frac{\p H}{\p \nu} \right ](y)  \mbox{ d} s_y  .
	\end{equation}
	We have
	\begin{equation*}
	\begin{split} \left|E_r^{(2)}\right| & = \left | \int_{\p D} ( \Gamma - \Gamma_K )(x_r-y)  (\lambda I - \mc{K}_D^* )^{-1} \left [ \frac{\p H_K}{\p \nu} \right ](y) \mbox{ d} s_y \right |
	\\ & \le  \| \Gamma(x_r- \, \cdot\,) - \Gamma_K(x_r-\cdot)\|_{L^{\infty}(\p D)} \int_{\p D}  \left | (\lambda I - \mc{K}_D^* )^{-1} \left [ \frac{\p H}{\p \nu} \right ](y) \right |  \mbox{ d} s_y 
	\\ & \le  \| \Gamma(x_r- \, \cdot\,) - \Gamma_K(x_r-\cdot)\|_{L^{\infty}(\p D)}  \left \| (\lambda I - \mc{K}_D^* )^{-1} \left [ \frac{\p H}{\p \nu} \right ] \right \|_{L^2(\p D)}  |\p D|^{1/2}
	\\ & \le  C \| \Gamma(x_r- \, \cdot\,) - \Gamma_K(x_r-\cdot)\|_{L^{\infty}(\p D)} \left \| \frac{\p H}{\p \nu} \right \|_{L^2(\p D)} |\p D|^{1/2}
	\\ & \le  C \| \Gamma(x_r- \, \cdot\,) - \Gamma_K(x_r-\cdot)\|_{L^{\infty}(\p D)} \left \| \frac{\p H}{\p \nu} \right \|_{L^\infty(\p D)} |\p D| .
	\end{split}
	\end{equation*}
	
	Since $\| \Gamma(x_r- \, \cdot\,) - \Gamma_K(\,\cdot\,; x_r,z)\|_{L^{\infty}(\p D)} \le C \left ( \frac{\delta}{\rho} \right )^{K+1}$, see \cite{Ammari2014}, and 
	\begin{equation}\left \| \frac{\p H}{\p \nu} \right \|_{L^\infty(\p D)} \le C \rho^{-1},
	\end{equation}
	we have
	\begin{equation*} 
	\left|E_r^{(2)}\right|  = O (\varepsilon^{K+2}).
	\end{equation*}
	
\end{proof}

%Now we derive asymptotic behavior of the reminder in the near-field regime, that is, when $\varepsilon \ll 1$ but $R \ll 1$. In this regime the reminder deteriorates faster than $O(R^{-(K+2)})$ (near the target) as $R$ decreases. This effect is due to the pure dipole term.
%
%
%\begin{prop} \label{prop:dist_dependent} For $\varepsilon \ll 1$ and $R \ll 1$ we have the following asymptotic behavior
%	\begin{equation} |E_{sr}| = O \left ({\varepsilon^{K+2}} R^{-1} \right ). \end{equation}
%\end{prop}
%
%
%\begin{proof} Estimate 	\eqref{eq:dnu_dnuH} yields 
%	
%	\[ \left  \|  \frac{\p H}{\p \nu} -  \frac{\p H_K}{\p \nu}   \right \|_{L^\infty(\p D)} = O \left ( \frac{\delta^{K}}{R^{K+2}} \right ) .\] 
%	
%	and therefore 
%	\begin{equation*}
%	\left  | (\star) \right |  =   \frac{\delta^{K+2}}{R^{K+3}} = O \left ( {\varepsilon^{K+2}} R^{-1}\right ) .
%	\end{equation*}	
%	
%	On the other hand, 		
%	\begin{equation*} \left \| \frac{\p H}{\p \nu} \right \|_{L^\infty(\p D)} \le C R^{-2},	\end{equation*}
%	from which we get
%	\begin{equation*} 
%	|(\star \star)|  = O (\varepsilon^{K+2} R^{-1}).
%	\end{equation*}	
%\end{proof}
%
%
%
%Proposition \ref{prop:dist_dependent} states that the asymptotic behavior of the truncation error in the expansion \eqref{eq:q_sr} deteriorates like $R^{-(K+3)}$ as the fish gets closer to the target. We perform numerical experiments to validate this predicted asymptotic behavior.

\section{Technical estimates}

\subsection{Uniqueness results}
\label{apx:uniqueness}

In this section, we show that the matrix of receptors $\bfG^{(s)}$ is full column rank. We begin by observing that $\bfG^{(s)}$ is closely related to a special $N_r \times 2K$ Vandermonde matrix of the form
%\[ \bfG^{(s)} = \begin{bmatrix} \cos(\theta_1) &  \sin(\theta_1) &  \cos(2\theta_1) &  \sin(2\theta_1) &   \dots &  \cos(K\theta_1) &  \sin(K\theta_1) \\
%\cos(\theta_2) &  \sin(\theta_2) &  \cos(2\theta_2) &  \sin(2\theta_2) &   \dots &  \cos(K\theta_2) &  \sin(K\theta_2) \\
%\vdots & & \vdots   & &   \ddots &  & \vdots \\
%\cos(\theta_{N_r}) &  \sin(\theta_{N_r}) &  \cos(2\theta_{N_r}) &  \sin(2\theta_{N_r}) &   \dots &  \cos(K\theta_{N_r}) &  \sin(K\theta_{N_r}) \\
%\end{bmatrix} ,\]
\begin{equation} \label{eq:vmon_matrix} \bfV_{K} \deff \begin{bmatrix} \zeta_1 & \overline{\zeta}_1 & \zeta_1^2 & \overline{\zeta}_1^2 &   \dots & \zeta_1^{K} & \overline{\zeta}_1^{K} \\
\zeta_2 & \overline{\zeta}_2 & \zeta_2^2 & \overline{\zeta}_2^2 &   \dots & \zeta_2^{K} & \overline{\zeta}_2^{K} 
\\
\vdots & & \vdots   & &   \ddots &  & \vdots \\
\zeta_{N_r} & \overline{\zeta}_{N_r} & \zeta_{N_r}^2 & \overline{\zeta}_{N_r}^2 &   \dots & \zeta_{N_r}^{K} & \overline{\zeta}_{N_r}^{K} 
\end{bmatrix} .\end{equation}
Matrices of this type are of interest when dealing with univariate polynomial interpolation on complex conjugate points. We introduce the map 
\[ \mc{K} : \CC \setminus \{0\} \arr \CC\setminus \{0\} \]
\[ \mc{K}(z) = \overline{z}\,^{-1} = \frac{z}{|z|^2}. \]
$\mc{K}$ is known in the literature as Kelvin transform. Then the entries of $\bfV_K$ \eqref{eq:vmon_matrix} are defined as follows: for any receptor $z_l \in \p \Omega$ we set $\zeta_l = \mc{K}(z_l) = e^{i \theta_l}/r_l$. 

To relate $\bfV_K$ with $\bfG^{(s)}$ we introduce other two matrices. By employing the product defined in Definition \ref{def:kronecker_mitra}, consider the $2K \times 2K$ diagonal scaling matrix
\[ \bfC \deff \bfI_K \kron \left\{ \ell^{-1} \bfI_2 \right\},\]
and, by setting $J = \frac{1}{2} \begin{bmatrix}
1 & -i \\ 1 & i 
\end{bmatrix}$, define the complex $2K \times 2K$ block diagonal matrix:
\[ \bfJ \deff \bfI_K \kron J. \]

It can be easily verified that $\bfG^{(s)} = -  \frac{1}{2\pi} \bfV_K \bfJ \bfC$.

Observe that there exists a permutation matrix $\bfP$ such that $\bfV_{K} \bfP = \begin{bmatrix}
\bfW_K & \overline{\bfW}_K
\end{bmatrix}$. 

\begin{lemma} \label{lem:vandermonde_minor} Let $\bfV_{(1,...,2K)}$ be the square sub-matrix  of $\bfV_{K} \bfP$ obtained by considering the first $2K$ rows. Then
	\begin{equation*}
	\mbox{det}(\bfV') = \Re(P(\zeta_1, .... , \zeta_{2K}, \overline{\zeta}_1, .... , \overline{\zeta}_{2K})),
	\end{equation*}
	where $P \in \CC[z_1, ..., z_{2K},\overline{z}_1, ..., \overline{z}_{2K}]$.	Since $\zeta_l = x_l + i y_l$ and $\overline{\zeta}_l = x_l - i y_l$, it is clear that 
	\begin{equation*}
	\mbox{det}(\bfV') = Q(x_1 , y_1 , .... , x_{2K}, y_{2K} ),
	\end{equation*}
	where $Q \in \RR [x_1 , y_1 , .... , x_{2K}, y_{2K}]$.
\end{lemma}

The set $\mc{H} = \{ Q = 0 \}$ is an hyper-surface in the affine space $\RR^{4K} \simeq \underbrace{\RR^2 \times ... \times \RR^2}_{2K}$. 
\begin{definition} \label{def:general_config}A finite set of points $\mc{S} = \{(x_1,y_1),...,(x_{2K},y_{2K})\}  \subset \RR^2$, $\#\mc{S} = 2K$, is called \emph{general configuration} if the point $(x_1,y_1,...,x_{2K},y_{2K})$ doesn't lie on $\mc{H}$.
\end{definition}

\begin{prop} Suppose there are $2K$ receptors $z_l$ such that $\zeta_l = \mathcal{K}(z_l)$ are a general configuration. Then $\bfV_{K}$ is of maximal rank.
	%	 $\rk_{\mathbb{C}}(\bfV_{K}) = 2K$.
\end{prop}
\begin{proof} Without loss of generality, let $\zeta_1 , ... , \zeta_{2K}$ are a general configuration. From Lemma \ref{lem:vandermonde_minor} it follows that $\mbox{det}(\bfV')$ is a real multivariate polynomial and it doesn't vanish when evaluated at the points $\zeta_1 , ... , \zeta_{2K}$.
\end{proof}

\begin{rem} \label{rem:abuseofdef}
	By an abuse of definition, a set of points $\{z_l\}_l$ such that $\{\mathcal{K}(z_l)\}_l$ is a general configuration, shall be called a general configuration likewise.
\end{rem}

%\begin{lemma} Given $t_1 , ... , t_K \in \{1 , ... , N_r \}$, $t_1 < ... < t_K$, let $\bfW_{(t_1 , ... , t_K) , \,: }$ be the minor of order $K$ obtained from $\bfW$ by keeping only the columns with indices $t_1 , ... , t_K$. Then
%	\begin{equation} \label{eq:weighted_det} |\mbox{det} (\bfW_{(t_1 , ... , t_K) , \,: })| = \left ( \prod_{i=1}^{K} |z_{t_i}| \right )^{-K} \left ( \prod_{1\le i< j\le K} |z_{t_i} - z_{t_j} | \right )  . \end{equation}
%\end{lemma}
%\begin{proof} Without loss of generality, we can assume $\bfW_{(t_1 , ... , t_K) , \,: } = \bfW_{(1 , ... , K) , \,: }$. Recall that
%	\begin{equation} \label{eq:det_vm} |\mbox{det}(\bfW_{(1 , ... , K) , \,: }) | = \left ( \prod_{i =1}^{K} |\zeta_{i} | \right ) \left ( \prod_{1\le i< j\le K} |\zeta_{i} - \zeta_{j} | \right ). \end{equation}
%	The Kelvin transform $\mc{K}$ is such that $|\mc{K}(z)| = |z|^{-1}$, and it has the following distortion property:
%	\[ \forall z,w \in \CC\setminus\{0\} ,\qquad | \mc{K}(z) - \mc{K}(w) | =   \frac{|z - w|}{|z||w|} . \]
%	Thus, recalling that $\zeta_i = \mc{K}(z_i)$, \eqref{eq:det_vm} becomes
%	\begin{equation}  |\mbox{det}(\bfW_{(1 , ... , K) , \,: }) | = \left ( \prod_{i =1}^{K} |z_{i} | \right )^{-1} \left ( \prod_{1\le i< j\le K} \frac{|z_{i} - z_{j} |}{|z_i||z_j|} \right ). \end{equation}
%	The result follows by noticing that \[ \prod_{1\le i< j\le K} |z_i| | z_j| =  \left ( \prod_{i=1}^K |z_i| \right )^{K-1} .\]
%\end{proof}

\subsection{Moore-Penrose inverse of $\bfS$}

For $\ell \in \mathbb{N}$, let $r_\ell(\gamma)$ be the following rotation matrix:
\[r_\ell(\gamma) := \begin{bmatrix}  \cos(\ell\gamma) & \sin (\ell\gamma) \\ -\sin(\ell\gamma) & \cos(\ell\gamma) \end{bmatrix}.\] 
%Observe that
%\[r_\ell(\gamma) = r_1(\ell\gamma), \qquad (r_1(\gamma))^\ell = r_1(\ell\gamma), \qquad (r_1(\gamma))^\top = r_1(-\gamma).\]
If we assume that the body of the fish lies in a thin annulus around the orbit of radius $\rho$ while swimming around the target, then the form of the design matrix can be simplified. In this case, if \eqref{eq:dipole_approx} is used, the information on the fish concerning its geometry and its electric field can be separated from the ``kinematics". Easy calculations show that
\[\bfG^{(s)} = \bfG^{(1)} (\bfI_K \kron \{r_1(\ell (s-1) \gamma)\}),\]
\[\bfZ_{s,\,:} = \bfZ_{1,\,:} \,(\bfI_K \kron \{r_1((1+\ell)(s-1)\gamma)\}),\]
where
\[ \bfZ_{1,:} = \begin{bmatrix}
\cos(2\overline{\theta}_1) & \sin(2\overline{\theta}_1) & \cos(3\overline{\theta}_1) & \sin(3\overline{\theta}_1) & \dots & \cos((K+1)\overline{\theta}_1) & \sin((K+1)\overline{\theta}_1)
\end{bmatrix}  . \]
Moreover, given $\alpha_s=\alpha_0 + (s-1) \gamma$ angle of $\textbf{p}_s$, we have 
\[\bfP^{(s)}_K = \bfI_K \kron \left\{(-1)^{\ell+1}r_1(\alpha_0 + (s-1) \gamma)\right\}.\]
Hereinafter, we assume that the fish moves with its electric organ along on a circular orbit of radius $\rho$. 

Recall that
\[\bfS = \bfS_{dip} + \bfS_{SL}.\]
By using the rotational symmetry of the configuration and the notation hereabove, we rewrite $\bfS_{dip}$ and $\bfS_{SL}$ as follows:
\[
\begin{aligned} 
(\bfS_{dip})_{s,\,:} & = \bfZ_{s, \,:} \,\bfP^{(s) \top}_K \bfD_{2,K+1} \\ & = \bfZ_{1,\,:}\, (\bfI_K \kron \{r_1((1+\ell)(s-1)\gamma)\}_\ell) (\bfI_K \kron \{(-1)^{\ell+1}r_1(-\alpha_0 - (s-1) \gamma)\}_\ell) \bfD_{2,K+1} \\ & 
= \bfZ_{1,\,:} \,(\bfI_K \kron \{(-1)^{\ell+1}r_1(-\alpha_0 + \ell (s-1) \gamma)\}_\ell) \bfD_{2,K+1},
\end{aligned}
\]
and
\[
\begin{aligned} 
(\bfS_{SL})_{s,\,:} & = -\textbf{u} \bfG^{(1)} (\bfI_K \kron \{r_1(\ell (s-1) \gamma)\}_{\ell}).
\end{aligned}
\]
Here, $\bfD_{2,K+1}$ is given as in Section 3.2 by 
\[
\bfD_{2,K+1} =\bfI_K \kron \{\rho^{-(\ell+1)}\bfI_2\}_\ell.
\]
Finally,
\[
\begin{aligned} 
%\bfS_{dip} & = (\bfI_M \kron \bfZ_{0,\,:})\begin{bmatrix}(\bfI_K \kron \{(-i)^{i+1}r_1(-\alpha_0)\}_i) \bfD_{2,K+1} \\ \vdots \\ (\bfI_K \kron \{(-i)^{i+1}r_1(-\alpha_0 + i(M-1) \gamma)\}_i) \bfD_{2,K+1}
%\end{bmatrix} \\
%& = 
\bfS_{dip} = (\bfI_M \kron \bfZ_{1,\,:})\begin{bmatrix}\bfI_K \kron \{(-1)^{\ell+1}r_1(-\alpha_0)\}_\ell \\ \vdots \\ \bfI_K \kron \{(-1)^{\ell+1}r_1(-\alpha_0 + \ell(M-1) \gamma)\}_\ell \end{bmatrix}
\bfD_{2,K+1},
\end{aligned}
\] 
%\[\bfS_{dip} = (\bfI_M \kron \bfZ_{0,\,:})\begin{bmatrix}\bfZ_{0,\,:} (\bfI_K \kron \{(-i)^{i+1}r_1(-\alpha_0)\}_i) \bfD_{2,K+1} \\ \vdots \\ \bfZ_{M-1,\,:} (\bfI_K \kron \{(-i)^{i+1}r_1(-\alpha_0 + i(M-1) \gamma)\}_i) \bfD_{2,K+1}
%\end{bmatrix},\] 
and 
\[\bfS_{SL} = (\bfI_M \kron -\textbf{u} \bfG^{(1)}) \begin{bmatrix} \bfI_K \kron r_1(0) \\ \vdots \\ \bfI_K \kron \{r_1(\ell(M-1) \gamma)\}_\ell \end{bmatrix}. \]

Hereinafter, we assume that $\partial^{\alpha} \mathcal{S}_{\Omega_s}(z)$ is non-zero for each $|\alpha|\leq K$. In the far field regime this assumption guarantees that the GPTs can be retrieved up to order $K$. It will become clear soon that this condition also makes $\bfS_{SL}$ full column rank, allowing to compute its Moore-Penrose inverse as
\[\bfS_{SL}^\dagger = (\bfS_{SL}^\top \bfS_{SL})^{-1} \bfS_{SL}^\top.\]
We have the following lemma. 

\begin{lemma}
	For $1 \leq k \leq 2K$ and $M \gg 1$, we have
	\begin{equation} 
	\begin{aligned} 
	\|(\bfS_{SL}^{\emph{\dagger}})_{k,\,:}\|_F \lesssim \frac{\left\lceil k/2 \right\rceil \rho^{\left\lceil k/2 \right\rceil}}{\sqrt{M}}. 
	\label{estim_1} 
	\end{aligned}
	\end{equation} 
\end{lemma}
\begin{proof} 
	Let us denote the $1 \times 2K$ row vector $-\textbf{u} \bfG^{(1)}$ by $\bf{w}$. We need to compute $\bfS_{SL}^\top \bfS_{SL}$, that is 
	\[
	\begin{bmatrix} \bfI_K \kron r_1(0) & \ldots & \bfI_K \kron \{r_1(-\ell (M-1) \gamma)\}_\ell \end{bmatrix} (\bfI_M \kron \textbf{w}^\top) (\bfI_M \kron \textbf{w})\begin{bmatrix} \bfI_K \kron r_1(0) \\ \vdots \\ \bfI_K \kron \{r_1(\ell (M-1) \gamma)\}_\ell \end{bmatrix}.
	\]
	Since $(\bfI_M \kron \textbf{w}^\top) (\bfI_M \kron \textbf{w}) = \bfI_M \kron \textbf{w}^\top \textbf{w}$, the product $\bfS_{SL}^\top \bfS_{SL}$ boils down to
	\[
	\bfS_{SL}^\top \bfS_{SL} = \sum_{s=1}^{M} (\bfI_K \kron \{r_1(-\ell (s-1)\gamma)\}_\ell) \textbf{w}^\top \textbf{w} (\bfI_K \kron \{r_1(\ell (s-1)\gamma)\}_\ell).
	\]
	Notice that the receptor $x_{l}^{(1)}$ is $x_{l}^{(1)} = r_l e^{i \theta_l^{(1)}}$. Since $\rho - \eta \leq r_i \leq \rho + \eta $, for some $\eta$ small, we assume that $x_{l}^{(1)} = \rho e^{i \theta_l^{(1)}}$. As a consequence, we can factorize $\bfG^{(1)}$ as
	\[
	\bfG^{(1)} = \widetilde{\bfG}^{(1)} (\bfI_K \kron \{\ell^{-1}\rho^{-\ell} \bfI_2 \}_\ell),
	\]
	where $\widetilde{\bfG}^{(1)}$ is a Vandermonde-type matrix with nodes on the unit disk $|z| = 1$.
	
	We obtain
	\[
	\textbf{w}^\top \textbf{w} = (\bfI_K \kron \{ \ell^{-1}\rho^{-\ell} \bfI_2 \}_\ell)  (\widetilde{\bfG}^{(1)})^\top \textbf{u}^\top \textbf{u} \widetilde{\bfG}^{(1)} (\bfI_K \kron \{ \ell^{-1}\rho^{-\ell} \bfI_2 \}_\ell). 
	\]
	It is immediate to see that
	\[
	\begin{aligned} 
	\bfS_{SL}^\top \bfS_{SL} = \bfD_{1,K} \bfC \left(\sum_{s=1}^{M} (\bfI_K \kron \{r_1(-\ell (s-1)\gamma)\}_\ell) \widetilde{\textbf{w}}^\top \widetilde{\textbf{w}}(\bfI_K \kron \{r_1(\ell (s-1)\gamma)\}_\ell)\right) \bfC \bfD_{1,K},
	\end{aligned}
	\]
	where
	\[
	\bfD_{1,K} := \bfI_K \kron \{ \rho^{-\ell} \bfI_2 \}_\ell,  \qquad \widetilde{\textbf{w}} = - \textbf{u} \widetilde{\bfG}^{(1)}.\]
	%Notice that $\widetilde{\textbf{w}}^\top\widetilde{\textbf{w}}$ is a rank-one matrix. 
	
	Let us denote $\bfS_{SL}\bfD_{1,K}^{-1}\bfC^{-1}$ by $\widetilde{\bfS}_{SL}$. 
	
	To prove estimate \eqref{estim_1}, we rely on the following inequality:
	\begin{equation}
	\begin{aligned} 
	\|(\bfS_{SL}^\dagger)_{k,\,:}\|_F & \leq  \|(\bfC^{-1}\bfD_{1,K}^{-1})_{k,\,:}\|_F \; \|(\widetilde{\bfS}_{SL}^\top \widetilde{\bfS}_{SL})^{-1}\|_F \; \|\widetilde{\bfS}_{SL}^\top\|_F. 
	\end{aligned}
	\label{lem:ineq}
	\end{equation}
	
	By a straightforward calculation, it is immediate to see that
	\begin{equation}
	\|(\bfC^{-1}\bfD_{1,K}^{-1})_{k,\,:}\|_F \leq \left\lceil k/2 \right\rceil \rho^{\left\lceil k/2 \right\rceil}, 
	\label{lem:estim_1} 
	\end{equation}
	and
	\begin{equation}
	\| \widetilde{\bfS}_{SL}^\top\|_F \leq \sqrt{M} \max_{1 \leq m \leq M} \| (\widetilde{\bfS}_{SL})_{m,\,:}\|_F  \lesssim \sqrt{M}.
	\label{lem:estim_2} 
	\end{equation}
	
	Finally, we investigate the Frobenius norm of $(\widetilde{\bfS}_{SL}^\top \widetilde{\bfS}_{SL})^{-1}$. For $1 \leq j \leq K$, we observe that
	\[
	\begin{aligned} 
	\left(\widetilde{\bfS}_{SL}^\top \widetilde{\bfS}_{SL}\right)_{2j-1:2j, \, 2j-1:2j} & = \sum_{s=1}^{M} r_1(-j (s-1) \gamma ) 
	\begin{bmatrix}
	\widetilde{w}_{2j-1} \\ \widetilde{w}_{2j} 
	\end{bmatrix}
	\begin{bmatrix}
	\widetilde{w}_{2j-1} & \widetilde{w}_{2j} 
	\end{bmatrix}
	r_1(j (s-1) \gamma ) \\ & = U^\mathsf{H}  \left( \sum_{s=1}^{M} \begin{bmatrix}
	e^{-ij(s-1) \gamma} & 0 \\ 0 & e^{ij(s-1) \gamma} 
	\end{bmatrix} \begin{bmatrix}
	{y}_{2j-1} \\ {y}_{2j} 
	\end{bmatrix}
	\begin{bmatrix}
	{y}_{2j-1} \\ {y}_{2j} 
	\end{bmatrix}^\mathsf{H}
	\begin{bmatrix}
	e^{ij(s-1) \gamma} & 0 \\ 0 & e^{-ij(s-1) \gamma} 
	\end{bmatrix} \right) U \\
	%& = V^\mathsf{H}  \left( \sum_{\ell=0}^{M-1} \begin{bmatrix}
	%e^{-ij\ell \gamma} & 0 \\ 0 & e^{ij\ell \gamma} 
	%\end{bmatrix} \begin{bmatrix}
	%|{y}_{2j-1}|^2 & {y}_{2j} {y}_{2j-1}^\mathsf{H} \\  {y}_{2j-1} {y}_{2j}^\mathsf{H} & |{y}_{2j}|^2
	%\end{bmatrix}
	%\begin{bmatrix}
	%e^{ij\ell \gamma} & 0 \\ 0 & e^{-ij\ell \gamma} 
	%\end{bmatrix}\right) V \\
	%& = 
	%V^\mathsf{H}  
	%\begin{bmatrix}
	%M |{y}_{2j-1}|^2 & ({y}_{2j} {y}_{2j-1}^\mathsf{H}) \sum_{\ell=0}^{M-1} e^{-2ij\ell \gamma} \\  ({y}_{2j-1} {y}_{2j}^\mathsf{H}) \sum_{\ell=0}^{M-1} e^{2ij\ell \gamma} & M |{y}_{2j}|^2
	%\end{bmatrix}
	%V\\
	& = 
	U^\mathsf{H}  
	\begin{bmatrix}
	M |{y}_{2j-1}|^2 & \frac{1-e^{-2ijM\gamma}}{1-e^{-2ij\gamma}} ({y}_{2j} \overline{y}_{2j-1}) \\  \frac{1-e^{2ijM\gamma}}{1-e^{2ij\gamma}}({y}_{2j-1} \overline{y}_{2j}) & M |{y}_{2j}|^2
	\end{bmatrix}
	U,
	\end{aligned} 
	\]
	where 
	\[
	\begin{bmatrix}
	{y}_{2j-1} \\ {y}_{2j} 
	\end{bmatrix} = U \begin{bmatrix}
	\widetilde{w}_{2j-1} \\ \widetilde{w}_{2j}
	\end{bmatrix}, \qquad
	U = 
	\frac{1}{\sqrt{2}}\begin{bmatrix}
	1 & i\\1 & -i 
	\end{bmatrix}.
	\]
	
	For $1 \leq j < k \leq K$, we have 
	\[
	\begin{aligned} 
	\hspace{-14mm} \left(\widetilde{\bfS}_{SL}^\top \widetilde{\bfS}_{SL}\right)_{2j-1:2j, \, 2k-1:2k} & = \sum_{s=1}^{M} r_1(-j (s-1) \gamma )
	\begin{bmatrix}
	\widetilde{w}_{2j-1} \\ \widetilde{w}_{2j} 
	\end{bmatrix}
	\begin{bmatrix}
	\widetilde{w}_{2k-1} & \widetilde{w}_{2k} 
	\end{bmatrix}
	r_1(k (s-1) \gamma ) \\ 
	%& = V^\mathsf{H}  \left( \sum_{\ell=0}^{M-1} \begin{bmatrix}
	%e^{-ij\ell \gamma} & 0 \\ 0 & e^{ij\ell \gamma} 
	%\end{bmatrix} \begin{bmatrix}
	%{y}_{2j-1} \\ {y}_{2j} 
	%\end{bmatrix}
	%\begin{bmatrix}
	%{y}_{2k-1} \\ {y}_{2k} 
	%\end{bmatrix}^\mathsf{H}
	%\begin{bmatrix}
	%e^{ik\ell \gamma} & 0 \\ 0 & e^{-ik\ell \gamma} 
	%\end{bmatrix} \right) V \\
	%& = V^\mathsf{H}  \left( \sum_{\ell=0}^{M-1} \begin{bmatrix}
	%e^{-ij\ell \gamma} & 0 \\ 0 & e^{ij\ell \gamma} 
	%\end{bmatrix} \begin{bmatrix}
	%{y}_{2j-1} y_{2k-1}^\mathsf{H} & {y}_{2j} {y}_{2k-1}^\mathsf{H} \\  {y}_{2j-1} {y}_{2k}^\mathsf{H} & {y}_{2j} y_{2k}^\mathsf{H}
	%\end{bmatrix}
	%\begin{bmatrix}
	%e^{ik\ell \gamma} & 0 \\ 0 & e^{-ik\ell \gamma} 
	%\end{bmatrix}\right) V \\
	%& = 
	%V^\mathsf{H}  
	%\begin{bmatrix}
	%({y}_{2j-1} y_{2k-1}^\mathsf{H}) \sum_{\ell=0}^{M-1} e^{-i(j-k)\ell \gamma} & ({y}_{2j} {y}_{2k-1}^\mathsf{H}) \sum_{\ell=0}^{M-1} e^{-i(j+k)\ell \gamma} \\  ({y}_{2j-1} {y}_{2k}^\mathsf{H}) \sum_{\ell=0}^{M-1} e^{i(j+k)\ell \gamma} & M y_{2j} y_{2k}^{\mathsf{H}} \sum_{\ell=0}^{M-1} e^{i(j-k)\ell \gamma} 
	%\end{bmatrix}
	%V\\
	& = 
	U^\mathsf{H}  
	\begin{bmatrix}
	\frac{1-e^{-i(j-k)M\gamma}}{1-e^{-i(j-k)\gamma}} ({y}_{2j-1} \overline{y}_{2k-1}) & \frac{1-e^{-i(j+k)M\gamma}}{1-e^{-i(j+k)\gamma}} ({y}_{2j} \overline{y}_{2k-1}) \\  \frac{1-e^{i(j+k)M\gamma}}{1-e^{i(j+k)\gamma}}({y}_{2j-1} \overline{y}_{2k}) & \frac{1-e^{i(j-k)M\gamma}}{1-e^{i(j-k)\gamma}} ({y}_{2j} \overline{y}_{2k})
	\end{bmatrix}
	U.
	\end{aligned} 
	\]
	Therefore, $\widetilde{\bfS}_{SL}^\top \widetilde{\bfS}_{SL}$ can be written as follows:
	\[
	\widetilde{\bfS}_{SL}^\top \widetilde{\bfS}_{SL} = \bfU^\mathsf{H} (M \bfD + \bfR) \bfU,
	\]
	where
	\[
	\bfD := \begin{bmatrix}
	|{y}_{1}|^2 & & 0 \\ 
	& \ddots & \\  
	0 &  & |{y}_{2K}|^2
	\end{bmatrix}
	, \qquad \bfR:= \widetilde{\bfS}_{SL}^\top \widetilde{\bfS}_{SL} - \bfD, \qquad 
	\bfU := \bfI_K \kron U.
	\]
	Notice that requiring $\bfD$ to be invertible is equivalent to saying that each $|y_j|$, and thus each $\widetilde{w}_j$, is non-zero. For this reason, the initial assumption on $\partial^{\alpha} \mathcal{S}_{\Omega_s}(z)$  implies that $\bfD$ is invertible (see \cite{Ammari11652}) and thus, for $M$ large enough, $\bfS_{SL}$ is full column rank. Since $\| \bfR \|_F = O(1)$, we have \cite{ONERRORBOUND}
	\[
	\bfU (\widetilde{\bfS}_{SL}^\top \widetilde{\bfS}_{SL})^{-1} \bfU^\mathsf{H}=  \frac{1}{M} \left(\bfD + \frac{1}{M}\bfR\right)^{-1} = \frac{1}{M}\bfD^{-1} + O\left(\frac{1}{M^2}\right). 
	\]
	Since $\bfU$ is unitary, $\| \bfU (\widetilde{\bfS}_{SL}^\top \widetilde{\bfS}_{SL})^{-1} \bfU^\mathsf{H}\|_F = \|(\widetilde{\bfS}_{SL}^\top \widetilde{\bfS}_{SL})^{-1}\|_F$. Hence
	\begin{equation}
	\|(\widetilde{\bfS}_{SL}^\top \widetilde{\bfS}_{SL})^{-1}\|_F \lesssim \frac{1}{M}. 
	\label{lem:estim_3}
	\end{equation}
	%TO BE CHANGED
	%\[
	%\sum_{k=1}^{K} \|\left(\widetilde{\bfS}_{SL}^\top \widetilde{\bfS}_{SL}\right)_{2k-1:2k,2k-1:2k}\|^2_F \leq \|\widetilde{\bfS}_{SL}^\top \widetilde{\bfS}_{SL}\|^2_F.
	%\]
	%Since $V$ is unitary, $\| V^\mathsf{H} AV\|^2_F = \| A\|^2_F$ for each $A \in M_{2,2}(\mathbb{C})$. Hence
	%\[
	%0<M\sqrt{\sum_{k=1}^{K} (|y_{2k-1}|^4 + |y_{2k}|^4)} \leq \|\widetilde{\bfS}_{SL}^\top \widetilde{\bfS}_{SL}\|_F.
	%\]
	%%We have
	%%\[
	%%\begin{aligned} 
	%%\|\bfS_{SL}^\dagger\|_F & = \| (\bfC \bfD_{1,K} \bfC^{-1} \bfD_{1,K}^{-1}\bfS_{SL}^\top \bfS_{SL}\bfD_{1,K}^{-1} \bfC^{-1} \bfC \bfD_{1,K})^{-1} \bfC \bfD_{1,k}\bfC^{-1}\bfD_{1,K}^{-1}\bfS_{SL}^\top \|_F \\ & = \| \bfC^{-1}\bfD_{1,K}^{-1} (\bfC^{-1}\bfD_{1,K}^{-1}\bfS_{SL}^\top \bfS_{SL}\bfD_{1,K}^{-1}\bfC^{-1})^{-1} \bfD_{1,K}^{-1}\bfC^{-1}\bfC \bfD_{1,K}\bfC^{-1}\bfD_{1,K}^{-1}\bfS_{SL}^\top \|_F \\ & = \| \bfC^{-1} \bfD_{1,K}^{-1} (\bfC^{-1}\bfD_{1,K}^{-1}\bfS_{SL}^\top \bfS_{SL}\bfD_{1,K}^{-1}\bfC^{-1})^{-1} (\bfC^{-1}\bfD_{1,K}^{-1}\bfS_{SL}^\top) \|_F \\ & \leq \frac{ \|\bfC^{-1}\bfD_{1,K}^{-1}\|_F \|\bfC^{-1}\bfD_{1,K}^{-1}\bfS_{SL}^\top\|_F}{\|\bfC^{-1}\bfD_{1,K}^{-1}\bfS_{SL}^\top \bfS_{SL}\bfD_{1,K}^{-1}\bfC^{-1}\|_F} \\& \lesssim \frac{\sqrt{\sum_{\ell=1}^K \ell^2 \rho^{2\ell}}}{M}. 
	%%\end{aligned}
	%%\] 
	Finally, inequality \eqref{lem:ineq} together with estimates \eqref{lem:estim_1}, \eqref{lem:estim_2} and \eqref{lem:estim_3} shows that $\|(\bfS_{SL}^{{\dagger}})_{k,\,:}\|_F$ satisfies \eqref{estim_1}. 
\end{proof}

We are now ready to estimate the dependency of $(\bfS^\dagger)_{k, \, :}$ on $M$ and $\rho$. For $\rho$ large enough, we observe that the leading order term of $(\bfS^\dagger)_{k, \, :}$ is  $(\bfS_{SL}^\dagger)_{k, \, :}$.

Since
\begin{equation*}
\bfS = \bfS_{dip} + \bfS_{SL} ,
\end{equation*}
we have \cite{ONERRORBOUND}
\begin{equation*}
\bfS^\dagger = (\bfS_{dip} + \bfS_{SL})^\dagger = \bfS_{SL}^\dagger - \bfS_{SL}^\dagger \bfS_{dip} \bfS_{SL}^\dagger + \sum_{\ell=2}^\infty (-1)^\ell (\bfS_{SL}^\dagger \bfS_{dip})^\ell \bfS_{SL}^\dagger.
\end{equation*}
%Observe that
%\[
%(\bfS^\dagger)_{k, \, :} = (\bfS_{SL}^\dagger)_{k, \, :} - (\bfS_{SL}^\dagger)_{k, \, :} \bfS_{dip} \bfS_{SL}^\dagger + \ldots
%\]
%We denote by $\bfA$ the matrix $(\bfC^{-1}\bfD_{1,K}^{-1}\bfS_{SL}^\top \bfS_{SL}\bfD_{1,K}^{-1}\bfC^{-1})^{-1} (\bfC^{-1}\bfD_{1,K}^{-1}\bfS_{SL}^\top)$ and by $\bfB$ the matrix $\bfS_{dip} \bfD_{2,K+1}^{-1}$. Notice that both $\bfA$ and $\bfB$ are independent of $\rho$. We have
%\[
%\begin{aligned} 
%\|(\bfS_{SL}^\dagger)_{k, \, :} \bfS_{dip} \bfS_{SL}^\dagger\|_F & = \|(\bfC^{-1}\bfD_{1,K}^{-1})_{k,\,:} \bfA \bfB \bfD_{2,K+1} \bfC^{-1} \bfD_{1,K}^{-1} \bfA\|_F \\
%& = \rho^{-1}  \|(\bfC^{-1}\bfD_{1,K}^{-1})_{k,\,:} \bfA \bfB \bfC^{-1}\bfA\|_F \\
%& \lesssim \rho^{\left\lceil k/2 \right\rceil-1}.
%\end{aligned} 
%\]
%Finally, we obtain
It is immediate to prove the following lemma.
\begin{lemma} \label{lem:estimate_mp_S} For $1 \leq k \leq 2K$, we have
	\[
	(\bfS^{\emph{\dagger}})_{k, \, :} = (\bfS_{SL}^{\emph{\dagger}})_{k, \, :} + O(\rho^{\left\lceil k/2 \right\rceil-1}).
	\]
\end{lemma} 

\begin{proof}
	Let us denote the matrix $\bfS_{dip} \bfD_{2,K+1}^{-1}$ by $\widetilde{\bfS}_{dip}$. Observe that both $\widetilde{\bfS}_{SL}$ and $\widetilde{\bfS}_{dip}$ are independent of $\rho$. We have
	\[
	\begin{aligned}
	\|(\bfS_{SL}^\dagger)_{k, \, :} \bfS_{dip} \bfS_{SL}^\dagger\|_F & = \|(\bfC^{-1}\bfD_{1,K}^{-1})_{k,\,:} \widetilde{\bfS}_{SL}^\dagger \widetilde{\bfS}_{dip} \bfD_{2,K+1} \bfC^{-1} \bfD_{1,K}^{-1} \widetilde{\bfS}_{SL}^\dagger\|_F \\
	& = \rho^{-1}  \|(\bfC^{-1}\bfD_{1,K}^{-1})_{k,\,:} \widetilde{\bfS}_{SL}^\dagger \widetilde{\bfS}_{dip} \bfC^{-1}\widetilde{\bfS}_{SL}^\dagger\|_F \\
	& \lesssim \rho^{\left\lceil k/2 \right\rceil-1}.
	\end{aligned}
	\]
\end{proof}

\section{Transferable Belief Model}
In this appendix we present some basic definitions from evidence theory. In particular, we shall consider the Transferable Belief Model (TBM), see \cite{SMETS1992292}. This theory doesn't require any underlying probability space.

%Let $\mc{L}$ be a finite propositional language supplemented by the tautology $\top$ and the contradiction $\perp$, and closed under the usual Boolean connectives $\neg$, $\vee$ and $\wedge$. The propositions of $\mc{L}$ aim at describing the uncertainty of the  problem. As usual, we shall identify the propositions of our language with  subsets of a set $\Omega = \{ \omega_1, ... , \omega_N \}$, which is called \textit{frame of discernment}. We assume that there exists $\omega_0 \in \Omega$ which corresponds to the \textit{actual world}.
%In other words, the elements of $\Omega$ corresponds to the interpretations of $\mc{L}$ \cite{SMETS1992292}.

%Given a partition $\Pi$ of $\Omega$, let $\mathfrak{R}$ be the Boolean algebra of sets generated by $\Pi$.

%The pair $(\Omega , \mathfrak{R})$ is called \textit{propositional space}.

In the context of dictionary classification, the frame of discernment is usually modeled as a finite set $\mathscr{C} = \{ c_1 , c_2 , ... , c_N \}$, which is called \textit{dictionary}.

A belief function is a function $\mathsf{bel} : 2^{\mathscr{C}} \arr [0,1]$ such that:

\begin{enumerate}
	\item $\mathsf{bel}(\emptyset) = 0$;
	\item for all $A_1 , A_2 , ... , A_n \in 2^{\mathscr{C}}$,
	\begin{equation}
		\bel(A_1 \cup A_2 \cup ... \cup A_n) \ge \sum_i \bel (A_i) - \sum_{i > j} \bel (A_i \cap A_j) - \;\;\dots\;\; - (-1)^n \bel(A_1 \cap A_2 \cap ... \cap A_n) \;;
	\end{equation} 
	\item $\bel (\mathscr{C}) \le 1$.
\end{enumerate}
%The triple $(\mathscr{C} , 2^{\mathscr{C}} , \mathsf{bel})$ is called a \textit{credibility space}. 

%
A basic belief assignment (BBA) $\mathsf{m}$ is a function $\mathsf{m} : 2^{\mathscr{C}} \arr [0,1]$ such that
\[ \sum_{A\in 2^{\mathscr{C}}} \mathsf{m} (A) = 1 . \]
The value $\mathsf{m}(A)$ for $A \in 2^{\mathscr{C}}$ is called the basic belief mass (bbm) given to $A$. This is a part of the agent's belief that supports $A$, and that, due to lack of information, does not support any strict subset of $A$.
If $\mathsf{m}(A) > 0$ then $A$ is called focal set. Observe that it is allowed to allocate positive BBM to $\emptyset$, i.e., $\mathsf{m}(\emptyset) > 0$.

The basic belief assignment (BBA) related to a belief function $\bel$ is the function $\mathsf{m} : 2^{\mathscr{C}} \arr [0,1]$ such that:
\begin{equation} \label{eq:bba-def}
	\begin{split}
		\mathsf{m}(A) &= \sum_{B \in 2^{\mathscr{C}} ,\; \emptyset \ne B \subseteq A} (-1)^{|A|-|B|} \bel(B) , \quad \mbox{for all } A \in 2^{\mathscr{C}} ,\; A \ne \emptyset,\\
		\mathsf{m}(\emptyset) &= 1 - \bel (\mathscr{C}).
	\end{split}
\end{equation}

Moreover, there is a one-to-one correspondence between $\mathsf{m}$ and $\bel$ via the following formula 
\begin{equation}
	\bel(A) = \sum_{B \in \mathscr{C} ,\; \emptyset \ne B \subseteq A}  \mathsf{m}(B) , \quad \mbox{for all } A \in 2^{\mathscr{C}} ,\; A \ne \emptyset,
\end{equation}
that is, $\bel(A)$ is obtained by summing all BBMs given to subsets $B \in 2^{\mathscr{C}} $ with $B \subseteq A$, and it quantifies the total amount of justified specific support given to $A$.

Let $\mathsf{m}_1$ and $\mathsf{m}_2$ be two BBAs, we define the TBM conjunctive combination of the two as
\begin{equation} \label{eq:TBM-conjunctive-rule}
(\mathsf{m}_{1} \myointersection \mathsf{m}_{2} ) (A) = \sum_{B \cap C = A} 
\mathsf{m}_1(B)\mathsf{m}_2(C).
\end{equation}

%Now we aim at introducing another representation for a nondogmatic BBA  $\mathsf{m}$.

%A BBA is called \textit{simple} (SBBA) if it has at most two focal sets and, if it has two, $\mathscr{C}$ is one of those. A SBBA $\mathsf{m}$ such that $\mathsf{m}(A) = 1 - w$ for some $A \ne \mathscr{C}$ and $ \mathsf{m}(\mathscr{C}) = w$ can be denoted by $A^w$. Observe that the \textit{vacuous} BBA can be denoted by $A^1$.
%
%Given $A^{w_1}$ and $A^{w_2}$ two SBBAs with the same focal element $A \ne \mathscr{C}$, the TBM conjunctive rule reads \begin{equation} \label{eq:TBM_Aw1Aw2} A^{w_1} \myointersection A^{w_2} = A^{w_1 w_2} .\end{equation}

In the TBM model, from a decision-making point of view we need to resort to a probability distribution in order to select the most reliable hypothesis in $\mathscr{C}$. Such function is called \textit{pignistic probability}, and is defined as
\begin{equation} \label{eq:pignisticprob}
\mbox{BetP}(c) = \sum_{A \subseteq \mathscr{C}, c \in A} \frac{\mathsf{m}(A)}{|A|} .
\end{equation}
The term \emph{pignistic} stresses the fact that the only purpose of these probabilities is that of forcing a decision.

\FloatBarrier

 \bibliographystyle{ieeetr}
 \bibliography{biblio.bib} 

\end{document}